\documentclass[aos]{imsart}

\RequirePackage{amsthm,amsmath,amsfonts,amssymb}
\RequirePackage[numbers]{natbib}
\RequirePackage[colorlinks,citecolor=blue,urlcolor=blue]{hyperref}
\RequirePackage{graphicx}

\startlocaldefs
\theoremstyle{plain}

\newtheorem{theorem}{Theorem}[section]

\theoremstyle{remark}

\usepackage{smile, bbm,relsize}
\newcommand{\T}{{\mathsmaller {\rm T}}}
\newcommand*{\supp}{\mathrm{supp}}
\newcommand*{\var}{\textnormal{var}}
\newcommand{\inte}{{{\rm int}}}
\newcommand{\PPi}{{\mathlarger \Pi}}
\newcommand{\nn}{\nonumber}

\def\##1\#{\begin{align}#1\end{align}}
\def\$#1\${\begin{align*}#1\end{align*}}
 
\def\sn{\sum_{i=1}^n}

\newcommand{\BB}{\mathbb{B}}
\newcommand{\1}{\mathbbm{1}}

\newcommand{\wt}{\widetilde}
\newcommand{\bfsym}[1]{\ensuremath{\boldsymbol{#1}}}
       \def \bbeta    {\bfsym{\beta}}
\def \bgamma   {\bfsym{\gamma}}

\endlocaldefs

\begin{document}

\begin{frontmatter}
\title{Scalable Estimation and Inference for \\  Censored Quantile Regression Process}
\runtitle{Censored quantile regression in high dimensions}

\begin{aug}
\author[A]{\fnms{Xuming} \snm{He}\ead[label=e1]{xmhe@umich.edu}},
\author[B]{\fnms{Xiaoou} \snm{Pan}\ead[label=e2,mark]{xip024@ucsd.edu}}
\author[A]{\fnms{Kean Ming} \snm{Tan}\ead[label=e3,mark]{keanming@umich.edu}}

\and
\author[B]{\fnms{Wen-Xin} \snm{Zhou}\ead[label=e4,mark]{wez243@ucsd.edu}}
\address[A]{Department of Statistics,
University of Michigan,
\printead{e1,e3}}

\address[B]{Department of Mathematics,
University of California, San Diego,
\printead{e2,e4}}
\end{aug}

\begin{abstract}
Censored quantile regression (CQR) has become a valuable tool to study the heterogeneous association between a possibly censored outcome and a set of covariates, yet computation and statistical inference for CQR have remained a challenge for large-scale data with many covariates. In this paper, we focus on a smoothed martingale-based sequential estimating equations approach, to which scalable gradient-based algorithms can be applied. Theoretically, we provide a unified analysis of the smoothed sequential estimator and its penalized counterpart in increasing dimensions. When the covariate dimension grows with the sample size at a sublinear rate, we establish the uniform convergence rate (over a range of quantile indexes) and provide a rigorous justification for the validity of a multiplier bootstrap procedure for inference. In high-dimensional sparse settings,  our results considerably improve the existing work on CQR by relaxing an exponential term of sparsity.  We also demonstrate the advantage of the smoothed CQR over existing methods with both simulated experiments and  data applications.
\end{abstract}



\end{frontmatter}

\section{Introduction}
\label{sec:intro}

Censored data are prevalent in many applications where the response variable of interest is partially observed, mostly due to loss of follow-up. For instance, in a lung cancer study considered by \citet{Setal2008}, 46.6\% of the lung cancer patients' survival time are censored, due to either early withdrawal from the study or death because of other reasons that are unrelated to lung cancer.  Commonly used methods to study the association between the censored response and explanatory variables (covariates) are through the use of Cox proportional hazards model and the accelerated failure time (AFT) model \citep{ABGK1993, KK2012}. Both models assume homogeneous covariate effects and are not applicable to cases in which the lower and upper quantiles of the conditional distribution of the censored response, potentially with different covariate effects, are of interest. Moreover, in many scientific studies, higher or lower quantiles of the response variable are more of interest than the mean. To capture heterogeneous covariate effects and to better predict the response at different quantile levels,  various censored quantile regression (CQR) methods have been developed under different assumptions on the censoring mechanism \citep{P1984,P1986, YJW1995,BH1998,CH2002,HKP2002,P2003,WW2009,LT2013,YNH2018, DEV2019, DEV2020}.
We refer the reader to Chapters~6 and 7 in \cite{KCHP2017} as well as \cite{P2021} for a comprehensive review of censored quantile regression.

We consider the random right censoring mechanism,  in which the censoring points are unknown for the uncensored observations.  Statsitcal methods for CQR were first proposed under the stringent assumption that the uncensored response variable (not observable due to censoring) is marginally independent of the censoring variable; see, for example \cite{YJW1995, HKP2002}. 
Under a more relaxed conditional independence assumption, conditioned on the covariates,  \cite{P2003} generalized the Kaplan-Meier estimator for estimating the (univariate) survival function to the regression setting, based on  \citet{E1967}'s redistribution-of-mass construction. From a different perspective, \cite{PH2008} employed a martingale-based approach for fitting CQR, and the resulting method has been shown to be closely related to \cite{P2003}'s method \citep{NBP2006,P2012}. Both \cite{P2003}'s and \cite{PH2008}'s methods,  along with their variants,  involve solving a series of quantile regression problems that can be reformulated as linear programs, solvable by the simplex or interior point method \citep{BR1974, PK1997, KM2014}. Statistical properties of the aforementioned methods have been well studied, assuming that the number of covariates, $p$, is fixed \citep{NBP2006, PH2008, PL2010, P2012}. To this date, the impact of dimensionality in the increasing-$p$ regime, in which $p$ is allowed to increase with the number of observations, has remained unclear in the presence of censored outcomes.

In the high-dimensional setting in which $p>n$, convex and nonconvex penalty functions are often employed to perform variable selection and to achieve a trade-off between statistical bias and model complexity. 
While penalized Cox proportional hazards and  AFT models have been well studied \citep{FL2002, HMX2006, CHT2009, BFJ2011}, existing work on penalized CQR under the framework of \cite{P2003} and \cite{PH2008} in the high-dimensional setting is still lagging. Large-sample properties of penalized CQR estimators were first derived under the fixed-$p$ setting ($p<n$),  mainly due to the technical challenges introduced by the sequential nature of the procedure \citep{SLZ2010,WZL2013,VWD2014}. More recently,  \cite{ZPH2018} studied a penalized CQR estimator, extending the method of \cite{PH2008} to the high-dimensional setting $(p>n)$.  
 They showed that the estimation error (under $\ell_2$-norm) of the $\ell_1$-penalized CQR estimator is upper bounded by $\cO\big(\exp(C s)   \sqrt{s \log (p) / n}\big)$ with high probability, where $C>0$ is a dimension-free constant. Compared to the $\ell_1$-penalized QR for uncensored data \citep{BC2011},  whose convergence rate is of order $\cO\big( \sqrt{s \log (p) / n}\big)$,  there is a substantial gap in terms of the impact of the sparsity parameter $s$.

In addition to the above theoretical issues, our study is motivated by the computational hardness of CQR under the framework of \cite{P2003} and \cite{PH2008} for problems with large dimension.
Recall that this framework  involves fitting a series of quantile regressions sequentially over a dense grid of quantile indexes, each of which is solvable by the Frisch-Newton algorithm with computational complexity that grows as a cubic function of $p$ \citep{PK1997}. Moreover, under the regime in which $p< n$, the asymptotic covariance matrix of the estimator is rather complicated and thus resampling methods are often used to perform statistical inference \citep{P2003, PH2008}.  A sample-based inference procedure (without resampling) for Peng-Huang's estimator \citep{PH2008} is available by adapting the plug-in covariance estimation method from \cite{SPHL2016}. In the high-dimensional setting ($p>n$), computation of the $\ell_1$-penalized QR is based on either reformulation as linear programs \citep{KN2005} or alternating direction method of multiplier algorithms \citep{YLW2017, Gu2018}. These algorithms are generic and applicable to a broad spectrum of problems but lack scalability.  Since the $\ell_1$-penalized CQR not only requires the estimation of the whole quantile regression process, but also relies on cross-validation to select the sequence of (mostly different) penalty levels,  the state-of-the-art methods  \citep{ZPH2018, FZHL2021} can be highly inefficient when applied to large-$p$ problems.

To illustrate the computational challenge for CQR,  we compare the $\ell_1$-penalized CQR proposed by \citet{ZPH2018} and our proposed method by analyzing a gene expression dataset studied in \cite{Setal2008}.
In this study,  22,283 genes from 442 lung adenocarcinomas are incorporated to predict the survival time in lung cancer,  with $46.6\%$ subjects that are censored.   We implement both methods with quantile grid set as $\{0.1, 0.11, \dots, 0.7\}$, and use a predetermined sequence of regularization parameters. For \citet{ZPH2018},  we use the \texttt{rqPen} package to compute the $\ell_1$-penalized QR estimator at each quantile level \citep{SM2020}.  The computational time and maximum allocated memory are reported  in Table~\ref{tb.time.space}. The reference machine for this experiment is  a worker node with 2.5 GHz 32-core processor and 512 GB of memory in a high-performance computing cluster.

\begin{table}[!ht] 
\begin{center}
\begin{tabular}{ c | c | c}
  \hline
Methods & Runtime & Allocated memory \\
\hline
$\ell_1$-penalized CQR & 170 hours+ & 38 GB \\
Proposed method & 2 minutes & 926 MB   \\
\hline
\end{tabular}
\caption{Computational runtime and maximum allocated memory for fitting $\ell_1$-penalized CQR and the proposed method  on the gene expression data with censored response in \cite{Setal2008}.
One gigabyte (GB)  equals 1024 megabytes (MB). } 
\label{tb.time.space}
\end{center}
\end{table}

In this paper, we develop a smoothed framework for CQR that is scalable to problems with large dimension $p$ in both low- and high-dimensional settings. Our proposed method is motivated by the smoothed estimating equation approach that has surfaced mostly in the econometrics literature \citep{W2006, WMY2015, KS2017, CGKL2019,FGH2021,HPTZ2020}, which can be applied to the stochastic integral based sequential estimation procedure proposed by \citet{PH2008} for CQR. We show in Section~\ref{sec:smooth} that the smoothed sequential estimating equations method can be reformulated as solving a sequence of optimization problems with (at least) twice-differentiable and convex loss functions for which gradient-based algorithms are available.  
Large-scale statistical inference can then be performed efficiently via multiplier/weighted bootstrap. 
In the high-dimensional setting, we propose and analyze $\ell_1$-penalized smoothed CQR estimators obtained by sequentially minimizing smoothed convex loss functions plus $\ell_1$-penalty,  which we solve using a scalable and efficient majorize-minimization-type algorithm, as evidenced in Table~\ref{tb.time.space}.

Theoretically, we provide a unified analysis for the proposed smoothed estimator in both low- and high-dimensional settings.  In the low-dimensional case where the dimension is allowed to increase with the sample size,   we establish the uniform rate of convergence and a uniform Bahadur-type representation for the smoothed CQR estimator.  We also provide a rigorous justification for the validity of a weighted/multiplier bootstrap procedure with explicit error bounds as functions of $(n, p)$.  To our knowledge, these are the first results for censored quantile regression in the increasing-$p$ regime with $p<n$.  The main challenges are as follows. To fit the QR process with censored response variables, the stochastic integral based approach entails a sequence of estimating equations that correspond to a prespecified grid of quantile indexes. A sequence of pointwise estimators can then be sequentially obtained by solving these equations. The sequential nature of this procedure poses technical challenges because at each quantile level, the objective function (or the estimating equation) depends on all of the previous estimates. To establish convergence rates for the estimated regression process, a delicate analysis beyond what is used in \cite{HPTZ2020} is required to deal with the accumulated estimation error sequentially. The mesh width of the grid should converge to zero at a proper rate in order to balance the accumulated estimation error and discretization error.
In the high-dimensional setting,  we show that with suitably chosen penalty levels and bandwidth,  the $\ell_1$-penalized smoothed CQR estimator has a uniform convergence rate of $\cO(\sqrt{s \log (p) /n})$, provided the sample size satisfies $n\gtrsim s^3 \log (p)$.  The technical arguments used in this case are also very different from those in \cite{ZPH2018} and subsequent work \cite{FZHL2021}, and as a result, our conclusion improves that of \citet{ZPH2018} by relaxing the exponential term $\exp(C s)$ in the convergence rate to a linear term in $s$. Such an improvement is significant when the effective model size $s$ is allowed to grow with $n$ and $p$ in the context of censored quantile regression.

The rest of the article is organized as follows.  In Section~\ref{sec:cqr}, we provide a formal formulation of the CQR.  We then briefly review the martingale-based estimating equation estimator proposed by \citet{PH2008} in Section~\ref{sec:penghuang}. The proposed smoothed CQR is detailed in Section~\ref{sec:smooth}, along with the multiplier bootstrap method for large-scale inference in Section~\ref{sec:mb}. 
We then provide a comprehensive theoretical analysis for the smoothed CQR estimator in Section~\ref{sec:theory} and its bootstrap counterpart. In Section~\ref{sec:highdim}, we generalize the smoothed CQR to the high-dimensional setting by incorporating a penalty function to the smoothed CQR loss  and study the theoretical properties of the regularized estimator. Extensive numerical studies and data applications are in Sections~\ref{sec:numerical} and~\ref{sec:real.data}. The \texttt{R} code that implements the proposed method is available at \href{https://github.com/XiaoouPan/scqr}{https://github.com/XiaoouPan/scqr}. 

\medskip 
\noindent
{\sc Notation}.   For any two real numbers $a$ and $b$, we write $a \wedge b = \min\{a, b\}$ and $a \vee b = \max\{a, b\}$. 
Given a pair of vectors $\bu, \bv \in \RR^p$, we use $\bu^{\T} \bv$ and $\langle \bu, \bv \rangle$ interchangeably to denote their inner product.
For a positive semi-definite matrix $\Sigma \in \RR^{p \times p}$, we define the $\Sigma$-induced $\ell_2$-norm $||\bu||_{\Sigma} = ||\Sigma^{1/2} \bu ||_2$ for any $\bu \in \RR^p$.
For every $r\geq 0$, we use $\BB^p(r) = \{\bbeta\in \RR^p : ||\bbeta||_2 \le r\}$ and $\SSS^{p - 1}(r) = \{\bbeta \in \RR^p : ||\bbeta||_2 = r\}$ to denote the Euclidean ball and sphere, respectively, with radius $r$. In particular, we write $\SSS^{p-1}= \SSS^{p-1}(1)$.
Given an event/subset $\cA$, $\1\{\cA\}$ or $\1_{\cA}$ represents the indicator function of this event/subset.
For two non-negative arrays $\{a_n\}_{n \geq 1}$ and $\{ b_n\}_{n \geq 1}$,  we write $a_n  \lesssim b_n$ if $a_n \le C b_n$ for some constant $C > 0$ independent of $n$,  $a_n \gtrsim b_n$ if  $b_n \lesssim a_n$, and $a_n  \asymp b_n$ if $a_n  \lesssim b_n$ and $a_n \gtrsim b_n$.

\section{Censored Quantile Regression}
\label{sec:cqr}

Let $z \in \RR$ be a response variable of interest,  and $\bx = ( x_1, \dots, x_p)^\T$ be a $p$-vector ($p \geq 2$) of random covariates with $x_1 \equiv 1$. 
In this work, we focus on a global conditional quantile model on $z$ described as follows.
Given a closed interval $[\tau_L, \tau_U] \subseteq (0,1)$,   assume that the $\tau$-th conditional quantile of $z$ given $\bx$ takes the form
\begin{align}
F^{-1}_{z | \bx}(\tau) = \bx^\T \bbeta^*(\tau) ~\mbox{ for any }~ \tau \in [\tau_L, \tau_U] ,   \label{cqr.model}
\end{align}
where $\bbeta^*(\tau) \in \RR^p$, formulated as a function of $\tau$, is the unknown vector of regression coefficients.

We assume that $z$ is subject to right censoring by $C$,  a random variable that is conditionally independent of $z$ given the covariates $\bx$.  Let $y = z \wedge C$  the censored outcome, and  $\Delta = \1(z \leq C)$ be an event indicator.  The observed samples $\{y_i, \Delta_i, \bx_i\}_{i = 1}^n$ consist of independent and identically distributed (i.i.d.)~replicates of the triplet $(y,  \Delta, \bx)$. In addition, we assume at the outset that the lowest quantile of interest $\tau_L$ satisfies $\PP\{ y \le   \bx^\T  \bbeta^*(\tau_L) , \Delta = 0 \} = 0$. 
This condition, interpreted as no censoring below the $\tau_L$-th quantile, is commonly imposed in the context of CQR; see, e.g., Condition~C in \cite{P2003} and Assumption~3.1 in \cite{ZPH2018}. Moreover, our quantiles of interest are confined up to $\tau_U < 1$ subject to some identifiability concerns, which is a subtle issue for CQR problems.
Briefly speaking, the model \eqref{cqr.model} may become non-identifiable as $\tau$ moves towards $1$, due to large amount of censored information in the upper tail.
In practice, determining $\tau_U$ is usually a compromise between inference range of interest and data censoring rate, and $\tau_L$ can be chosen to be close to 0 if censoring occurs at early stages. Theoretically, the above assumption on $\tau_L$ helps us simplify the technical analysis.

The above model is broadly defined,  yet it is inspired by approaching survival data with quantile regression \citep{KG2001}. To briefly illustrate, let $T$ be a non-negative random variable representing the failure time to an event.   The conditional quantile model \eqref{cqr.model} on $z = \log (T)$ can be viewed as a generalization of the standard AFT model in the sense that coefficients not only shift the location  but also affect the shape and dispersion of the conditional distributions.

\subsection{Martingale-based estimating equation estimator}
\label{sec:penghuang}

Under the global linear model \eqref{cqr.model}, two well-known methods are the recursively re-weighted estimator of \cite{P2003} and the stochastic integral based estimating equation estimator of \cite{PH2008}. Both methods are grid-based algorithms that iteratively solve a sequence of (weighted) check function minimization problems over a predetermined grid of $\tau$-values.  Motivated by the recent success of smoothing methods for  uncensored quantile regressions \citep{FGH2021,HPTZ2020,TWZ2022},  we propose a smoothed estimating equation approach for CQR in the next subsection.  We start with a brief introduction of \cite{PH2008}'s method that is built upon the martingale structure of randomly censored data.

To this end,  denote by $\Lambda_{z | \bx}(t) = -\log \{1 - \PP(z \le t | \bx)\}$ the cumulative conditional hazard function of $z$ given $\bx$, and  define the counting processes $N_{ i}(t) = \1\{y_i \le t, \Delta_i = 1\}$ and $N_{0i}(t) = \1\{y_i \le t, \Delta_i = 0\}$ for $i = 1, \dots, n$, where $\Delta_i = \1(z_i \leq C_i)$.
Define $\cF_{i }(s) = \sigma \{ N_{ i}(u),  N_{0i}(u):  u \le s\}$ as the $\sigma$-algebra generated by the foregoing processes.  Note that $\{\cF_i(s) :  s \in \RR \}$ is an increasing family of sub-$\sigma$-algebras, also known as filtration, and $N_i(t)$ is an adapted sub-martingale.
By the unique Doob-Meyer decomposition, one can construct an $\cF_i(t)$-martingale $M_i(t) = N_i(t) - \Lambda_{ z_i | \bx_i}(y_i \wedge t)$ satisfying $\EE\{ M_i(t) | \bx_i \} = 0$; see Section~1.3 of \cite{FH1991} for details.
Taking $t =  \bx_i^\T \bbeta^*(\tau)$ for each  $i$, the martingale property implies
\$
\EE \Bigg[ \sn \big\{N_i \big(  \bx_i^\T \bbeta^*(\tau)  \big) - \Lambda_{ z_i | \bx_i} \big(y_i \wedge  \bx_i^\T \bbeta^*(\tau)  \big) \big\} \bx_i \Bigg] = \bm{0} .
\$
This lays the foundation for the stochastic integral based estimating equation approach.
The monotonicity of the function $\tau \mapsto  \bx^\T \bbeta^*(\tau)$, implied by the global linearity in \eqref{cqr.model}, leads to
\begin{align*}
\Lambda_{ z_i | \bx_i} \big(y_i \wedge   \bx_i^\T \bbeta^*(\tau)  \big) = H(\tau) \wedge H\big(\PP( z_i \le y_i | \bx_i ) \big) = \int_0^\tau \1\{y_i \ge   \bx_i^\T \bbeta^*(u)   \} \,\mathrm{d}H(u) 
\end{align*}
for $\tau \in [\tau_L, \tau_U]$, where  $H(u) := -\log (1 - u)$ for $0< u < 1$. This motivates Peng and Huang's estimator \citep{PH2008}, which solves the following estimating equation
\begin{align*}
\frac{1}{n} \sn \Bigg[ N_i \big(  \bx_i^\T \bbeta(\tau)  \big) - \int_0^\tau \1\{y_i \ge  \bx_i^\T \bbeta(u)  \} \,\mathrm{d}H(u) \Bigg] \bx_i  = \bm{0} ,   ~\mbox{ for every }~ \tau_L \leq \tau \leq \tau_U.
\end{align*}

However, the exact solution to the above equation is not directly obtainable. By adapting Euler's forward method for ordinary differential equation, \cite{PH2008} proposed a grid-based sequential estimating procedure as follows. Let $\tau_L = \tau_0 < \tau_1 < \cdots < \tau_m = \tau_U$ be a grid of quantile indices. Noting that $\PP \{  y \le \bx^\T \bbeta^*(\tau_0)  , \Delta = 0 \} = 0$, we have $\EE  \int_0^{\tau_0} \1\{y_i \ge  \bx_i^\T \bbeta^*(u)  \} \,\mathrm{d}H(u) = \tau_0$, and hence $\bbeta^*(\tau_0)$ can be estimated by solving the usual quantile equation $(1/n) \sn\{N_i(  \bx_i^\T \bbeta   ) - \tau_0 \}\bx_i = \bm{0}$.
Denote $\wt{\bbeta}(\tau_0)$ as the solution to the above equation. At grid points $\tau_k$, $k = 1, \dots, m$, the estimators $\wt{\bbeta}(\tau_k)$ are sequentially obtained by solving 
\begin{align}
\frac{1}{n} \sn \Bigg[ N_i(  \bx_i^\T \bbeta   ) - \sum_{j = 0}^{k - 1} \int_{\tau_j}^{\tau_{j + 1}} \1\{y_i \ge   \bx_i^\T \wt{\bbeta}(\tau_j)   \} \, \mathrm{d}H(u) -\tau_0 \Bigg] \bx_i = \bm{0}.  \label{eq.phk}
\end{align}
The resulting estimated function $\wt{\bbeta}(\cdot) : [\tau_L,\tau_U] \mapsto \RR^p$ is  right-continuous and piecewise-constant that jumps only at each grid point. 
Computationally, solving the above equation is equivalent to  minimizing an $\ell_1$-type convex objective function after introducing a sufficiently large pseudo point to the data. The minimizer, however,  is not always uniquely defined.  To avoid this lack of uniqueness as well as grid dependence,  \cite{H2010} introduced a more general (population) integral equation, and then proposed a Progressive Localized Minimization (PLMIN) algorithm to solve its empirical version exactly. This algorithm automatically determines the breakpoints of the solution and thus is grid-free.  Under a continuity condition on the density functions (see,  e.g. condition (C2) in \cite{H2010}), the estimating functions used in \cite{PH2008} and \cite{H2010} are asymptotically equivalent.

\subsection{A smoothed estimating equation approach}
\label{sec:smooth}
Due to the discontinuity stemming from the indicator function in the counting process $N_i(\cdot)$, exact solutions to the estimating equations \eqref{eq.phk} may not exist.
In fact, $\wt{\bbeta}(\tau_j)$ for $j = 0, \dots, m$ are defined as the general solutions to generalized estimating equations \citep{FR1994}, which correspond to subgradients of some convex yet non-differentiable functions. Computationally, one may reformulate these equations as a sequence of linear programs,  solvable by the Frisch-Newton algorithm described in \cite{PK1997}. The computation complexity grows rapidly when the dimensionality $p$ increases. To mitigate the computational burden of the existing methods, we employ a smoothed estimating equation (SEE) approach for fitting large-scale censored quantile regression models.

Let $K(\cdot)$ be a symmetric and non-negative kernel function and let $\bar{K}(u) = \int_{-\infty}^u K(x)\,\mathrm{d}x$, which is a non-decreasing function that is between 0 and 1. The non-smooth indicator function $\mathbbm{1}(u \geq 0)$ can thus be approximated by $\bar K(u/h)$ for some $h>0$ in the sense that  as $h\to 0$, $\bar K (u/h) \to 1$ for $u\geq 0$ and $\bar K (u/h) \to 0$ for $u<0$. Hereinafter, $h>0$ will be referred to as a bandwidth. As the aforementioned, let $\tau_L = \tau_0 <  \tau_1 < \dots < \tau_m = \tau_U$ be a grid of quantile indices for some $m\geq 1$.  Given a kernel function $K(\cdot)$ and a bandwidth $h>0$,  write 
\$
	K_h(u) = h^{-1} K(u/h) ~~\mbox{ and }~~ \bar K_h( u ) = \bar K(u/h) = \int_{-\infty}^{u/h}  K(v) {\rm d}v , \ \ u \in \RR , 
\$
so that $\bar K_h'(u) = K_h(u)$.
We now propose a smooth SEE approach for CQR. 

\begin{enumerate}
\item At $\tau= \tau_0$, we estimate $\bbeta^*(\tau_0)$ by $  \hat \bbeta(\tau_0)$, obtained from solving $\hat{Q}_0(\bbeta) = \bm{0}$, where
\begin{align}
\hat{Q}_0(\bbeta) := \frac{1}{n} \sn \big\{\Delta_i \bar{K}_h  ( - r_i(\bbeta)  ) - \tau_0 \big\}\bx_i ~~\mbox{ and }~~r_i(\bbeta) = y_i - \bx_i^\T \bbeta .  \label{eq.tau0}
\end{align}
\item At grid points $\tau_k$ for $k = 1, \dots, m$, set $\hat \bbeta(\tau) = \hat \bbeta(\tau_{k - 1})$ for any $\tau \in (\tau_{k - 1}, \tau_k)$, and then obtain estimators $  \hat \bbeta(\tau_k)$ of $\bbeta^*(\tau_k)$ by solving $\hat{Q}_k(\bbeta) = \bm{0}$, where
\begin{align}
\hat{Q}_k(\bbeta) := \frac{1}{n} \sn \Bigg[ \Delta_i \bar{K}_h ( - r_i(\bbeta)   ) - \sum_{j = 0}^{k - 1}  \bar K_h  ( r_i(\hat{\bbeta}(\tau_j)   ) \{ H(\tau_{j+1} ) - H(\tau_j) \}  - \tau_0 \Bigg] \bx_i.  \label{eq.tauk}
\end{align}
\end{enumerate}
Note that the resulting estimator $\hat {\bbeta}(\cdot) : [\tau_L,\tau_U] \mapsto \RR^p$ is right-continuous and piecewise-constant with jumps only at grids.
For notational convenience, throughout the remainder of this paper we write
\$
	\bbeta^*_k = \bbeta^*(\tau_k) ~~\mbox{ and }~~ \hat{\bbeta}_k = \hat{\bbeta}(\tau_k) , \quad k=0,1,\ldots, m .
\$

Before proceeding,  it is worth noticing that the above smoothed estimating equations method is closely related to the convolution smoothing approach studied in \cite{FGH2021} and \cite{HPTZ2020}.  Consider the check function $\rho_\tau(u) =  \tau \{ u - \mathbbm{1}(u<0) \}$, and  its convolution smoothed counterpart 
\$
	\ell_{\tau, h} (u) = (\rho_\tau * K_h) (u)  = \int_{-\infty}^\infty \rho_\tau(v) K_h(v-u) \, {\rm d}v ,
\$
where $*$ denotes the convolution operator. Given censored data $\{ (y_i, \Delta_i, \bx_i)\}_{i=1}^n$, define the empirical smoothed loss
\#
	\hat L_0 (\bbeta) = \frac{1}{n} \sn  \bigl\{  \Delta_i \ell_{\tau_0 , h} (y_i - \bx_i^\T \bbeta)   + \tau_0 (\Delta_i - 1) \bx_i^\T \bbeta  \bigr\}  , \label{def:L0}   
\#
whose gradient and Hessian are 
$$
	\nabla \hat L_0(\bbeta) = \hat Q_0(\bbeta)~~\mbox{ and }~~ \nabla^2 \hat L_0(\bbeta) =\frac{1}{n} \sn \Delta_i K_h(r_i(\bbeta)) \bx_i \bx_i^\T ,
$$
respectively.  Hence, the foregoing estimator $\hat \bbeta_0$ can be equivalently defined as the solution to the (unconstrained) optimization problem $\min_{\bbeta \in \RR^p}  \hat L_0 (\bbeta)$.
When a non-negative kernel is used, the objective function $\hat L_0 (\cdot)$ is convex, and thus any minimizer satisfies the first-order condition. At subsequent grid points $\tau_k$ for $k=1,\ldots, m$, the estimator $\hat \bbeta_k$ can also be viewed as an $M$-estimator that solves
\#
	\min_{\bbeta \in \RR^p} ~ \Bigg\{  \hat L_k(\bbeta) :=  \hat L_0(\bbeta)  -  \Bigg\langle  \frac{1}{n} \sn   \sum_{j = 0}^{k - 1} \bar K_h   (  y_i -\bx_i^\T \hat{\bbeta}_j )  \{ H(\tau_{j+1} )- H(\tau_j) \}  \bx_i    ,   \bbeta   \Bigg\rangle  \Bigg\}. \label{def:Lk} 
\#
Notably, kernel smoothing produces continuously differentiable  estimating functions $\hat{Q}_k(\cdot)$ ($k=0,\ldots, m$), or equivalently, convex and twice-differentiable loss functions $\hat L_k(\cdot)$, which have the same  positive semi-definite Hessian matrix $\nabla^2 \hat L_k(\bbeta) =  (1/n)  \sn \Delta_i K_h(    \bx_i^\T \bbeta - y_i ) \bx_i \bx_i^\T$.
As we shall see, the empirical loss functions $\hat L_k(\cdot)$ are not only globally convex but also locally strongly convex (with high probability).
This property ensures the existence of global solutions to the sequential estimation problems, which can be efficiently solved by a quasi-Newton algorithm described in Section~A.1 of the supplementary material.

\subsection{Inference with bootstrapped process}
\label{sec:mb}

In this subsection, we construct component-wise confidence intervals for $\bbeta^*(\tau)$ at some quantile index $\tau$ of interest by bootstrapping the quantile process. 
Recall that $\hat{\bbeta}_k$'s are the solutions to the equations $\hat{Q}_k(\bbeta) = \bm{0}$, where $\hat{Q}_k(\cdot)$ ($k=0,1,\ldots,m$) are defined in \eqref{eq.tau0} and \eqref{eq.tauk}.
Analogously, we construct bootstrap estimators $\hat \bbeta^\flat_k$ following a sequential procedure based on the bootstrapped SEEs obtained by perturbing $\hat{Q}_k(\cdot)$ with random weights. 
Independent of the observed data $\{y_i, \Delta_i, \bx_i\}_{i = 1}^n$,  let $W_1, \ldots, W_n$ be exchangeable non-negative random variables, satisfying $\EE(W_i) = 1$ and $\var(W_i)>0$.
The bootstrap estimators can be constructed as follows:  

\begin{enumerate}
\item Set $\hat \bbeta^\flat_0$ as the solution of $\hat{Q}^\flat_0(\bbeta) = \bm{0}$, where
\#
\hat{Q}^\flat_0(\bbeta) &:= \frac{1}{n} \sn W_i \big\{\Delta_i \bar{K}_h  ( - r_i(\bbeta)  ) - \tau_0 \big\}\bx_i ~\mbox{ with }~ r_i(\bbeta) = y_i - \bx_i^\T \bbeta.   \label{eq.tau0.boot}
\# 

\item For $k = 1, \dots, m$,  compute $\hat \bbeta^\flat_k$ sequentially by solving $\hat{Q}^\flat_k(\bbeta) = \bm{0}$, where
\#
\hat{Q}^\flat_k(\bbeta) &:= \frac{1}{n} \sn W_i \Bigg[  \Delta_i \bar{K}_h ( -  r_i(\bbeta) ) - \sum_{\ell = 0}^{k - 1}  \bar K_h (  r_i(\hat{\bbeta}^\flat_\ell) ) \{ H(\tau_{\ell+1 } - H(\tau_\ell)  \} - \tau_0 \Bigg]  \bx_i. \label{eq.tauk.boot}
\#

\item  Define the bootstrap estimate of the coefficient process $\hat \bbeta^\flat(\cdot) : [\tau_L, \tau_U ] \mapsto \RR^p$ as  $\hat \bbeta^\flat(\tau) = \hat \bbeta^\flat_{k - 1}$ for $\tau \in [\tau_{k - 1}, \tau_k)$ and $k=1,\ldots, m$.
\end{enumerate}

For a prescribed nominal level,  we can construct component-wise percentile or normal-based confidence intervals for $\bbeta_j^*(\tau)$ ($j=1,\ldots, p$).
The above multiplier bootstrap estimator $\hat {\bbeta}^\flat(\cdot) : [\tau_L,\tau_U] \mapsto \RR^p$ of the coefficient process behaves similarly as $\hat{\bbeta}(\cdot)$, in the sense that they are both right-continuous and piecewise-constant with jumps only at the grids.  
The multiplier bootstrap method, which dates back at least to \cite{BB1995}, is motivated by the following simple yet important observation. Let $\EE^*(\cdot)$ be conditional expectation given the data, i.e., $\EE^*(\cdot) = \EE(\cdot | \, \{y_i, \Delta_i, \bx_i\}_{i = 1}^n)$.  Since $\EE(W_i) = 1$, we have $\EE^*\{\hat{Q}^\flat_0(\bbeta)\} = \hat{Q}_0(\bbeta)$ and $\EE^*\{\hat{Q}^\flat_k(\bbeta)\} \approx \hat{Q}_k(\bbeta)$ for $k = 1, \dots, m$. This means that in the bootstrap world, $\hat{Q}_k^\flat(\bbeta)$ can be viewed as an empirical version of $\hat{Q}_k(\bbeta)$, and thus $\hat{\bbeta}_k^\flat$ can be regarded as the bootstrap estimator of $\hat{\bbeta}_k$.

We complete this section with a  brief discussion of other resampling methods for quantile regression. Given the random weights $\{ W_i \}_{i = 1}^n$ independent of data,  another available approach is to minimize the randomly perturbed objective functions \citep{JYW2001, PH2008}.  In the current setting,  it seems more natural to directly bootstrap the estimating equations. In terms of bootstrapping estimating equations with uncensored data,  \cite{PWY1994}'s method is based on the assumption that the estimating equation is exactly or asymptotically pivotal,  and \cite{HK2000}'s proposal is based on resampling with replacement.
A generalized weighted bootstrap and its asymptotic theory has been rigorously studied in \cite{CB2005} and \cite{MK2005}.  For censored quantile regression,  the sequential SEEs \eqref{eq.tauk} are not directly formulated as empirical averages of independent random quantities,  nor do they satisfy the required assumptions in the literature; see Section 2 of \cite{PWY1994}, Section 2 of \cite{HK2000}, and Section 3 of \cite{CB2005}.
Hence,  the validity of weighted bootstrap for CQR is of independent interest,  and will be examined in Section~\ref{sec:fst.boot}.

\begin{remark}
In practice, random weights $\{ W_i \}_{i=1}^n$ can be generated from one of the following distributions.  (i) $(W_1, \ldots, W_n) \sim {\rm Multinomial}(n, 1/n, \ldots, 1/n)$.  This leads to Efron's nonparametric bootstrap, for which the random weights are exchangeable but not independent; (ii) $W_1, \ldots, W_n \sim {\rm Exp}(1)$ are i.i.d.~exponentially distributed random variables; and (iii) $W_i = e_i+1$, where $e_i$'s are  i.i.d.~Rademacher random variables, defined by $\PP(e_i = 1) = \PP(e_i = 0) = 1/2$. We refer to this as the Rademacher multiplier bootstrap. Its theoretical properties will be investigated in Section~\ref{sec:fst.boot}.
\end{remark}

\section{Theoretical Analysis}
\label{sec:theory}

\subsection{Regularity conditions}
\label{sec:not.cond}

We first impose some technical assumptions required for the results in Sections~\ref{sec:fst} and \ref{sec:fst.boot}.

\begin{cond}[Kernel function] 
\label{cond:kernel} 
Let $K(\cdot)$ be a symmetric, Lipschitz continuous and non-negative kernel function, that is, $K(u)= K(-u)$, $K(u) \geq 0$ for all $u \in \RR$ and $\int_{-\infty}^\infty K(u) \, {\rm d}u  = 1$. Moreover, $\kappa_u = \sup_{u \in \RR } K(u)<\infty$, $\kappa_l = \min_{|u| \le c}K(u) > 0$ for some $c>0$. 
We define its higher-order absolute moments as $\kappa_\ell = \int_{-\infty}^\infty |u|^\ell K(u)   {\rm d} u$ for  any positive integer $\ell$.
\end{cond}

 \begin{cond}[Random design] 
 \label{cond:x.sg}
 The random covariate vector $\bx = (x_1, \ldots, x_p)^\T \in \cX \subseteq  \RR^p$ is compactly supported with 
$\zeta_p := \sup_{\bx \in \cX} \| \Sigma^{-1/2} \bx \|_2   < \infty$,  where $\Sigma = \EE(\bx \bx^\T)$ is positive definite. 
\end{cond}

\begin{cond}[Conditional densities] 
 \label{cond:density}
 Assume $(z, \bx)$ follows the global conditional quantile model \eqref{cqr.model}.
Define the conditional cumulative distribution functions $F_z(u|\bx) = \PP(z \le u | \bx)$, $F_y(u| \bx ) = \PP(y \le u | \bx)$ and $G(u| \bx) = \PP(y \le u, \Delta = 1 | \bx)$, where $y=z \wedge C $ and $C$ is independent of $z$ given $\bx$.
Assume that the conditional densities $f_z(u| \bx)=F_z'(u| \bx)$,  $f_y(u| \bx)=F_y'(u| \bx)$ and $g(u| \bx)=G'(u| \bx)$ exist, and satisfy almost surely (over $\bx$) that
\begin{align*}
   \inf_{ \tau \in [\tau_L, \tau_U] }  \min\big\{ f_y( \bx^\T  \bbeta^*(\tau) | \bx ) , f_z( \bx^\T  \bbeta^*(\tau) | \bx )   \big\} \geq   \underline{f} >0, \quad   \sup_{ u \in \RR} f_y(u | \bx) \le \overline f, \\
  0<  \underline{g}  \le \inf_{  |u-  \bx^\T    \bbeta^*(\tau)  | \leq 1/2 , \tau \in  [\tau_L, \tau_U]} g( u | \bx ) \le \sup_{  u \in \RR} g(u | \bx) \le \overline g. 
\end{align*}
Moreover, there exists a constant $l_1 > 0$ such that for any $u \in \RR$,
\begin{align*}
&\sup_{\bx \in \RR^p, \tau \in [\tau_L, \tau_U] } | f_y( \bx^\T  \bbeta^*(\tau)  + u | \bx) - f_y(\bx^\T  \bbeta^*(\tau) | \bx)| \le l_1 |u|, \\
&\sup_{\bx \in \RR^p, \tau \in [\tau_L, \tau_U] } |g( \bx^\T  \bbeta^*(\tau)   +u | \bx) - g(\bx^\T  \bbeta^*(\tau) | \bx )| \le l_1 |u|.
\end{align*}
\end{cond}

\begin{cond}[Grid size]
\label{cond:grid}
The grid of quantile levels $\tau_L = \tau_0 < \tau_1 < \dots < \tau_m = \tau_U$ satisfies $n^{-1} \leq  \delta_* \le \delta^* \lesssim n^{-1/2}$, where $\delta^* = \max_{1 \le k \le m}(\tau_k - \tau_{k - 1})$ and $\delta_* = \min_{1 \le k \le m}(\tau_k - \tau_{k - 1})$.
\end{cond}

Condition~\ref{cond:kernel} holds for most commonly used kernel functions, including: (a) uniform kernel $K(u) = (1/2) \mathbbm{1} (|u|\leq 1)$, (b) Gaussian kernel  $K(u) = (2\pi)^{-1/2} e^{-u^2/2}$, (c)  logistic kernel $K(u) = e^{-u}/(1+ e^{-u})^2 $, (d) Epanechnikov/parabolic kernel $K(u) = (3 / 4) (1 - u^2) \mathbbm{1}(|u|\leq 1)$, and (e) triangular kernel $K(u) = (1 - |u|) \mathbbm{1}(|u|\leq 1)$. 
To simplify the analysis, we take $c=1$ in Condition~\ref{cond:kernel}; otherwise if $c<1$ and $K(\pm 1) =0$, we can simply use a re-scaled kernel $K_c(u) := c K(cu)$, so that $\min_{|u|\leq 1} K_c(u) = c \min_{|u|\leq c} K(u)$. The compactness of $\cX$ in Condition~\ref{cond:x.sg} is a common requirement for a global linear quantile regression model (quantile regression process) \citep{K2005}. If the support of the covariate space---the set of $x_j$'s that occur with positive probability---is unbounded,  at some points there will be 	``crossings'' of the conditional quantile functions, unless these functions are parallel, which corresponds to a pure location-shift model. The quantity $\zeta_p$ plays an important role in the theoretical results. Alternatively, one may assume $\| \Sigma^{-1/2} \bx \|_\infty\leq C_0$ (almost surely) as in \cite{ZPH2018}, which in turn implies $\zeta_p \leq C_0 p^{1/2}$ in the worst-case scenario. In general, it is reasonable to assume that $\zeta_p \asymp p^{1/2}$. In addition to $\zeta_p$,   define the moment parameters
\#
	m_q = \sup_{\bu\in\mathbb{S}^{p-1}} \EE \big(  | \bu^\T \Sigma^{-1/2} \bx |^q \big)  ~~\mbox{ for }~~ q=3, 4 ,  \label{def:m34}
\# 
which satisfy the worst-case bounds $m_3 \leq \zeta_p$ and $m_4\leq \zeta_p^2$.

 Conditions~\ref{cond:x.sg} and \ref{cond:density} ensure that the coefficient function $\bbeta^*(\cdot)$ is Lipschitz continuous.  Since $\bbeta^*(\tau)$ solves the equation $ \EE[  \{  \tau - \1(z \leq \bx^\T \bbeta ) \} \bx ] = \textbf{0}$,  we have $\frac{d}{d \tau } \bbeta^*(\tau)  =  \EE \{  f_z( \bx^\T \bbeta^*(\tau) | \bx ) \bx \bx^\T  \}^{-1}  \EE(\bx)$.  Under Condition~\ref{cond:x.sg}, it holds
$$\max_{ \tau \in [\tau_L, \tau_U]} \bigg\|  \frac{d}{d \tau } \Sigma^{1/2} \bbeta^*(\tau)   \bigg\|_2 \leq  \underline{f}^{-1} \max_{ \tau \in [\tau_L, \tau_U]}  \| \EE(\Sigma^{-1/2} \bx ) \|_2 \leq \underline{f}^{-1}, $$
which, together with the mean value theorem, implies 
\#
 \|  \bbeta^*(\tau) - \bbeta^* (\tau') \|_{\Sigma} \le  \underline{f}^{-1}  |\tau - \tau' |  ~~~\mbox{for any}~~ \tau, \tau' \in [\tau_L, \tau_U]. \label{beta.lip}
\#
By the definitions in Condition~\ref{cond:density}, $G(u | \bx) \le F(u| \bx)$ for any $u > 0$.  Recall that we have assumed no censored observations at low quantile levels $\tau \leq \tau_L$. Hence, $G( \bx^\T  \bbeta^*(\tau_L) | \bx ) = F(\bx^\T  \bbeta^*(\tau_L)| \bx ) = \tau_L$, and $G(\bx^\T  \bbeta^*(\tau) | \bx) \le \tau \le F(\bx^\T  \bbeta^*(\tau)| \bx)$ for $\tau_L < \tau \leq \tau_U$.  Condition~\ref{cond:grid} assures a fine grid by controlling the gap between two contiguous points,  so that the approximation/discretization error does not exceed the statistical error.

\subsection{Uniform rate of convergence and Bahadur representation}
\label{sec:fst}

In this section, we characterize the statistical properties of the SEE estimators for censored quantile regression with growing dimensions. That is,  the dimension $p = p_n$ is subject to the growth condition $p \asymp n^a$ for some $a \in (0,1)$. Our first result provides the uniform rate of convergence for the estimated coefficient function $\hat \bbeta(\cdot)$ under mild bandwidth constraints.

\begin{theorem}[Uniform consistency]
\label{thm:concen}
Assume Conditions~\ref{cond:kernel}--\ref{cond:grid} hold, and choose the bandwidth $h=h_n  \asymp \{ (p+\log n) / n \}^\gamma$ for some $\gamma \in [1/4, 1/2)$. Further let $n \gtrsim \{ \zeta_p^2 (p + \log n)^{1/2-\gamma} \}^{1/(1-\gamma)}$.  Then, the SEE estimator $\hat{\bbeta}(\cdot) :[\tau_L, \tau_U] \mapsto \RR^p$ satisfies 
\#
\sup_{\tau \in [\tau_L, \tau_U] }  \| \hat{\bbeta}(\tau) - \bbeta^*(\tau) \|_\Sigma \lesssim  \bigg( \frac{1-\tau_L}{1-\tau_U} \bigg)^{C_0 \overline f/\underline{g}}  \underline{g}^{-1}  \sqrt{\frac{p + \log n}{n}}    \label{eq:thm,concen}
\#
with probability at least $1 - C_1 n^{-1}$, where $C_0, C_1> 0$ are constants independent of $(n, p)$. 
\end{theorem}

Since the deviation bound in \eqref{eq:thm,concen}  depends explicitly on $n, p$ as well as other model parameters, this non-asymptotic result implies the classical asymptotic consistency by letting $n\to \infty$ with $p$ fixed. From an asymptotic perspective, Theorem~\ref{thm:concen} implies that the smoothed estimator with a bandwidth $h=h_n \asymp \{ \log(n) /n \}^\gamma$ for some $\gamma \in [1/4, 1/2)$ satisfies $\sup_{\tau_L \leq \tau \leq \tau_U } \| \hat   \bbeta(\tau) - \bbeta^*(\tau) \|_2 \rightarrow  0$ in probability as $n\to \infty$.

Recall that in the sequential estimation procedure described in Section~\ref{sec:smooth},  the $j$-th estimator $\hat{\bbeta}_j$ $(j \ge 1)$ depends implicitly on its predecessors through the estimating function  \eqref{eq.tauk}. 
In other words, the accumulative estimation errors of  $\hat\bbeta(\tau)$ for $\tau_L \le \tau < \tau_j$ may have a non-negligible impact on $\hat{\bbeta}_j = \hat \bbeta(\tau_j )$. The next result explicitly quantifies this accumulative error.
For $ \tau \in [\tau_L, \tau_U]$,  define $p\times p$ matrices
\begin{align}
	\Jb(\tau) = \EE    \big\{  g (  \bx^\T  \bbeta^*(\tau) | \bx  ) \bx \bx^\T  \big\} ~~\mbox{ and }~~ \Hb(\tau) =  \EE    \big\{  f_y (  \bx^\T  \bbeta^*(\tau) | \bx  ) \bx \bx^\T  \big\}  ,   \label{def.JH}
\end{align}
both of which are positive definite under Conditions~\ref{cond:x.sg} and \ref{cond:density}.
Moreover,  define the integrated covariate effect and its estimate 
\$
  \bbeta^*_{\inte}(\tau)  &:=  \Jb(\tau) \bbeta^*(\tau) + \int_{\tau_L}^{\tau} \Hb(u) \bbeta^*(u) \, \mathrm{d} H(u)   \\
   \mbox{and }~~ \hat{\bbeta}_{\inte}(\tau) &:=  \Jb(\tau) \hat \bbeta(\tau) + \int_{\tau_L}^{\tau} \Hb(u) \hat \bbeta(u)  \,\mathrm{d} H(u)  ,
\$
respectively,  so that $\hat \be(\tau) := \hat{\bbeta}_{\inte}(\tau) - \bbeta^*_{\inte}(\tau)$ can be interpreted as the accumulated error in the sequential estimation procedure up to $\tau$.
That is,
\#
\hat \be(\tau)=  \underbrace{ \Jb(\tau) \{ \hat \bbeta(\tau) - \bbeta^*(\tau) \}}_\text{current step}  + \underbrace{\int_{\tau_L}^{\tau} \Hb(u) \{ \hat \bbeta(u) - \bbeta^*(u) \}   \, \mathrm{d} H(u)}_\text{preceding steps}.  \label{def:accu.error}
\#
The following theorem provides a uniform Bahadur representation for $\hat{\be}(\cdot)$.

\begin{theorem}[Uniform Bahadur representation]
\label{thm:bahadur}   
Assume that the same set of conditions in Theorem~\ref{thm:concen} hold.  Moreover,  assume $\delta^* \asymp   n^{-(1/2 + \alpha)}$ for some $\alpha \in (0, 1/2)$. Then,  the SEE estimator $\hat{\bbeta}(\cdot) :[\tau_L, \tau_U] \mapsto \RR^p$ satisfies
\#
\hat \be(\tau) =  \hat{\bbeta}_{\inte}(\tau) - \bbeta^*_{\inte}(\tau) = \frac{1}{n} \sn \bU_i(\tau) + \br_{n}(\tau),   \label{eq:thm,bahadur}
\#
where
\#
\bU_i(\tau) := \Bigg\{ \tau_L +  \int_{\tau_L}^{\tau} \bar K_h  ( y_i - \bx_i^\T  {\bbeta}^*(u) ) \, \mathrm{d}H(u) - \Delta_i \bar K_h (  \bx_i^\T \bbeta^*(\tau) - y_i ) \Bigg\} \bx_i   \label{bahadur.ui}
\#
satisfies $\sup_{\tau \in [\tau_L, \tau_U]}  \|   \EE \bU_i(\tau) \|_{\Sigma^{-1}} \lesssim h^2$,
and  the remainder process $ \br_{n}(\cdot) : [\tau_L, \tau_U] \mapsto \RR^p$ is such that
\#
\sup_{\tau \in [\tau_L, \tau_U]}  \| \br_{n}(\tau) {\|}_{\Sigma^{-1}} \lesssim m_4^{1/2} \frac{p + \log n}{ nh^{1/2}} + m_3  \frac{p + \log n}{ n} + h \sqrt{\frac{p + \log n}{n}} + n^{-1/2 - \alpha}   \label{eq:prob,bahadur}
\#
with probability at least $1-C_2 n^{-1}$ for some absolute constant $C_2 > 0$, where $m_q$ ($q = 3, 4$) are given in \eqref{def:m34}. 
\end{theorem}

\begin{remark}  \label{remark:thm3.2}
Together, the above uniform Bahadur representation and the production integration  theory \citep{GJ1990}  establish the asymptotic distribution of $ \hat{\bbeta}(\cdot)$.
Define 
$$
	\btheta^*(\tau) = \Jb(\tau) \bbeta^*(\tau), \quad \hat \btheta(\tau) = \Jb(\tau) \hat \bbeta(\tau) ~~\mbox{ and }~~ \bPsi(\tau) = \frac{1}{1-\tau} \Hb(\tau) \Jb(\tau)^{-1} , \quad \tau \in [\tau_L, \tau_U].
$$
Then, equation \eqref{def:accu.error} reads $ \hat \be(\tau) =  \hat \btheta(\tau) -  \btheta^*(\tau)  + \int_{\tau_L}^\tau  \bPsi(u)  \{  \hat \btheta(u) -  \btheta^*(u) \}  {\rm d}u$.  Combined with Theorem~\ref{thm:bahadur},  this implies
\#
  n^{1/2} \{ \hat \btheta(\tau) -  \btheta^*(\tau) \}   \, + & \int_{\tau_L}^\tau  \bPsi(u)  \, n^{1/2}\{  \hat \btheta(u) -  \btheta^*(u) \}  \,{\rm d}u \nn  \\
  &  = \frac{1}{n^{1/2}} \sn \{  \bU_i(\tau) - \EE \bU_i(\tau)  \} +  \bar \br_n(\tau) , \quad \tau \in [\tau_L, \tau_U],  \label{sde}
\#
where the rescaled remainder $\bar \br_n(\cdot)$ satisfies $\sup_{\tau\in [\tau_L, \tau_U] } \| \bar \br_n (\tau)  \|_2= o_{\PP}(1)$, with a properly chosen bandwidth that will be discussed in Remark~\ref{remark:bandwidth}.
Note that equation \eqref{sde} is a stochastic differential equation for $n^{1/2} \{ \hat \btheta(\tau) - \btheta^*(\tau)\}$ \citep{PH2008}.  From the classical production integration theory (\cite{GJ1990} and Section~II.6 of \cite{ABGK1993}),  it follows that
\#
n^{1/2} \{ \hat \btheta(\tau) -  \btheta^*(\tau) \}   = \bphi\Bigg( \frac{1}{n^{1/2}} \sn \{  \bU_i(\tau) - \EE \bU_i(\tau)  \}  \Bigg)  + o_{\PP}(1) , \label{sde.solution}
\#
where $ \bphi$ is a linear operator from $\cF$ to $\cF$ defined as 
\#
	  \bphi(\bg) (\tau) =    \PPi_{u \in [\tau_L, \tau]}    \big\{ \Ib_p - \bPsi(u)  {\rm d} u  \big\} \bg(\tau_L)  + \int_{\tau_L}^\tau  \PPi_{u \in (s, \tau ]}    \big\{ \Ib_p - \bPsi(u)  {\rm d} u  \big\} \, {\rm d} \bg (s)    \label{def.linear.ope}
\#
for $\bg \in \cF :=   \{ \bm{f}: [\tau_L, \tau_U] \to \RR^p \,|\, \bm{f} \mbox{ is left-continuous with right limit} \}$, and $\PPi$ denotes the product-limit; see Definition~1 in \cite{GJ1990}. 
After careful proofreading, we believe  that the above form of $\bphi(\cdot)$ corrects an error (possibly a typo) in the proof of Theorem~2 in \cite{PH2008}; see the arguments between (B.1) and (B.3) therein. Specifically, the linear operator $ \bphi$ in \cite{PH2008} reads 
$$
 \bphi(\bg) (\tau) =    \PPi_{u \in [\tau_L, \tau]}    \big\{ \Ib_p + \bPsi(u)  {\rm d} u  \big\} \bg(\tau_L)  + \int_{\tau_L}^\tau  \PPi_{u \in (s, \tau ]}    \big\{ \Ib_p + \bPsi(u)  {\rm d} u  \big\} \, {\rm d} \bg (s)  .
$$

The asymptotic distribution of $n^{1/2} \{ \hat \btheta(\tau) -  \btheta^*(\tau) \}$ or its linear functional is thus determined by that of 
\$
 \bphi\Bigg( \frac{1}{n^{1/2}} \sn \{  \bU_i(\tau) - \EE \bU_i(\tau)  \}  \Bigg)   ~~\mbox{ and }~~    \frac{1}{n^{1/2}} \sn \{  \bU_i(\tau) - \EE \bU_i(\tau)  \}  .
\$
\end{remark}

\begin{remark}[Order of bandwidth]
\label{remark:bandwidth}
We further discuss the order of bandwidth $h$,  as a function of $(n,p)$,  required in Theorem~\ref{thm:bahadur} and Remark~\ref{remark:thm3.2}.
Following \eqref{sde}, if the moment parameters $m_3$ (absolute skewness) and $m_4$ (kurtosis) are dimension-free,  the Bahadur linearization remainder $\bar \br_n(\cdot)$ satisfies with high probability that  $\sup_{\tau \in [\tau_L, \tau_U]} \| \bar \br_n(\tau ) \|_{\Sigma^{-1}} \lesssim   n^{1/2}h^2 + (p+ \log n)/(n h)^{1/2}  + n^{- \alpha}$.
Set the bandwidth $h \asymp \{(p + \log n) / n\}^{\gamma}$ for some $\gamma \in [1/4, 1/2)$,  this implies
\#
\sup_{\tau \in [\tau_L, \tau_U]} \| \bar \br_n(\tau ) \|_{\Sigma^{-1}} \lesssim \frac{(p + \log n)^{2\gamma}}{n^{2\gamma - 1/2}} + \frac{(p + \log n)^{1 - \gamma / 2}}{n^{1/2 - \gamma / 2}} + \frac{1}{n^{\alpha}} = o_{\PP}(1), \nn
\#
provided that $p = o(n^{1 - 1 / (4 \gamma)} \wedge n^{(1 - \gamma) / (2 - \gamma)})$.
In particular,  letting $1 - 1 / (4 \gamma) = (1 - \gamma) / (2 - \gamma)$ yields $\gamma=2/5$. We therefore choose the bandwidth $h \asymp \{(p + \log n) / n\}^{2/5}$, so that all the asymptotic results (from uniform rate of convergence to Bahadur representation) hold under the growth condition $p=o(n^{3/8})$ of dimensionality $p$ in sample size $n$.
\end{remark}

Theorem~\ref{thm:bahadur} explicitly characterizes the leading term  of the integrated estimation error \eqref{def:accu.error},  along with a high probability bound on the remainder process. 
As discussed in Remark~\ref{remark:thm3.2}, the asymptotic distributions of $n^{1/2}  \{ \hat{\bbeta}(\tau) - \bbeta^*(\tau) \}$ or its linear functional can be established based on the  stochastic integral representation \eqref{sde.solution}, which further depends on the centered random process $n^{-1/2} \sn \{ \bU_i(\cdot) - \EE \bU_i (\cdot) \}$.
Let  $\{\ba_n\}_{n = 1}^\infty$ be a sequence of deterministic vectors in $\RR^p$, and define 
\#
\GG_n(\tau) :=  \frac{1}{ n^{1/2} } \sn \langle  \ba_n / || \ba_n||_{\Sigma} ,  \bU_i(\tau)  - \EE \bU_i(\tau) \rangle ,   \quad \tau \in [\tau_L, \tau_U].  \label{def.Gn}
\#
The asymptotic behavior of $\{\GG_n(\tau)  :  \tau \in [\tau_L, \tau_U]\}$ is provided in the following result.

\begin{theorem}[Weak convergence]
\label{thm:weak.convergence}
Assume Conditions~\ref{cond:kernel}--\ref{cond:grid} hold with $\delta^* \asymp  n^{-(1/2 + \alpha)}$ for some $\alpha \in (0, 1/2)$.  Moreover,  assume $h \asymp \{(p+\log n)/n\}^{2/5}$ and $p = o(n^{3/8})$ as $n\to \infty$.  For any deterministic sequence of vectors $\{\ba_n\}_{n\geq 1}$,  if the following limit 
\#
H(\tau ,  \tau') := \lim_{n \to \infty}  \frac{1}{ \| \ba_n \|_{\Sigma}^2} \ba_n^\T \, \EE \big\{ \bU_i(\tau ) \bU_i(\tau' )^\T \big\} \ba_n   \label{limit.cov}
\#
exists for any $\tau, \tau'  \in [\tau_L, \tau_U]$ with $\bU_i(\cdot)$ defined in \eqref{bahadur.ui}, then
\#
 \GG_n(\cdot) \rightsquigarrow \GG(\cdot) ~~~~ \mbox{in } \ell^\infty([\tau_L, \tau_U]), \label{eq:weak.convergence}
\#
where $\GG_n(\cdot)$ is given in \eqref{def.Gn},  and $\GG(\cdot)$ is a tight zero-mean Gaussian process with covariance function $H(\cdot, \cdot)$ and has almost surely continuous sample paths.
\end{theorem}

Regarding the relative efficiency of the proposed SEE estimator compared to its non-smoothed counterpart \citep{PH2008},  note that the (integrated) kernel $\bar K_h(u) $ converges to $\mathbbm{1}(u\geq 0)$ as $h\to 0$. Hence, the smoothed process $n^{-1/2}\sn \bU_i(\tau)$ with $\bU_i(\tau)$ given in \eqref{bahadur.ui} has the same asymptotic distribution as 
$$
	\frac{1}{\sqrt{n}} \sn \Bigg\{ \tau_L +  \int_{\tau_L}^{\tau} \mathbbm{1} ( y_i \geq  \bx_i^\T  {\bbeta}^*(u) ) \, \mathrm{d}H(u) - \Delta_i \mathbbm{1}( y_i \leq   \bx_i^\T \bbeta^*(\tau)   ) \Bigg\} \bx_i .
$$
As a result, the covariance function $H(\cdot, \cdot)$ defined in \eqref{limit.cov} coincides with that in \cite{PH2008}; see the proof of Theorem~2 therein. In other words, the SEE estimator and Peng and Huang's estimator converge to the same Gaussian process as $n\to \infty$ with $p$ fixed, and  thence the asymptotic relative efficiency is 1. The technical devices required to deal with the fixed-$p$ and growing-$p$ cases are quite different.  For the former, the consistency follows from the Glivenko-Cantelli theorem, and the weak convergence is a consequence of Donsker's theorem. To establish non-asymptotic results, we rely on a localized analysis as well as a (local) restricted strong convexity of the smoothed objective function that holds with high probability. The weak convergence is based on the non-asymptotic uniform Bahadur representation (Theorem~\ref{thm:bahadur}),  complemented by showing the convergence of finite-dimensional marginals and the asymptotic tightness.

\subsection{Rademacher multiplier bootstrap inference}
\label{sec:fst.boot}

In this section, we establish the theoretical guarantees of the Rademacher multiplier/weighted bootstrap for censored quantile regression as described in Section~\ref{sec:mb}.  In this case, $W_i = e_i+1$ and $e_i$'s are i.i.d. Rademacher random variables. 
For the random covaviate vector $\bx\in \RR^p$,  we assume that the moment parameters $m_3$ and $m_4$ defined in \eqref{def:m34} are dimension-free.
We first present the (conditional) uniform consistency of the bootstrapped process $\{ \hat{\bbeta}^\flat(\tau) : \tau \in [\tau_L, \tau_U]\}$ given the observed data $\mathbb{D}_n= \{ (y_i, \Delta_i, \bx_i \}_{i=1}^n$. Let $\PP^*(\cdot) = \PP(\cdot\, | \mathbb{D}_n)$ be the conditional probability given $\mathbb{D}_n$.

\begin{theorem}[Conditional uniform consistency]
\label{thm:concen.boot}
Assume Conditions~\ref{cond:kernel}--\ref{cond:grid} hold, and let the bandwidth satisfy $h=h_n  \asymp \{ (p+\log n) / n \}^\gamma$ for some $\gamma \in [1/4, 1/2)$.
Then, there exists an event $\cE = \cE(\mathbb{D}_n)$ with $\PP(\cE) \ge 1 - C_3 n^{-1}$ such that conditional on $\cE$, the bound \eqref{eq:thm,concen} holds, and the bootstrapped process $\hat{\bbeta}^{\flat}(\cdot):[\tau_L, \tau_U] \mapsto \RR^p$ satisfies
\#
\sup_{\tau \in [\tau_L, \tau_U] } ||\hat{\bbeta}^\flat(\tau) - \hat \bbeta(\tau) ||_{\Sigma} \lesssim \sqrt{\frac{p + \log n}{n}},    \label{eq:thm.concen.boot}
\#
with $\PP^*$-probability at least $1 - C_3 n^{-1}$, provided $\zeta_p^2 (p + \log n)^{1/2-\gamma} (p \log n)^{1/2} \lesssim n^{1-\gamma}$. Here $C_3 > 0$ is an absolute constant.
\end{theorem}

Analogously to \eqref{def:accu.error},  define the bootstrapped integrated error as
\#
\hat{\be}^{\,\flat}(\tau)   :=  \Jb(\tau) \{ \hat{\bbeta}^\flat(\tau) - \hat \bbeta(\tau) \} + \int_{\tau_L}^{\tau} \Hb(u) \{ \hat{\bbeta}^\flat(u) - \hat \bbeta(u) \}\, \mathrm{d} H(u) , \label{accu.error.boot}
\#
where $\Jb(\cdot)$ and $ \Hb(\cdot)$ are given in \eqref{def.JH}.
We then develop a linear representation for $\hat{\be}^{\,\flat}(\tau) $, which can be viewed as a parallel version of Theorem~\ref{thm:bahadur} in the bootstrap world.

\begin{theorem}[Conditional uniform Bahadur representation]
\label{thm:bahadur.boot}
Assume the conditions in Theorem~\ref{thm:concen.boot} hold, and that the kernel $K(\cdot)$ in Condition~\ref{cond:kernel} is Lipschitz continuous. Moreover, assume $\delta^* \lesssim n^{-(1/2 + \alpha)}$ for some $\alpha >0$.
Then, there exists an event $\cF = \cF(\mathbb{D}_n)$ with $\PP(\cF) \ge 1 - C_4 n^{-1}$ such that conditional on $\cF$,  \eqref{eq:thm,bahadur}--\eqref{eq:prob,bahadur} hold, and the bootstrapped process $\hat{\bbeta}^\flat(\cdot) :[\tau_L, \tau_U] \mapsto \RR^p$ satisfies
\#
\hat{\be}^{\,\flat}(\tau)   = \frac{1}{n} \sn \bU^\flat_i(\tau) + \br^\flat_n(\tau),   \label{eq:thm,bahadur.boot}
\#
where $\bU_i^\flat(\tau) =  e_i \bU_i(\tau)$ with $\bU_i(\tau)$ defined in \eqref{bahadur.ui}, and 
\#
& \sup_{\tau \in [\tau_L, \tau_U]}  \| \br^\flat_n(\tau) {\|}_{\Sigma^{-1}} \nn \\
 & \lesssim m_4^{1/2} \frac{p + \log n}{ nh^{1/2}} +  h \sqrt{\frac{p + \log n}{n}} + \zeta_p^2 \frac{(p + \log n) (p \log n)^{1/2}}{n^{3/2}h} +   n^{-1/2 - \alpha}    \label{eq:prob,bahadur.boot}
\#
with $\PP^*$-probability at least $1- C_4 n^{-1}$. 
\end{theorem}

Theorem~\ref{thm:bahadur.boot} shows that the bootstrap integrated error $\hat{\be}^{\,\flat}(\cdot) $  can be approximated,   up to a higher order remainder, by the linear process $\{ (1/n) \sn e_i \bU_i(\tau) : \tau\in[\tau_L, \tau_U]\}$, where $e_i$'s are independent Rademacher random variables, and $\EE^* \bU_i^\flat(\tau) = \bm{0}$. 
Provided that $h \asymp \{(p+\log n)/n\}^{2/5}$ and $p$ satisfies the growth condition $p = o(n^{3/8})$ as in Theorem~\ref{thm:weak.convergence}, then applying the same analysis in Remark~\ref{remark:thm3.2} gives us the following stochastic integral representation: with probability (over $\mathbb{D}_n$) approaching one, $\sup_{\tau \in [\tau_L, \tau_U]}  \| \br^\flat_n(\tau) {\|}_{\Sigma^{-1}} = o_{\PP^*}(1)$, and
\#
n^{1/2}  \Jb(\tau) \{ \hat \bbeta^\flat(\tau) -  \hat \bbeta (\tau) \}   = \bphi\Bigg( \frac{1}{n^{1/2}} \sn  \bU^\flat_i(\tau)   \Bigg)  + o_{\PP^*}(1) , \label{boot.sto.inte}
\#
where $\bphi$ is the linear operator defined in \eqref{def.linear.ope}.
Note that $\EE^* \{ \bU_i^\flat(s)  \bU_i^\flat(t)^\T \} = \bU_i(s) \bU_i(t)^\T$ for any $s, t \in [\tau_L, \tau_U]$. 
It can be shown that on $[\tau_L, \tau_U]$,  $n^{-1/2}\sn   \{  \bU_i(\cdot)- \EE \bU_i(\cdot)\}$ has the same asymptotic distribution as $n^{-1/2} \sn \bU^\flat_i(\cdot)$ conditionally on the data $\mathbb{D}_n$; see Theorem~\ref{thm:weak.convergence} and Theorem~\ref{thm:weak.convergence.boot} below.
This, together with \eqref{sde.solution} and \eqref{boot.sto.inte},  validates to some level the use of the bootstrap process $\hat \bbeta^\flat(\cdot)$ in the inference.
To illustrate this, consider the following bootstrap counterpart of the process $\GG_n(\cdot)$ defined in \eqref{def.Gn}:
\#
\GG^\flat_n(\tau) :=  \frac{1}{ n^{1/2} } \sn \langle  \ba_n / || \ba_n||_\Sigma ,  \bU^\flat_i(\tau)\rangle ,   \quad \tau \in [\tau_L, \tau_U].  \label{def.Gn.boot}
\#

\begin{theorem}[Validation of bootstrap process] 
\label{thm:weak.convergence.boot}
Assume Conditions~\ref{cond:kernel}--\ref{cond:grid} hold with $\delta^* \lesssim n^{-(1/2 + \alpha)}$ for $\alpha \in (0, 1/2)$, $h \asymp \{(p+\log n)/n\}^{2/5}$ and $p = o(n^{3/8})$.
In addition, assume the kernel $K(\cdot)$ is Lipschitz continuous. Then, for any sequence of (deterministic) vectors $\{\ba_n\}_{n = 1}^\infty$, there exists a sequence of events $\{\cF_n = \cF_n(\mathbb{D}_n)\}_{n = 1}^\infty$ such that $\PP(\cF_n) \to 1$, and conditional on $\{\cF_n\}_{n = 1}^\infty$, \eqref{eq:thm,bahadur.boot} holds and the conditional distribution of $\GG^\flat_n(\cdot)$ given $\mathbb{D}_n$ is asymptotically equivalent to the unconditional distribution of $\GG_n(\cdot)$ established in \eqref{eq:weak.convergence}.
\end{theorem}

\section{Regularized Censored Quantile Regression}
\label{sec:highdim}
We extend the proposed SEE approach to high-dimensional sparse QR models with random censoring.
The goal is to identify the set of relevant predictors, defined as
\#
	\cS^* = \bigcup_{\tau \in [\tau_L, \tau_U]} \supp\big( \bbeta^*(\tau) \big),  \label{true.set}
\#
assuming that its cardinality $s:=|\cS^*|$ is much smaller than the ambient dimension $p$---the total number of predictors, but may grow with sample size $n$.
Recall the sequentially defined smoothed loss functions $\hat L_k(\cdot)$ $(k=0,1,\ldots, m)$  in \eqref{def:L0} and \eqref{def:Lk}. When $p< n$, finding the solution to the SEE $\hat Q_k(\bbeta)=0$ is equivalent to solving the optimization problem $\min_{\bbeta\in \RR^p} \hat L_k(\bbeta)$. For fitting sparse models in high dimensions, we start with the $\ell_1$-penalized approach \citep{Tibs1996,BC2011}. At quantile levels $\tau_L=\tau_0   < \tau_1 <  \cdots <\tau_m=\tau_U$, we define  $\ell_1$-penalized smoothed CQR estimators $\hat \bbeta_k := \hat \bbeta(\tau_k)$ sequentially as
\#
\hat \bbeta(\tau_k)  \in \argmin_{\bbeta \in \RR^p}   \big\{   \, \hat L_k(\bbeta) + \lambda_k \cdot \| \bbeta \|_1  \big\}, \label{def:l1cqr}
\#
for $k \in \{0, \dots, m\}$,  where $0<\lambda_1 \leq \cdots \leq \lambda_m$ are regularization parameters. Define $\hat \bbeta(\tau) = \hat \bbeta(\tau_{k - 1})$ for $\tau \in (\tau_{k-1} , \tau_k)$.   It is worth noticing that for each $k\geq 1$, $\hat L_k(\cdot)$ is essentially a shifted or perturbed version of $\hat L_0(\cdot)$, that is,
$ \hat L_k(\bbeta) = \hat L_0(\bbeta)  - (1/n) \sn   \sum_{j = 0}^{k - 1} \bar K_h   (  y_i -\bx_i^\T \hat{\bbeta}_j )  \{ H(\tau_{j+1} )- H(\tau_j) \}  \bx_i^\T   \bbeta   $, where $H(u) = -\log(1-u)$.  All these empirical loss functions are convex, and have the same Hessian matrix.

 \begin{cond}[Random design in high dimensions] 
 \label{cond:hd.covariates}
 The (random)   covariate vector $\bx = (x_1, \ldots, x_p)^\T \in \cX \subseteq  \RR^p$ ($x_1 \equiv 1$) is compactly supported with $\max_{1\leq j\leq p} |x_j | \leq C_0$ almost surely for some $C_0\geq 1$. For convenience,  assume $C_0=1$.  The normalized vector $\Sigma^{-1/2} \bx$ has uniformly bounded kurtosis, that is, $m_4$ defined in \eqref{def:m34} is a dimension-free constant, where $\Sigma = \EE(\bx \bx^\T)$ is positive definite.
\end{cond}

\begin{theorem} 
\label{thm:hd}
Assume Conditions~\ref{cond:kernel}, \ref{cond:density}, \ref{cond:grid} and Condition~\ref{cond:hd.covariates} hold.  Under the sample size scaling $n\gtrsim s^3 \log p$,  let the bandwidth $h $ and penalty levels $\lambda_k$'s satisfy $s\sqrt{\log(p)/n } \lesssim h \lesssim \{ s \log(p)/n \}^{1/4}$ and  $\lambda_k \asymp   \{ 1+  \log(\tfrac{1-\tau_L}{1-\tau_k}) \}  \sqrt{\log(p)/n}$ for $k=0, 1,\ldots, m$. 
Then, there exist constants $C_0, C_1,C_2>0$ independent of $(s, p, n)$ such that 
\$
 \sup_{\tau_L \leq \tau \leq \tau_U } \| \hat \bbeta(\tau) - \bbeta^*(\tau) \|_\Sigma  \leq C_1  \bigg( \frac{1-\tau_L}{1-\tau_U} \bigg)^{C_0 \overline f / \underline{g}}\underline{g}^{-1}   \log\big( \tfrac{1-\tau_L}{1-\tau_U} \big) \sqrt{\frac{s \log p}{\underline{\gamma} n }}
\$
with probability at least $1-C_2 p^{-1}$, where $\underline{\gamma} = \lambda_{\min}(\Sigma)$ is the minimal eigenvalue of $\Sigma$.
\end{theorem}

Theorem~\ref{thm:hd} provides the rate of convergence for the $\ell_1$-penalized smoothed CQR estimator $\hat \bbeta(\cdot)$ uniformly in the set of quantile indices $\tau \in [\tau_L, \tau_U]$.
Under a similar set of assumptions,  \cite{ZPH2018} established the uniform convergence rate for the $\ell_1$-penalized (non-smoothed) CQR estimator, which is of order $\exp(C s) \sqrt{s \log (p \vee n) / n}$.
We conjecture that the additional exponential term $\exp(C s)$ is a consequence of the marginal smoothness condition posed in \cite{ZPH2018} (see Condition~(C4) therein),  and can be relaxed as in our Theorem~\ref{thm:hd}.
In fact,  our analysis relies on the global Lipschitz property \eqref{beta.lip}, which follows directly from the model assumption \eqref{cqr.model} and a lower bound on the conditional density.

\begin{remark}[Comments on the tuning parameters $h$ and $\{ \lambda_k \}_{k=0}^m$]
\label{rmk:highd.tuning}
To achieve the same convergence rate $\sqrt{s\log(p)/n}$ for the $\ell_1$-penalized QR estimator with non-censored data \citep{BC2011}, the bandwidth $h$ is required to be in the range specified in Theorem~\ref{thm:hd}; for example, one may choose $h\asymp \{ s \log (p)/ n \}^{1/4}$. Since such a choice depends on the unknown sparsity,  in practice we simply choose $h$ to be of order $\{\log(p)/n \}^{1/4}$.  Since the numerical performance is rather insensitive to the choice of bandwidth, we use the default value $h=\max\{ 0.05, 0.5\{ \log(p)/n \}^{1/4} \}$ as suggested in \cite{TWZ2022} although  it can also be tuned by cross-validation. 
 
The penalty levels $\lambda_k$'s play a more pivotal role in obtaining a reasonable fit for the whole CQR process.  Our theoretical analysis suggests that $\{ \lambda_k\}_{k=0}^m$ should be chosen as a slowly growing sequence along the $\tau$-grid.  Numerical results also confirm that a single $\lambda$ value, even after proper tuning,  cannot guarantee  a quality estimation of the entire regression process.  On the other hand,  it is computationally prohibitive to determine each $\lambda_k$ ($k=0, 1, \ldots, m$) via cross-validation.  By examining the proof of Theorem~\ref{thm:hd},  we see that once $\lambda_0$ is specified, the subsequent $\lambda_k$'s satisfy $\lambda_k = \{ 1 + \log (\frac{1-\tau_L}{1-\tau_k} ) \} \lambda_0$ for $k=1,\ldots, m$. Therefore, to implement the proposed sequential procedure, we only treat $\lambda_0$ as a tuning parameter, and use the above formula to determine the rest of $\lambda_k$'s.
\end{remark}

\begin{remark}[Adaptive $\ell_1$-penalization]
\label{nonconvexpenalty}
It has been recognized that the $\ell_1$-penalized estimator, with the penalty level determined via cross-validation,  typically has small prediction error but has a non-negligible estimation bias and tends to overfit with many false discoveries.  To reduce the estimation error and false positives,  a popular strategy is to use reweighted $\ell_1$-penalization via either adaptive Lasso \citep{Zou2006} or the local linear approximation (LLA) method for folded-concave penalties \citep{FL2001, ZL2008, TWZ2022}.   Let $w(\cdot)$ be a non-increasing and non-negative function defined on $[0, \infty)$.  Fix $k$,  let $\hat \bbeta^{(0)}_k = \hat \bbeta(\tau_k)$ be the $\ell_1$-penalized censored QR estimator at quantile level $\tau_k$.  For $t=1,  \ldots, T$, we iteratively update the previous estimate $\hat \bbeta^{(t-1)}_k$ by solving
\$
	\hat \bbeta^{(t)}_k = (\hat \beta^{(t)}_{k,1}  , \ldots, \hat \beta^{(t)}_{k,p})^\T  \in \argmin_{\bbeta \in \RR^p }  \Bigg\{  \hat L_k(\bbeta) +  \lambda_k \cdot \sum_{j=1}^p 	w\big( | \hat \beta_{k,j}^{(t -1)} | /\lambda_k  \big) | \beta_j |  \Bigg\} .
\$ 
When $T=1$ and $w(u) =u^{-1}$ for $u>0$ (or $(u+\epsilon)^{-1}$ for a small constant $\epsilon>0$) ,  this corresponds to an adaptive Lasso-type estimator \citep{Zou2006}; when $w(u) = \mathbbm{1}(u\leq 1) + \frac{(a-u)_+}{a-1} \mathbbm{1}(u>1)$ for $u\geq0$ and some $a>2$, this corresponds to the LLA method using the smoothly clipped absolute deviation (SCAD) penalty \citep{FL2001};  when $w(u) = (1-u/a)_+$ for $u\geq 0$ and some $a\geq 1$,  this corresponds to the LLA method using the minimax concave penalty (MCP) \citep{Z2010}.
\end{remark}

\section{Numerical Studies}
\label{sec:numerical}
We apply the proposed methods in Sections~\ref{sec:cqr} and \ref{sec:highdim} on simulated datasets and compare to that of \citet{PH2008} and \citet{ZPH2018} for both low- and high-dimensional settings in Sections~\ref{sec:numerical.est} and~\ref{sec:numerical.highd}, respectively.   
The proposed method involves selecting a smoothing parameter $h$: for $p<n$, we set $h = \{(p + \log n) / n\}^{2/5} \vee 0.05$; for $p>n$, guided by Remark~\ref{rmk:highd.tuning}, we set $h=\{0.05 \vee 0.5\{ \log(p)/n \}^{1/4}\}$.  We found that the performance of our proposed method is insensitive to the choice of bandwidth, as also observed in \cite{FGH2021} and \cite{HPTZ2020}.
We implemented \citet{PH2008} using the \texttt{crq} function with \texttt{method = "PengHuang"} from the  \texttt{quantreg} package  \citep{K2008}.   
On the other hand, \citet{ZPH2018} is implemented using the barebones function \texttt{LASSO.fit} from \texttt{rqPen} \citep{SM2020} instead of the function \texttt{rq(..., method = "lasso")} in the package \texttt{quantreg}.  This is because the function \texttt{rq(..., method = "lasso")} reports some numerical issues (e.g., singular design error) frequently in our numerical studies. 
All of the numerical studies are performed on a worker node with 32 CPUs, 2.5 GHz processor, and 512 GB of memory in a high-performance computing cluster.

\subsection{Censored quantile regression: estimation and inference}
\label{sec:numerical.est}
We assess the performance of our proposed method in the low-dimensional setting with $n=5000$ and $p=100$.  We start with generating the random covariates $\tilde\bx_i \in \RR^p$ from a mixture of different distributions to represent different types of variables commonly encountered in many datasets.  In particular, we generate the first 45 covariates from $\cN(\bm{0}, \Sigma = (\sigma_{jk})_{1\leq j, k \leq 45})$, where $\sigma_{jk} = 0.5^{|j - k|}$ for  $1 \le j, k \le 45$, the second 45 covariates from a multivariate uniform distribution on the  cube $[-2, 2]^{45}$ with the same covariance matrix $\Sigma$ using the \texttt{R} package \texttt{MultiRNG}, and the last 10 covariates from a Bernoulli distribution.  
Note that the three blocks of covariates generated are independent across the blocks.   
The response variables $z_i \in \RR$ are then generated from the following models, both of which satisfy the global assumption in \eqref{cqr.model}.

\begin{enumerate}
\item[(i)] Homoscedastic model: $z_i = \langle \tilde\bx_i, \bgamma \rangle + \varepsilon_i$ for $i = 1, \dots, n$, where $\gamma_j \sim \mathrm{Uniform}(-2, 2)$ for $j = 1, \dots, p$.  Let $Q_{t_2}(\tau)$ be the $\tau$-quantile of the $t_2$-distribution,  and let $\bx_i = (1, \tilde\bx_i^\T)^\T$. Then, the above model can be equivalently formulated as 
\#
z_i = \langle \bx_i,  \bbeta^*(\tau) \rangle, \ \ i= 1,\ldots, n , ~~\mbox{where}~~\bbeta^*(\tau) = (Q_{t_2}(\tau) , \bgamma^\T)^\T \in \RR^{p + 1} .  \label{model.homo}
\#
Under the above model, the covariate effects remain the same across all quantile levels.

\item[(ii)] Heteroscedastic model:  $z_i = \langle \wt \bx_i, \bgamma \rangle + |\tilde x_{i, 1}| \cdot \varepsilon_i$ for $i = 1, \dots, n$,  where $\gamma_1 = 0$ and $\gamma_j \sim \mathrm{Uniform}(-2, 2)$ for $j = 2, \dots, p$.
Let $\bx_i = (1, |\tilde x_{i,1}|,\wt\bx_{i,-1}^{\,\T})^\T$,  where $\wt\bx_{i,-1}\in \RR^{p-1}$ is  obtained by removing the first element of $\wt\bx_{i}$. The model is equivalent to 
\#
z_i = \langle \bx_i, \bbeta^*(\tau) \rangle, \ \ i= 1,\ldots, n , ~~\mbox{where}~~ \bbeta^*(\tau) = (0,  Q_{t_2}(\tau), \gamma_2, \dots, \gamma_p)^\T \in \RR^{p + 1} .  \label{model.hetero}
\#
In this model, the first covariate has varying marginal effects for different quantile levels.  Specifically,  the effect of $| \tilde x_1|$ on the $\tau$-th quantile of $z$ is $F^{-1}_{t_2}(\tau)$,  which is negligible when $\tau \approx 0.5$,  but grows stronger as $\tau$ moves towards 0 or 1.
\end{enumerate}
For both types of models, the random censoring variables are generated from a Gaussian mixture distribution, that is,
\#
C_i \sim \1\{w_i = 1\} \cN(0, 16) + \1\{w_i = 2\} \cN(5, 1) +  \1\{w_i = 3\} \cN(10, 0.25)  \label{model.censor}
\#
for $i = 1, \ldots, n$, where $w_i$ is sampled from $\{1, 2, 3\}$ with equal probability, and $y_i = z_i \wedge C_i$ is the censored outcome.  The corresponding censoring rate varies from $25\%$ to $50\%$.

\begin{figure}[!htp]
  \centering
  \subfigure[$\ell_2$-error under model \eqref{model.homo}]{\includegraphics[scale=0.32]{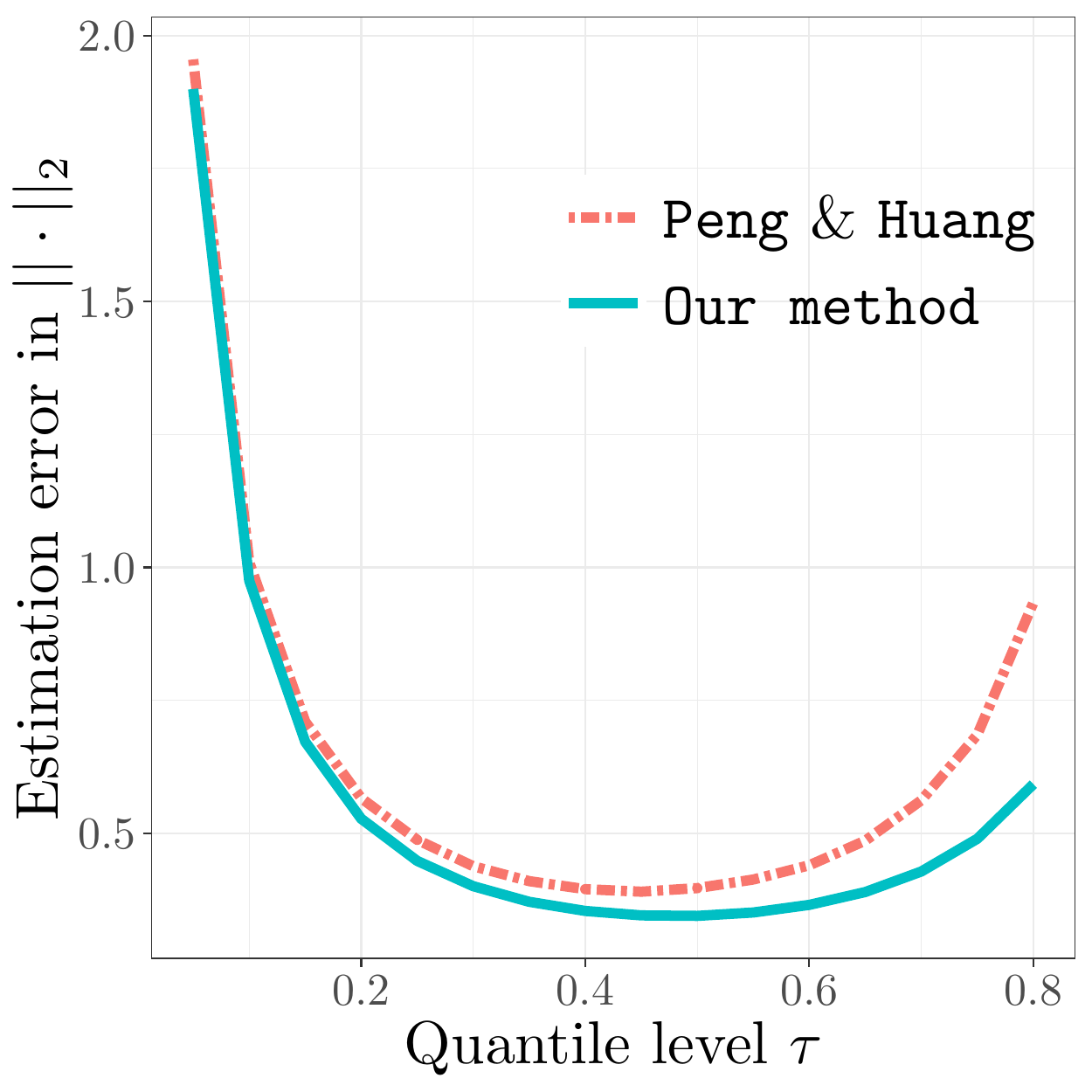}}  \qquad\quad 
  \subfigure[Estimated quantile effects under model \eqref{model.homo}]{\includegraphics[scale=0.32]{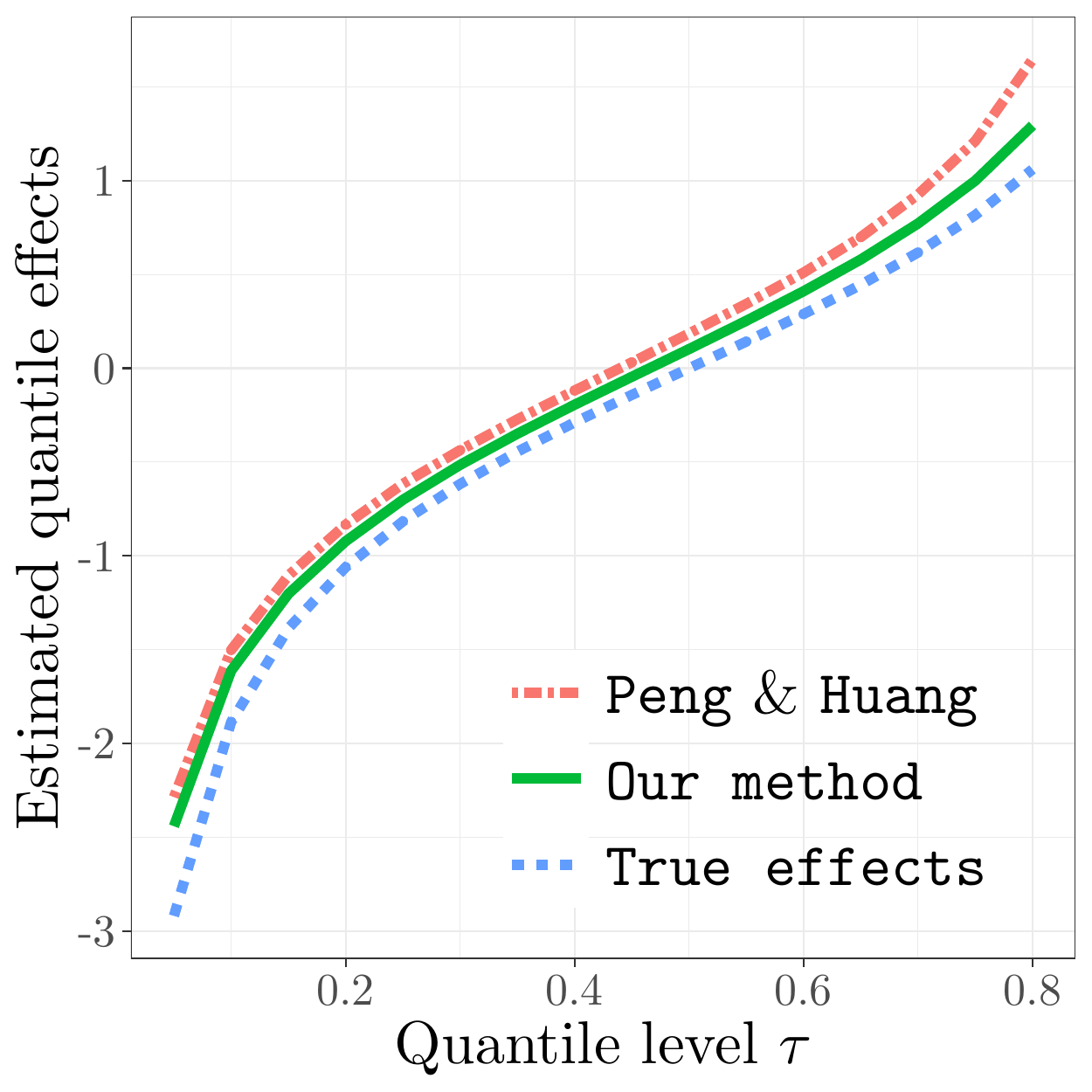}}  \qquad\quad 
  \subfigure[Runtime under model \eqref{model.homo}]{\includegraphics[scale=0.32]{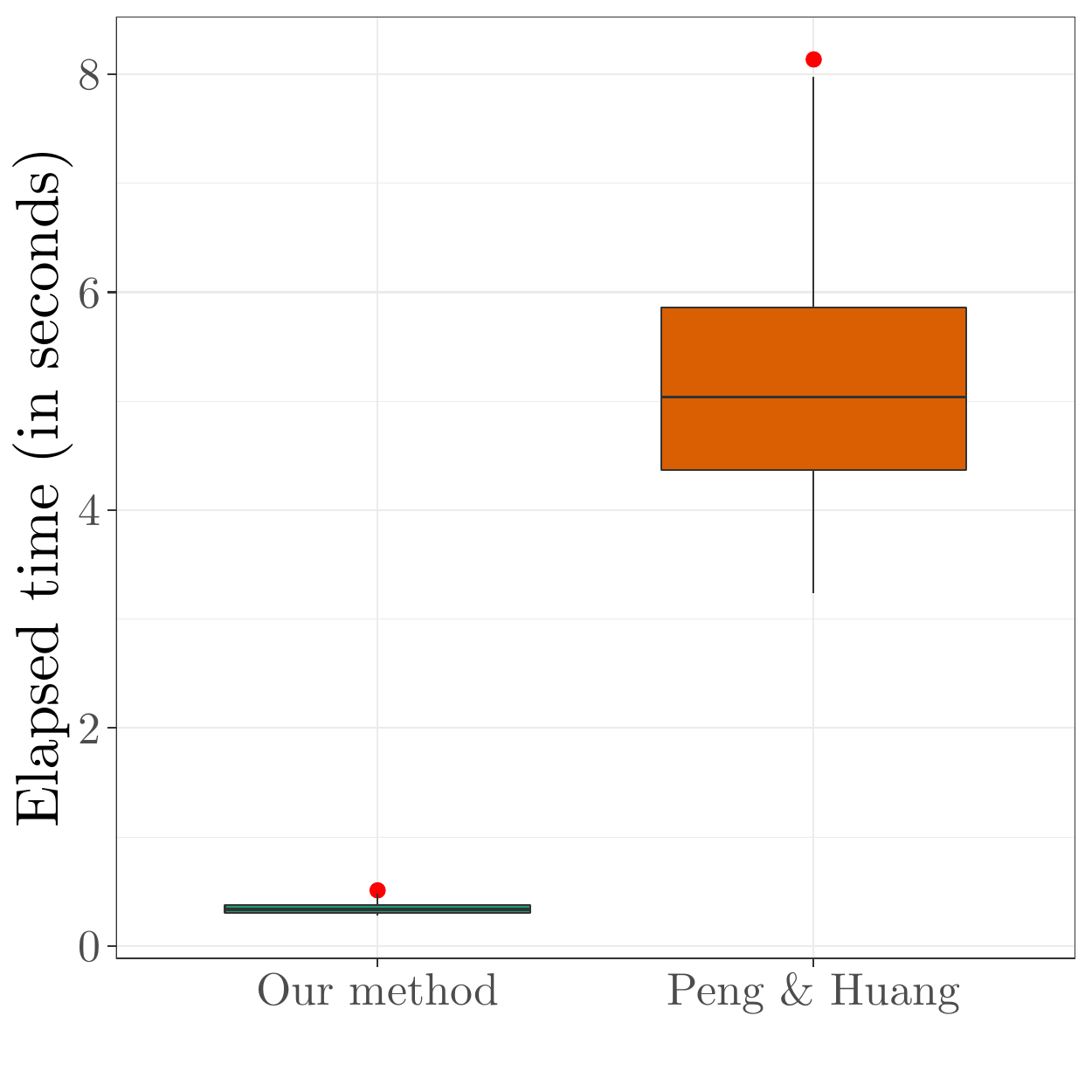}} 
  \subfigure[$\ell_2$-error under model \eqref{model.hetero}]{\includegraphics[scale=0.32]{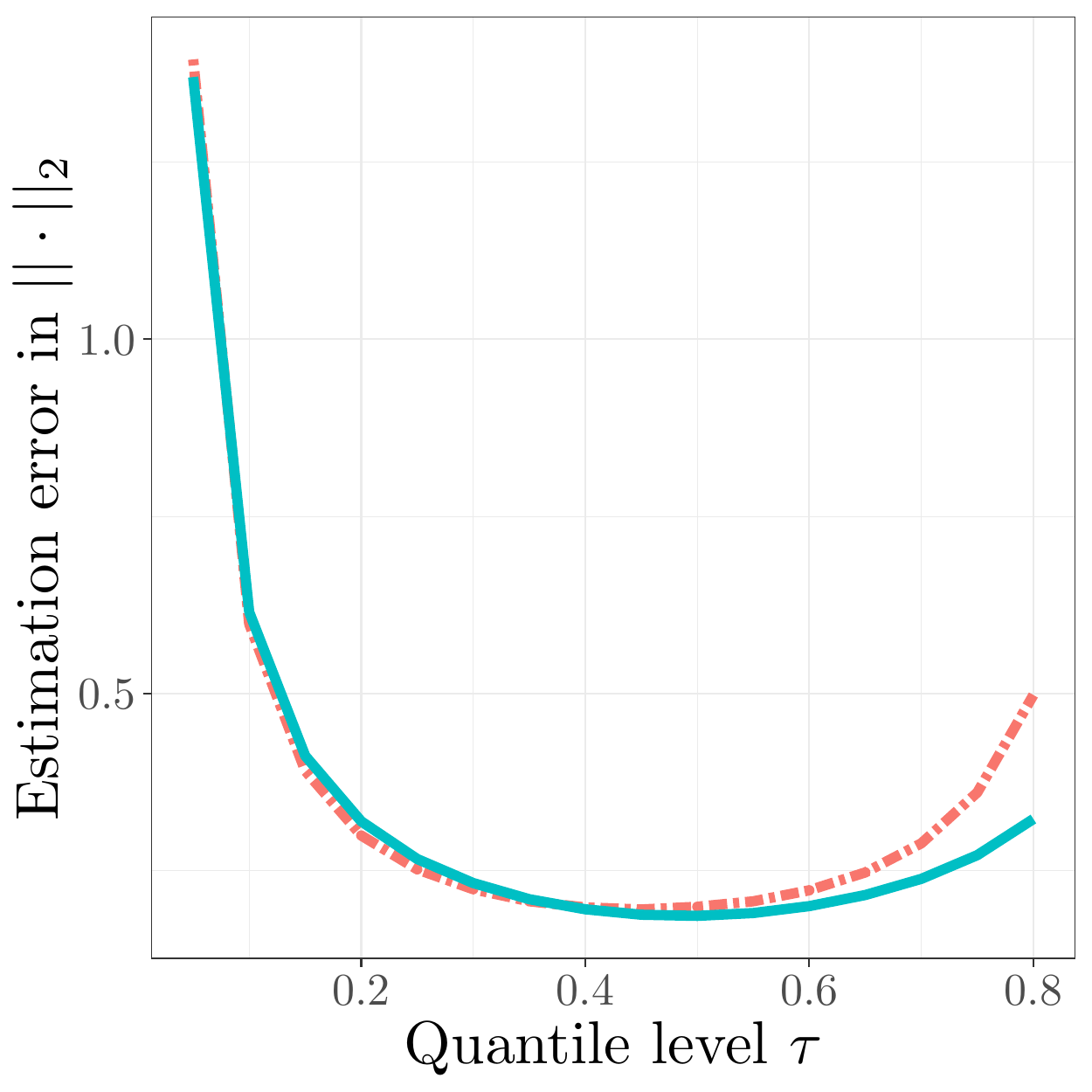}}  \qquad\quad 
  \subfigure[Estimated quantile effects under model \eqref{model.hetero}]{\includegraphics[scale=0.32]{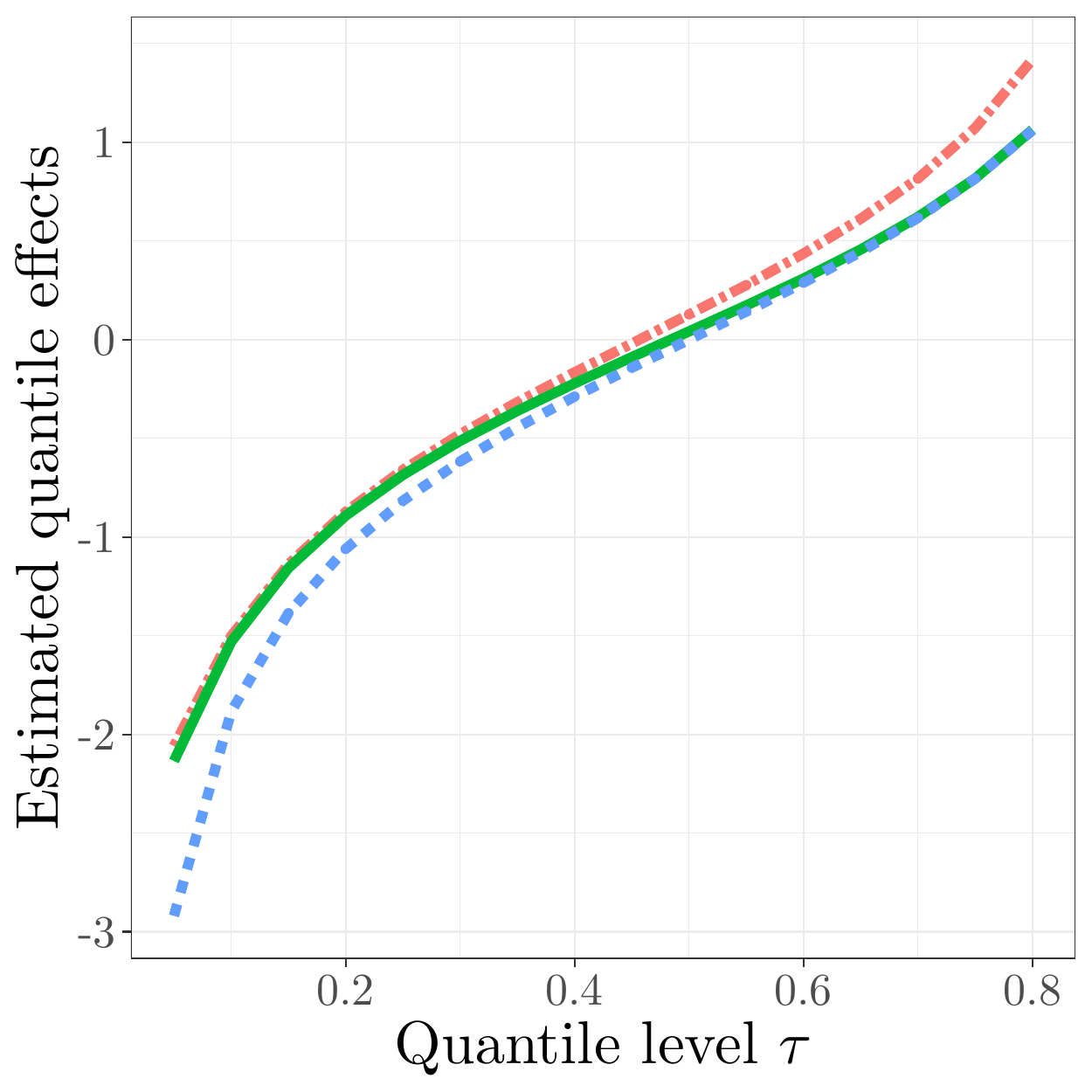}}  \qquad\quad 
  \subfigure[Runtime under model \eqref{model.hetero}]{\includegraphics[scale=0.32]{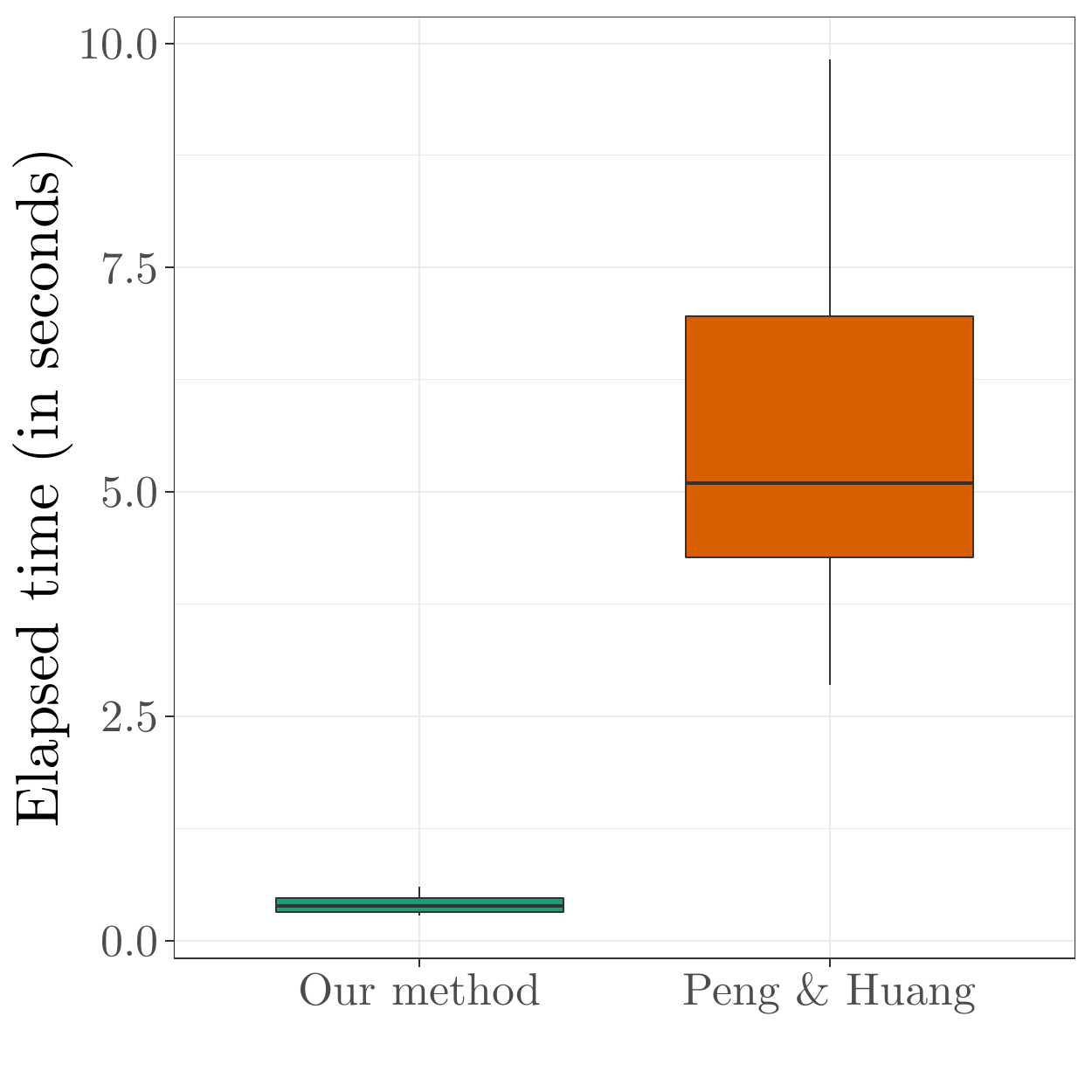}}  
\caption{Numerical comparisons among CQR and our smoothed CQR for models \eqref{model.homo}--\eqref{model.hetero} along the quantile grid. The left panels (a) and (d) display the $\ell_2$-induced estimation errors $\|\hat \bbeta(\tau_k) - \bbeta^*(\tau_k)\|_2$. The middle panels (b) and (e) present the estimated quantile effects, which are $\hat\beta_0(\tau_k)$ in model \eqref{model.homo} and $\hat\beta_1(\tau_k)$ in model \eqref{model.hetero} accordingly. The blue dashed lines in the middle panels represent the true quantile effects $Q_{t_2}(\tau)$. 
The right panels (c) and (f) record the empirical running time of the processes along the grid points.}
  \label{fig:process}
\end{figure}

We implement both methods with a quantile grid of $\{\tau_k\}_{k = 0}^m = \{0.05, 0.1, \dots, 0.75, 0.8\}$.  
At each quantile level $\tau_k$, we use the estimation error under the $\ell_2$ norm,  $\|\hat \bbeta(\tau_k) - \bbeta^*(\tau_k)\|_2$, as a general measure of accuracy.   We also calculate the run-time in seconds for both methods. Results, averaged across 500 independent replications, are reported in Figure~\ref{fig:process}.
Figures~\ref{fig:process}(a) and (d) contain the estimation error under the $\ell_2$ norm across all quantile levels;  Figures~\ref{fig:process}(b) and (d) contains the regression coefficient that varies across quantile levels, i.e., $\{ \beta_0(\tau_k) \}_{k=0}^m$ for  model \eqref{model.homo} and $\{ \beta_1(\tau_k) \}_{k=0}^m$ for model \eqref{model.hetero}; and  Figures~\ref{fig:process}(c) and (f) contain the computation time for fitting the entire QR process.   We see that the two methods perform very closely at low quantile levels, and the smoothed approach is particularly advantageous at high quantile levels. Computationally, our implementation of the smoothed method is about 10 to 20 times faster than \citet{PH2008}'s method, implemented by  the \texttt{crq} function in \texttt{quantreg}. The numerical results on smaller-scale datasets are presented in  Appendix~F.1 of the supplementary material.

\begin{figure}[!htp]
  \centering
  \subfigure[Coverage under model \eqref{model.homo}]{\includegraphics[scale=0.32]{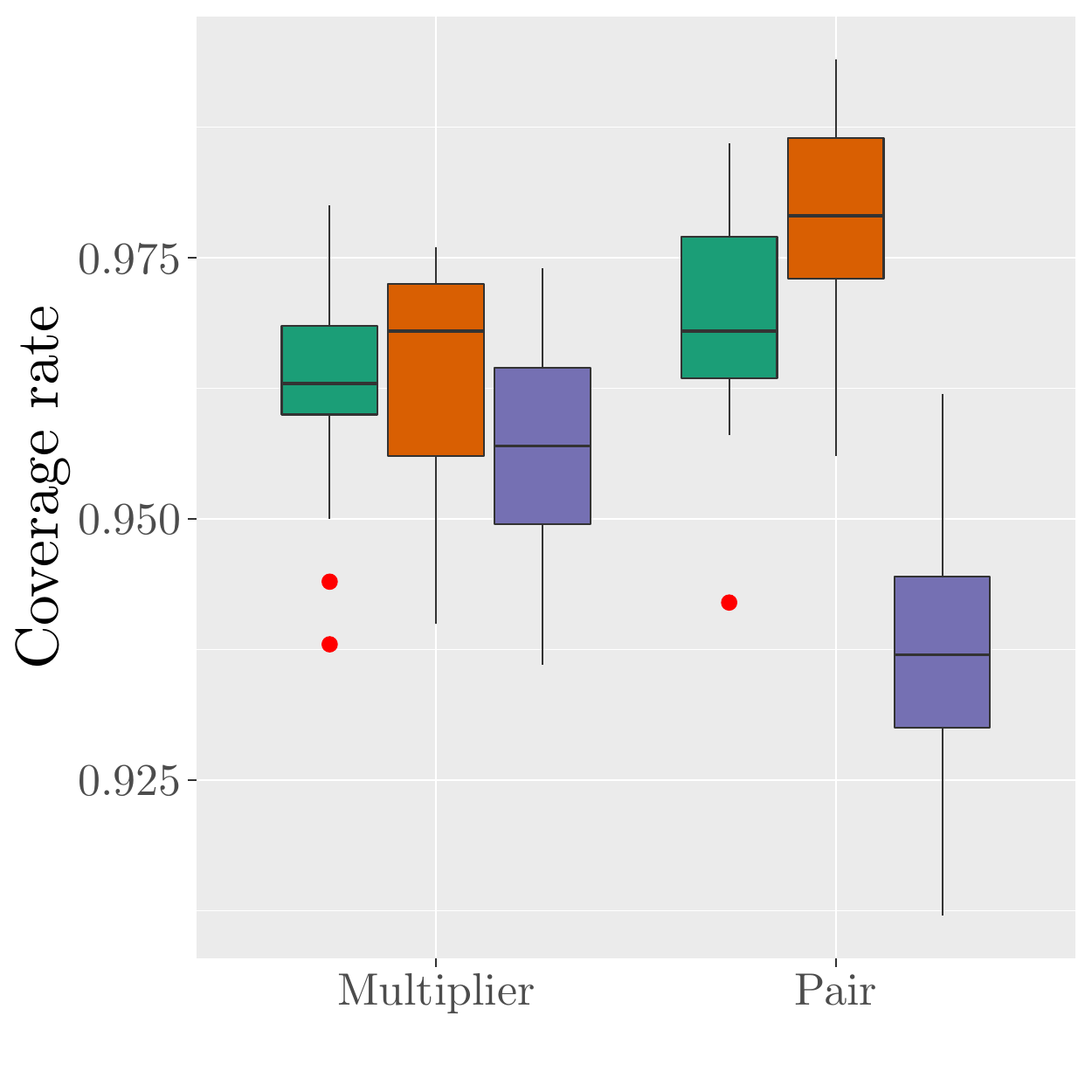}}  \qquad\quad 
  \subfigure[CI width under model \eqref{model.homo}]{\includegraphics[scale=0.32]{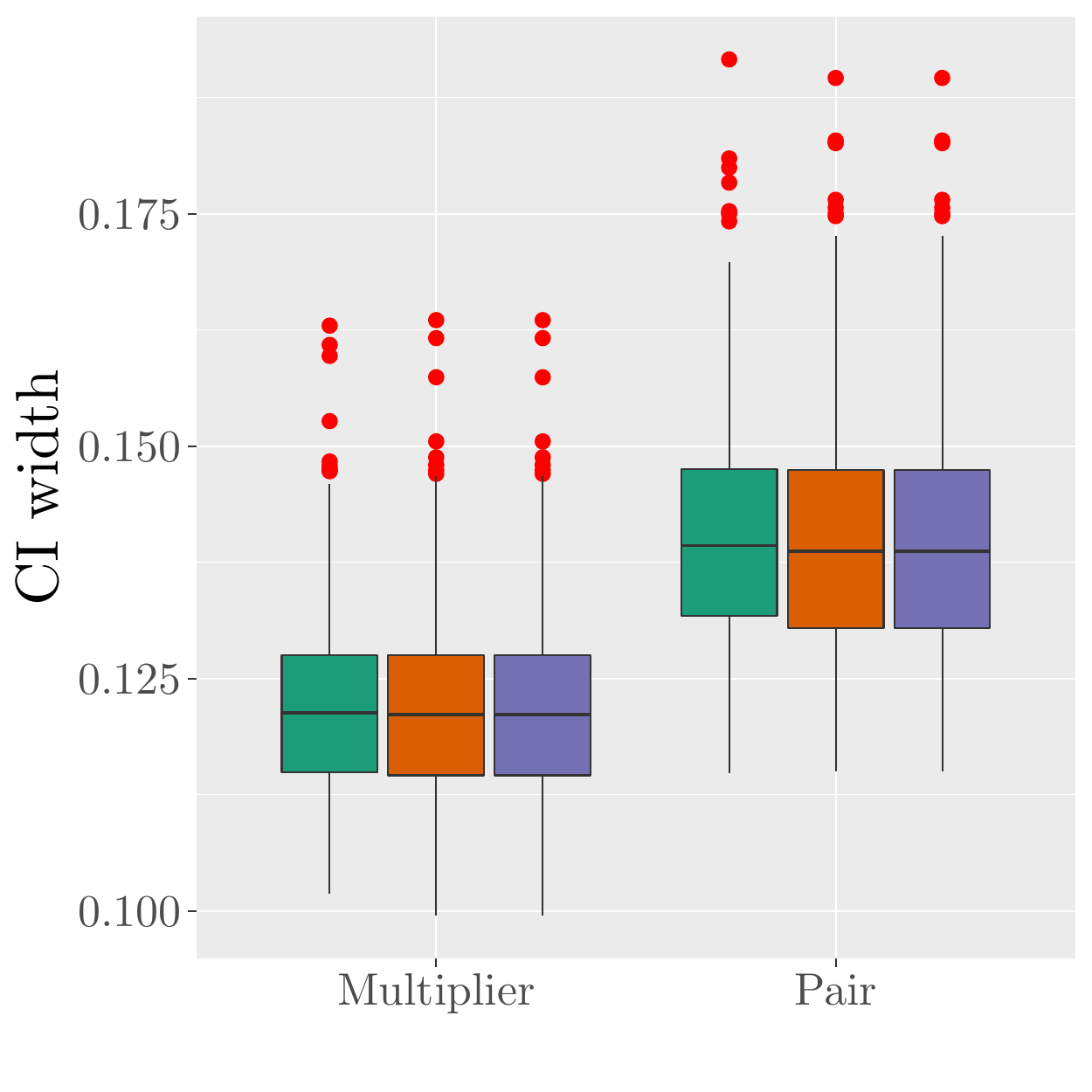}} \qquad\quad 
  \subfigure[Runtime under model \eqref{model.homo}]{\includegraphics[scale=0.32]{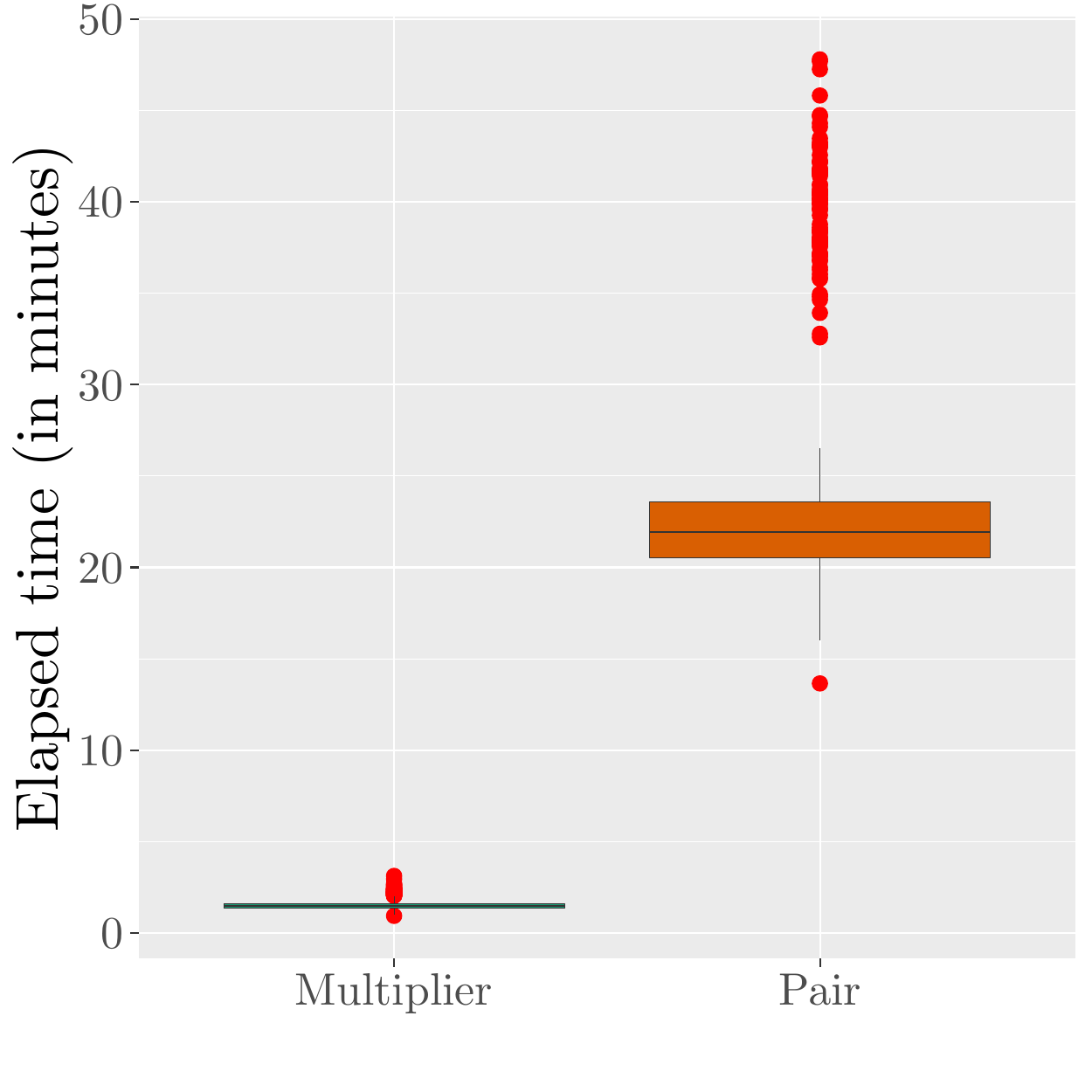}}
  \subfigure[Coverage under model \eqref{model.hetero}]{\includegraphics[scale=0.32]{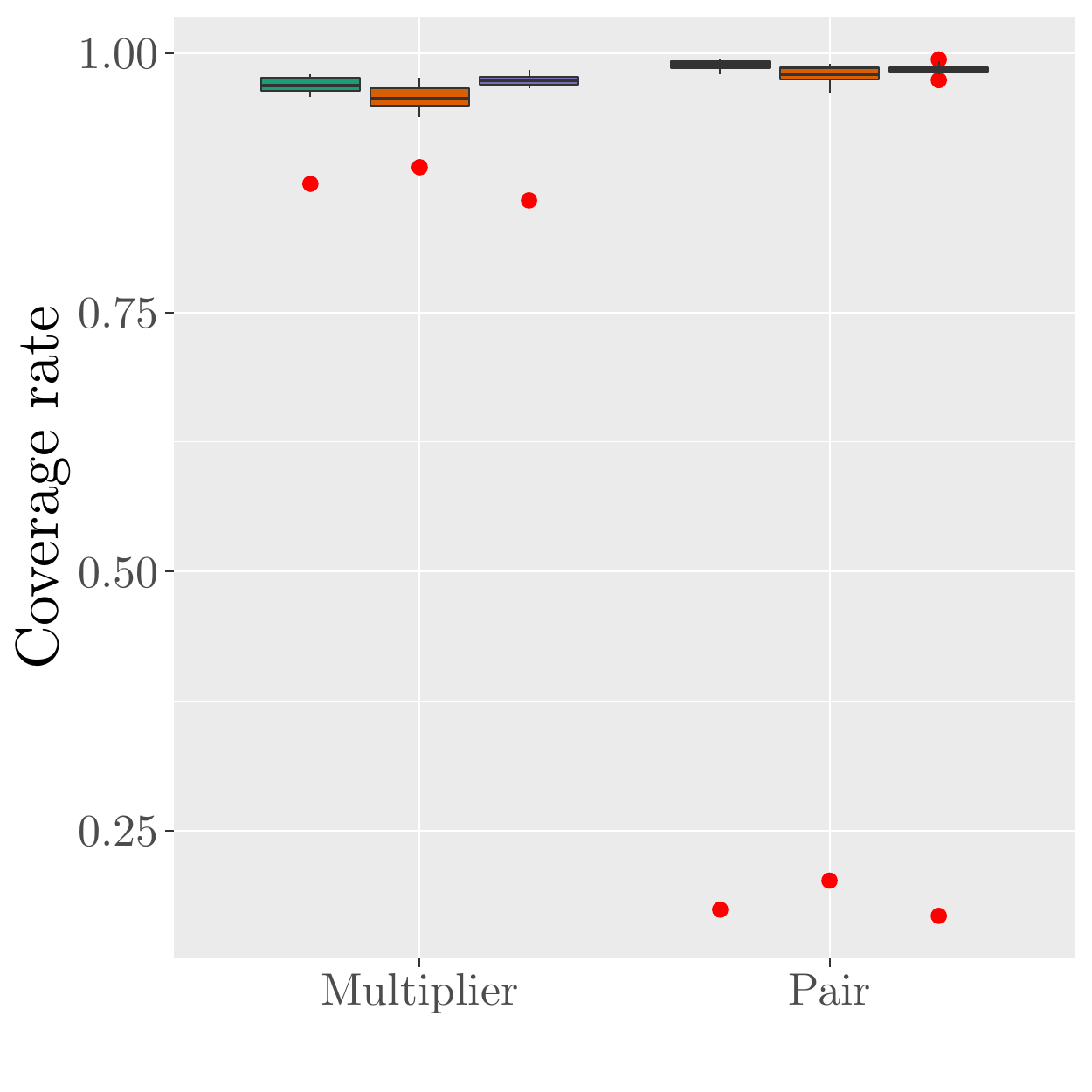}}  \qquad\quad 
  \subfigure[CI width under model \eqref{model.hetero}]{\includegraphics[scale=0.32]{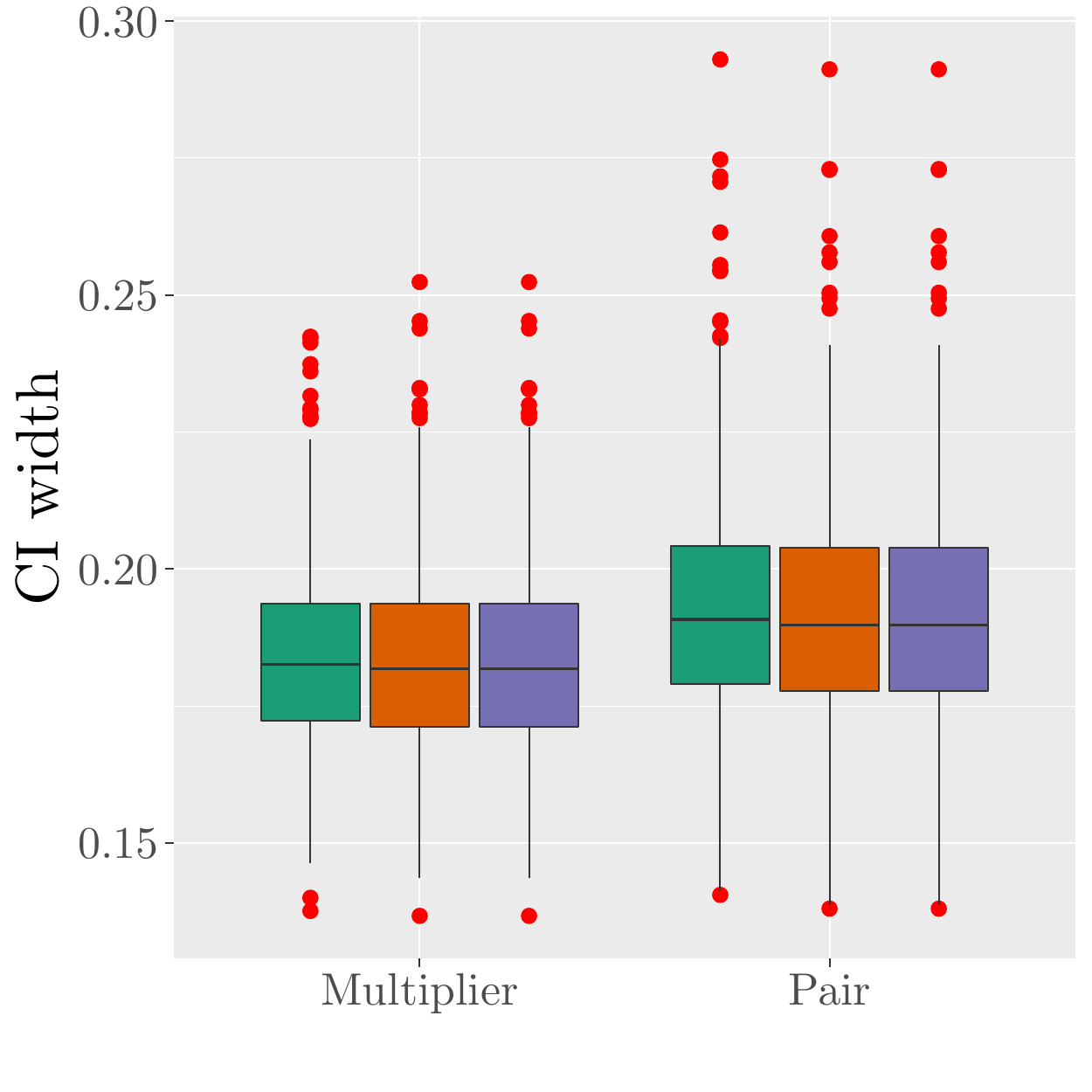}} \qquad\quad 
  \subfigure[Runtime under model \eqref{model.hetero}]{\includegraphics[scale=0.32]{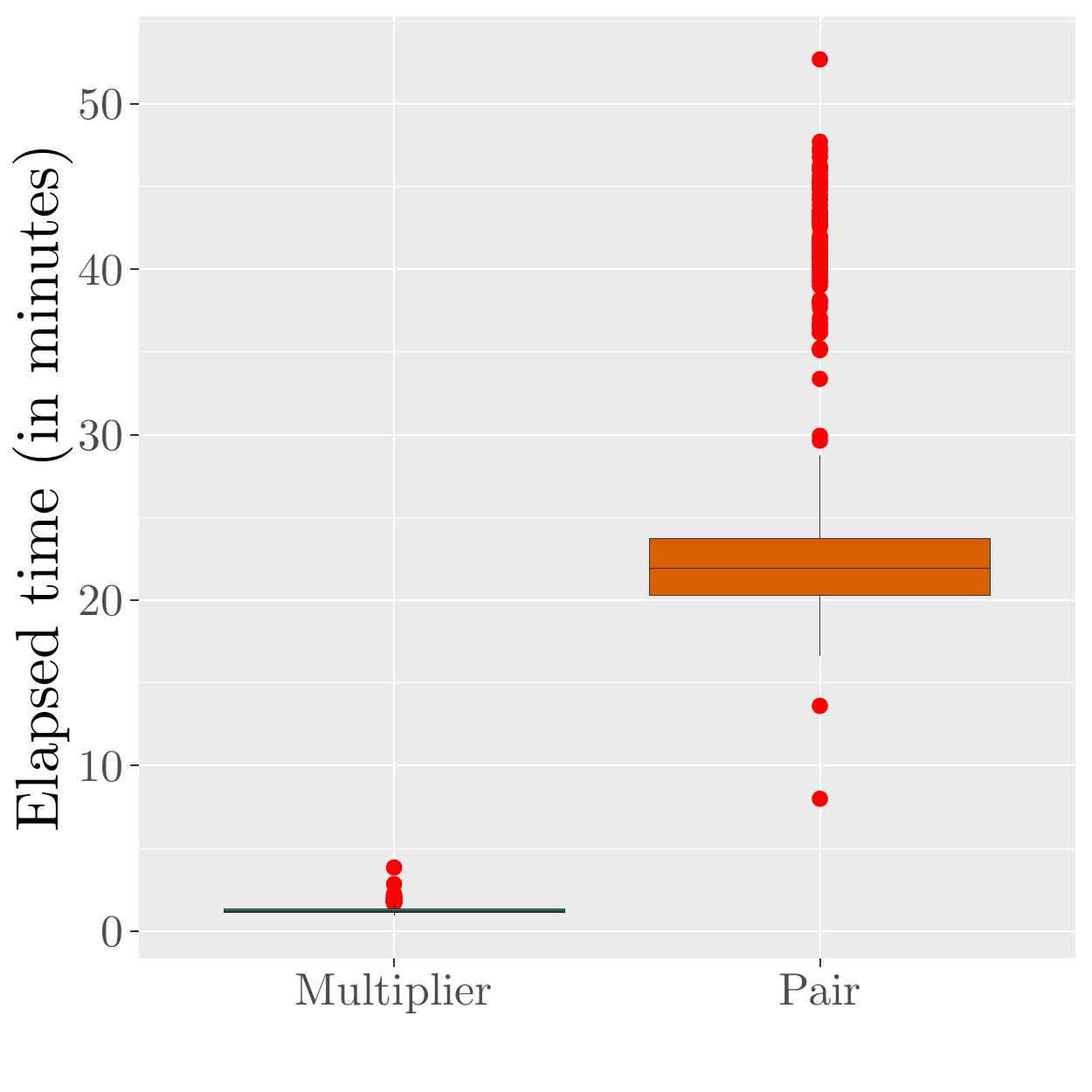}}
\caption{Box plots of the empirical coverage, confidence interval width, and running time for two resampling-based methods. ``Multiplier'' refers to the proposed multiplier bootstrap method, and ``Pair'' refers to pair resampling with replacement in the regression setting. In panels (a), (b), (d), and (e), within each method, different colors of boxes represent different types of confidence interval: (i) percentile interval \protect\includegraphics[height=0.85em]{Graphs/green.pdf}, (ii) pivotal interval \protect\includegraphics[height=0.85em]{Graphs/orange.pdf}, and (iii) normal interval \protect\includegraphics[height=0.85em]{Graphs/purple.pdf}.}
  \label{fig:simu.ci}
\end{figure}

Next, we consider both the proposed multiplier bootstrap detailed in Section~\ref{sec:mb} and the classical paired bootstrap for performing statistical inference at $\tau = 0.5$.
Three types of 95\% confidence intervals (CIs) are constructed with $B=1000$ bootstrap samples:  the percentile CI, the pivotal CI, and the normal CI.  Coverage proportions for all of the covariates, confidence interval width for the first covariate, and computational time for the entire bootstrap process, averaged over 500 replications, are plotted in Figure~\ref{fig:simu.ci}.
Under the homogeneous setting \eqref{model.homo}, all types of confidence intervals produced by multiplier bootstrap maintain the nominal level, while the normal intervals by pair resampling suffer from under coverage.
In the heterogeneous setting \eqref{model.hetero}, although outliers that correspond to the confidence intervals for the first covariate exist for both methods, multiplier bootstrap manages to mitigate this issue.
Furthermore, compared to pair resampling, multiplier bootstrap constructs narrower confidence intervals with slightly smaller standard deviations. Finally, the computational advantage of multiplier bootstrap for smoothed CQR is evident in Figures~\ref{fig:simu.ci}(c) and (f).

\begin{figure}[!htp]
  \centering
  \subfigure[{\scriptsize $\ell_2$-error under model \eqref{model.homo}}]{\includegraphics[scale=0.32]{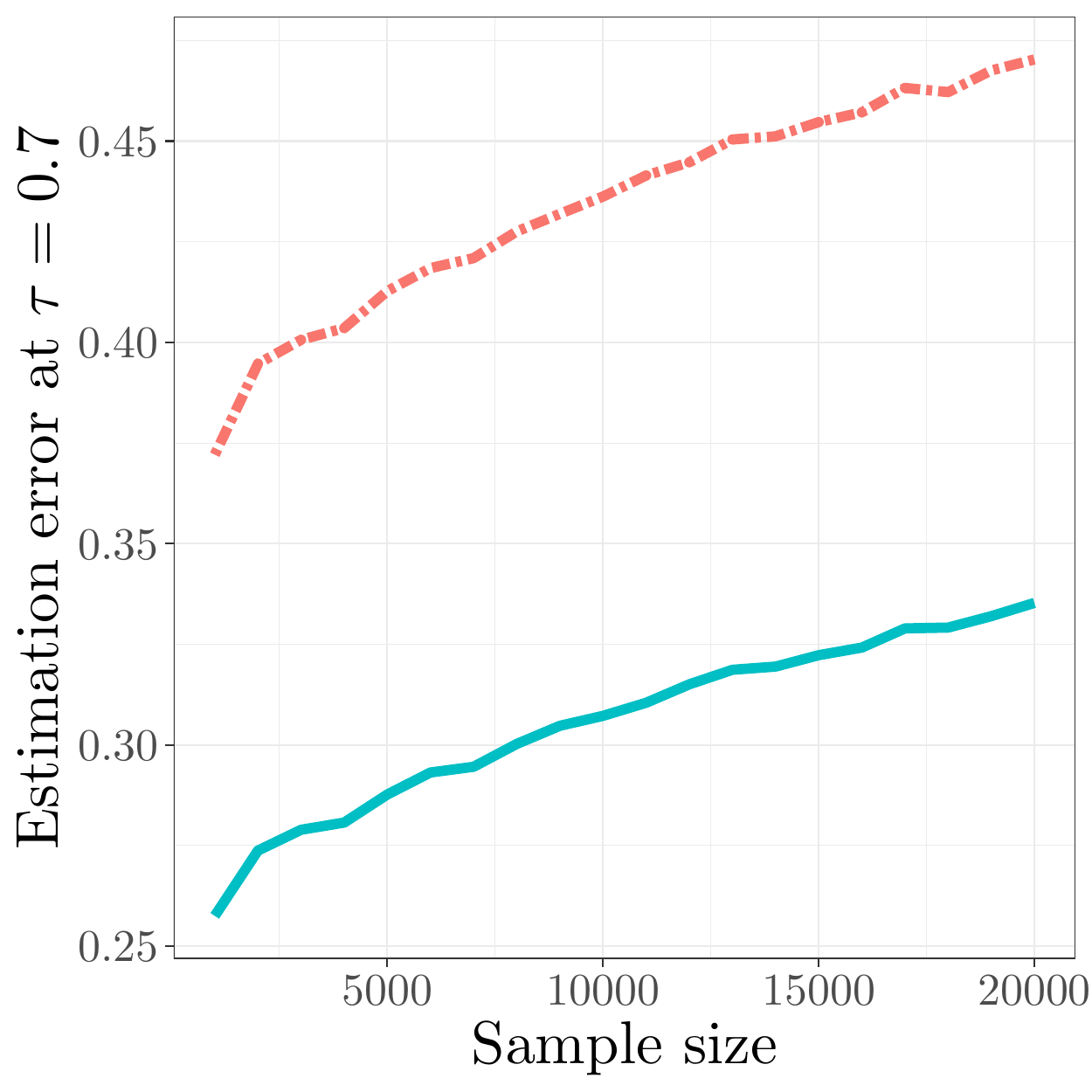}} \qquad\quad 
  \subfigure[{\scriptsize Runtime under model \eqref{model.homo}}]{\includegraphics[scale=0.32]{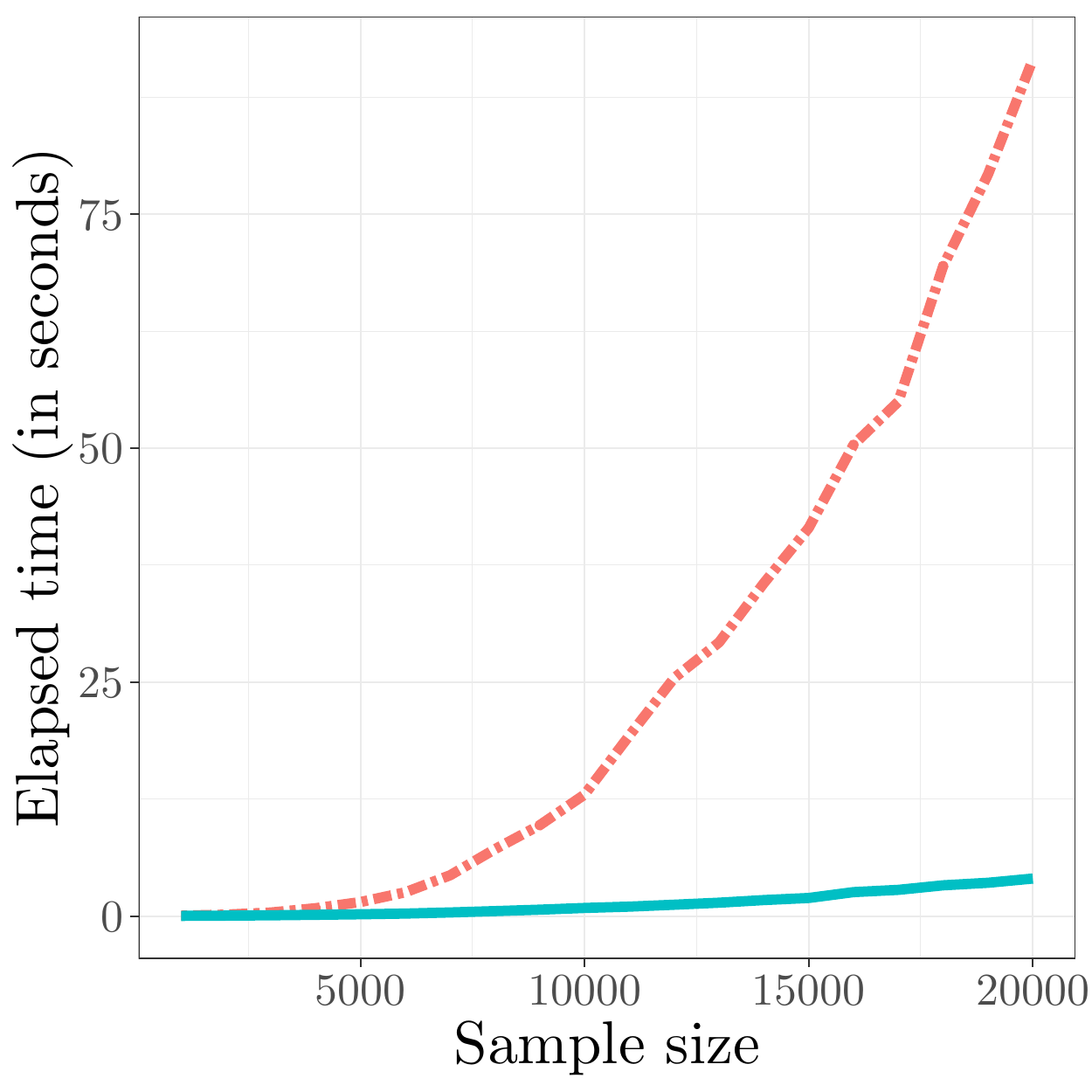}}  \\
  \subfigure[{\scriptsize $\ell_2$-error under model \eqref{model.hetero}}]{\includegraphics[scale=0.32]{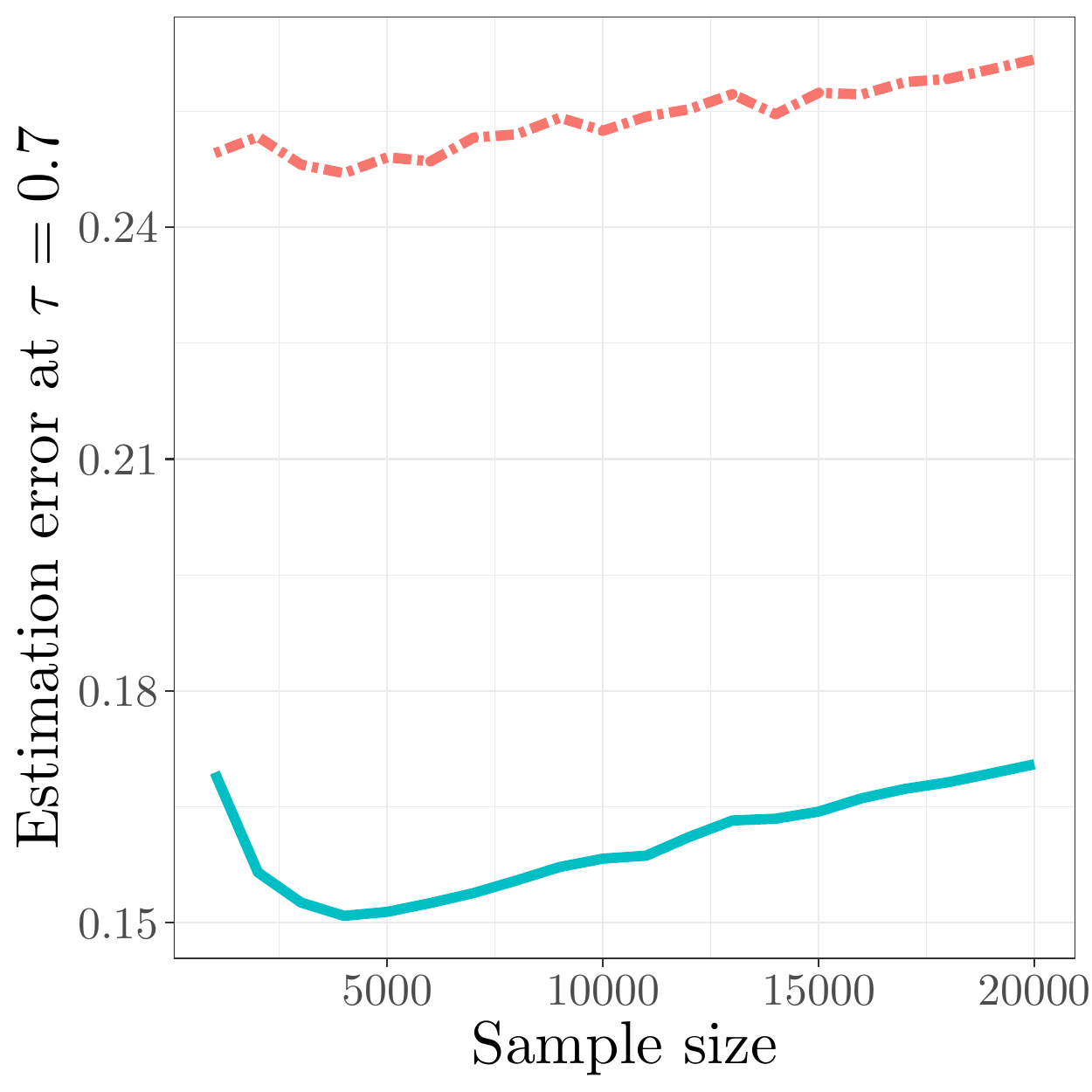}} \qquad\quad 
  \subfigure[{\scriptsize Runtime under model \eqref{model.hetero}}]{\includegraphics[scale=0.32]{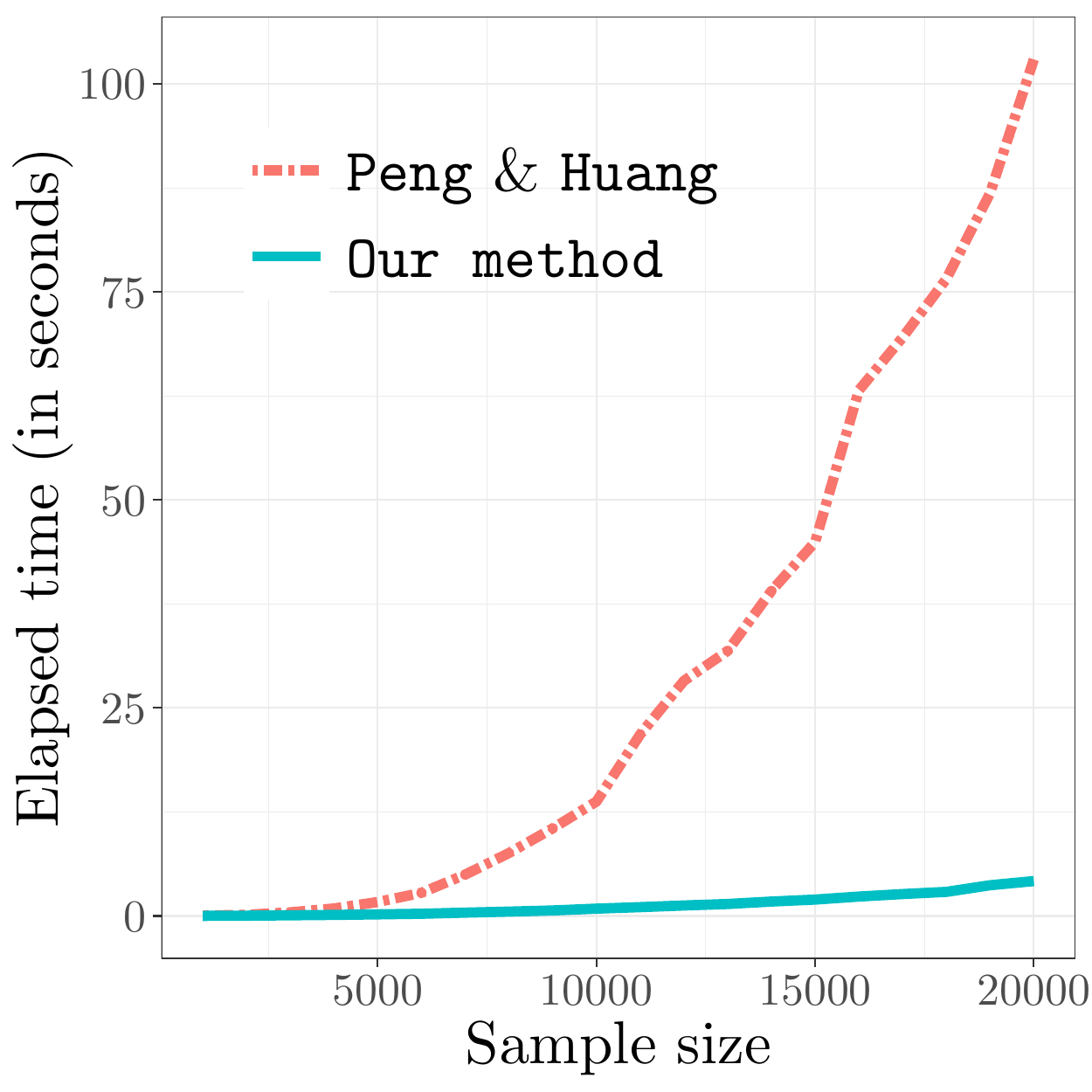}} 
\caption{Numerical comparisons between CQR and smoothed CQR under models \eqref{model.homo} and \eqref{model.hetero} with increasing $(n,p)$ subject to $p=n/100$. The left panels (a) and (c) display the $\ell_2$-error at $\tau = 0.7$ versus sample size. The right panels (b) and (d) present the runtime (in second) versus sample size.  
}
  \label{fig:growing}
\end{figure}

To better appreciate the computational advantage of smoothed CQR,   we further consider large-scale simulation settings by setting $n \in \{1000, 2000, \dots, 20000\}$ and $p = n/100$.  We use the same data generating processes as in \eqref{model.homo}--\eqref{model.censor}, except that the covariates $\tilde\bx_i$ are now generated from  $\cN(\bm{0}_p, \Sigma  )$ with $\Sigma = ( 0.5^{|j - k|})_{1\leq j, k \leq p}$.   The censoring rate varies from $30\%$ to $45\%$.  In this case, we restrict attention to the estimation error and runtime of the two methods when $\tau=0.7$.  The results,  averaged over $500$ repetitions, are presented in Figure~\ref{fig:growing}. We see from Figure~\ref{fig:growing} that the computation gain of the proposed method over \citet{PH2008} is dramatic, without compromising the statistical accuracy. The estimation errors at $\tau \in \{0.3, 0.5\}$, as functions of the sample size, are displayed in Figure F.3 in the supplementary material.

\subsection{High-dimensional censored quantile regression}
\label{sec:numerical.highd}
In this section, we examine the numerical performance of the regularized smoothed CQR method with different penalties, which will also be compared with its non-smoothed counterpart \citep{ZPH2018}. For the smoothed method,  we consider both the $\ell_1$ and folded-concave penalties (SCAD and MCP). The latter is implemented by the LLA algorithm as described in Remark~\ref{nonconvexpenalty}. The computational details are described in Section~A.2 of the supplementary material.

Penalized CQR involves selecting a sequence of regularization parameters $\{\lambda_k \}_{k=0}^m$ that correspond to the predetermined $\tau$-grid $\{ \tau_k \}_{k=0}^m$.  Guided by Theorem~\ref{thm:hd} and Remark~\ref{rmk:highd.tuning}, we adopt a sequence of dilating $\lambda_k$'s with $\lambda_k = \{ 1 + \log (\frac{1-\tau_L}{1-\tau_k} ) \} \lambda_0$ for $k=1,\ldots, m$, where $\lambda_0$ is chosen via the $K$-fold cross-validation ($K=3$ in our studies). 
To accommodate censoring, the cross-validation criterion is based on the the empirical mean of deviance residuals \citep{TGF1990} 
\#
R(\lambda) := \frac{1}{n}\frac{1}{m + 1}\sum_{i = 1}^n \sum_{k = 0}^m \sqrt{-2\{M_i(\tau_k, \lambda) + \Delta_i \log(\Delta_i - M_i(\tau_k, \lambda))\}}    \label{deviance}
\#
on the validation set, where 
\$
M_i(\tau_k, \lambda) = \mathbbm{1} \{y_i \le    \bx_i^\T \hat{\bbeta}(\tau_k, \lambda)  , \Delta_i = 1\} - \int_{\tau_0}^{\tau_k} \mathbbm{1} \{y_i \ge  \bx_i^\T \hat{\bbeta}(u, \lambda)   \}\,\mathrm{d} H(u) - \tau_0
\$
for $k = 0, \dots, m$ are the martingale residuals and $\hat{\bbeta}(\tau, \lambda)$ refers to the estimated $\bbeta(\tau)$ with a dilating $\lambda_k$'s starting with $\lambda_0 = \lambda$.
The deviance \eqref{deviance} produces a more symmetric distribution through a transformation on the skewed martingale residuals, and is also used in \cite{ZPH2018} and \cite{FZHL2021}.
In our simulations, we choose $\lambda_0$ from 50 candidates equally spaced on the interval $[0.01, 0.2]$.

In all of our numerical studies, we generate covariates $\tilde\bx_i \in \RR^p$ from $\cN(\bm{0}, \Sigma)$, where  $\Sigma$ is as defined in Section~\ref{sec:numerical.est}, and the random errors  $\varepsilon_i \sim t_2$. The response variables $z_i$ are generated from Models  \eqref{model.homo}--\eqref{model.hetero}, but with different $\bgamma$.
For Model~\eqref{model.homo}, we consider a sparse $\bgamma$ with global sparsity $s=10$ by setting $\gamma_j \sim  \mathrm{Uniform}(1, 1.5)$ for $j = 1, \dots, 10$, and the rest to be zero.  For Model~\eqref{model.hetero}, $\bgamma$ is generated similarly except with  $\gamma_1 = 0$.  
The random censoring variables are generated from \eqref{model.censor}, with overall censoring rates approximately $25\%$--$30\%$.

Since the estimated active set depends on the entire quantile process, all  numerical experiments are conducted via an estimation-after-selection procedure \citep{ZPH2018}.
That is, in stage one, we perform regularized smoothed CQR to obtain the  set $\hat{\cS} = \cup_{\tau \in \{\tau_0, \dots, \tau_m\}} \supp\, \big( \hat{\bbeta}(\tau) \big)$. In stage two, we perform smoothed CQR using the covariates in $\hat{\cS}$.
Recall that $\cS$ is the true active set defined in \eqref{true.set}, and let $\cS^{\mathrm{\, c}}$ be its complement.
To assess the numerical performance of our proposed method, we report (1) the true positive rate (TPR), $\mathrm{TPR} = | \cS \cap \hat{\cS} \, | \, \big/ \, | \cS |$; (2) the false discovery rate (FDR), $\mathrm{FDR} = |\cS^{\mathrm{\, c}} \cap \hat{\cS} \, | \, \big/ \, | \hat{\cS} \, |$; (3) average $\ell_2$-error, $(1 / m) \sum_{k = 0}^m \| \hat{\bbeta}(\tau_k) - \bbeta(\tau_k) \|_2$; and (4) elapsed time for running the estimation-after-selection process, including cross-validation.

\begin{figure}[!ht]
  \centering
  \subfigure[{\scriptsize FDR under homoscedastic model}]{\includegraphics[scale=0.32]{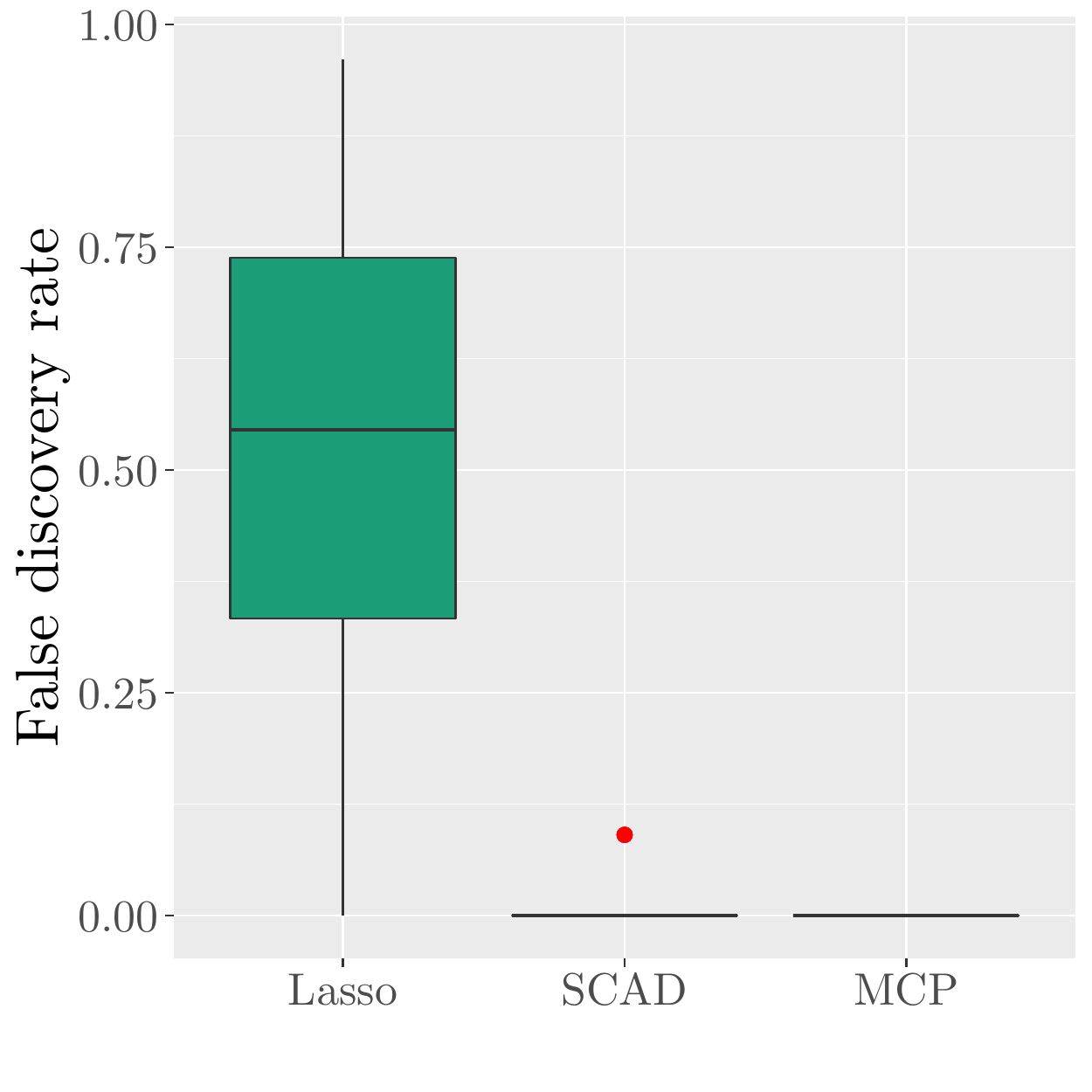}}  \qquad\quad 
  \subfigure[{\scriptsize $\ell_2$-error under homoscedastic model}]{\includegraphics[scale=0.32]{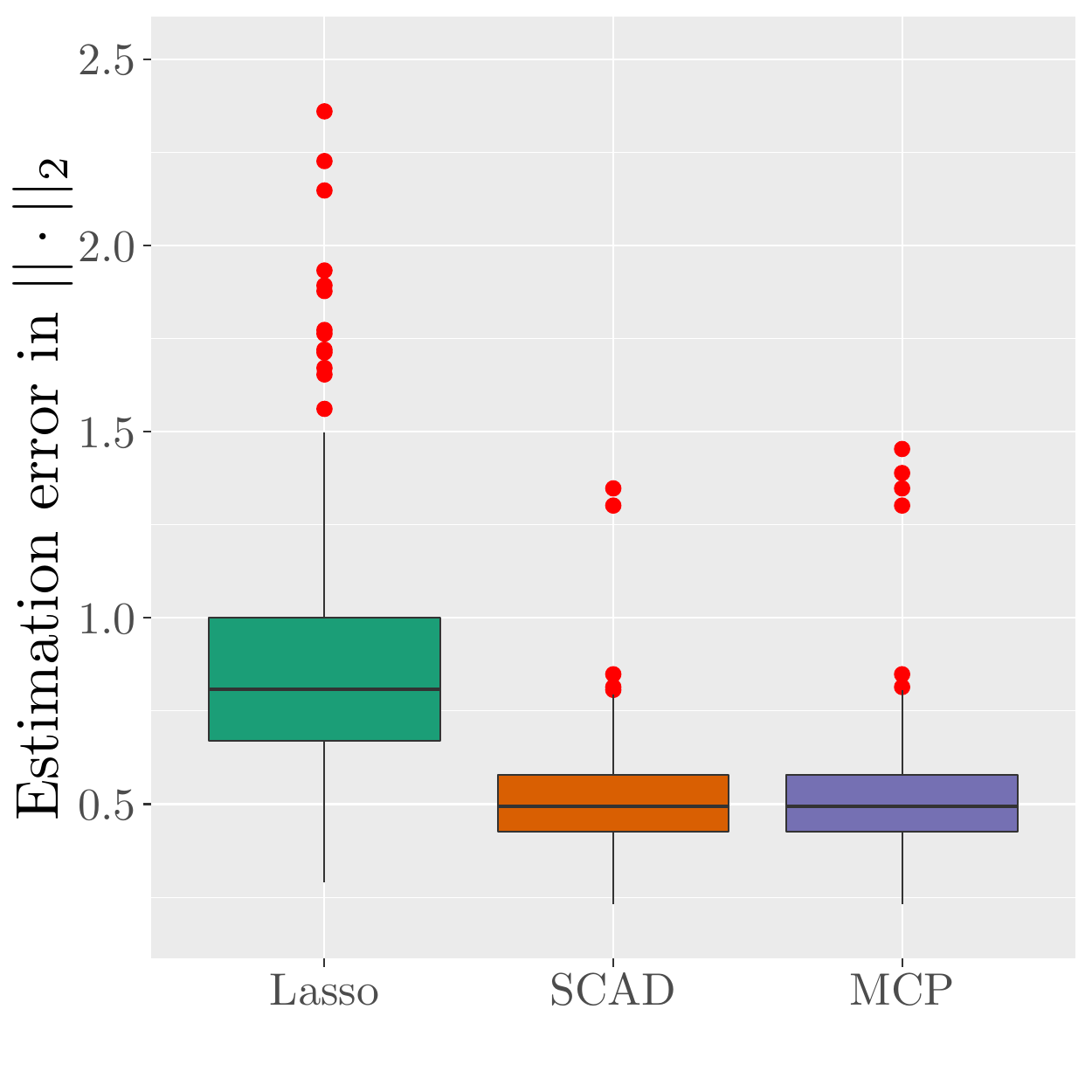}} \qquad\quad 
  \subfigure[{\scriptsize Runtime under homoscedastic model}]{\includegraphics[scale=0.32]{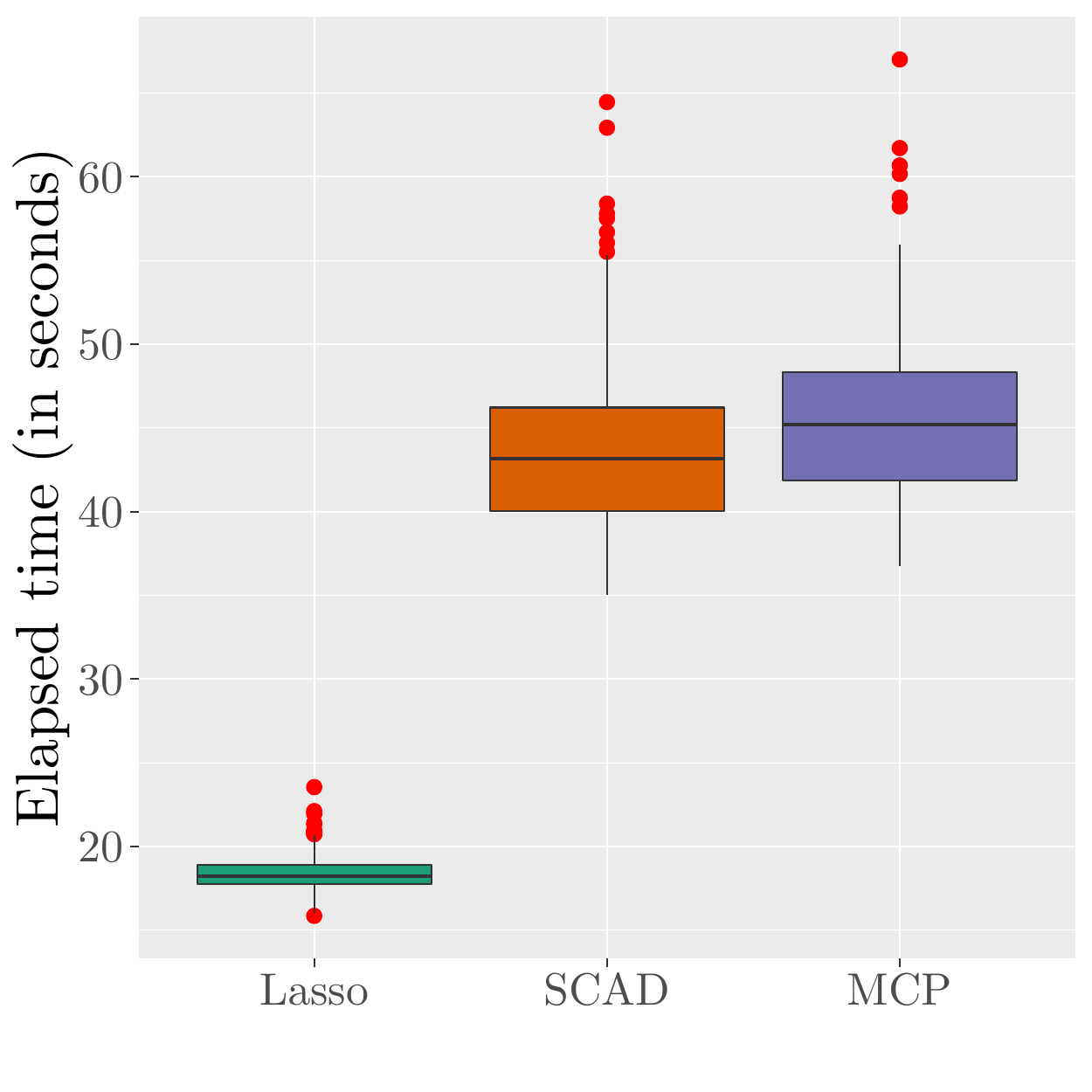}}
  \subfigure[{\scriptsize FDR under heteroscedastic model}]{\includegraphics[scale=0.32]{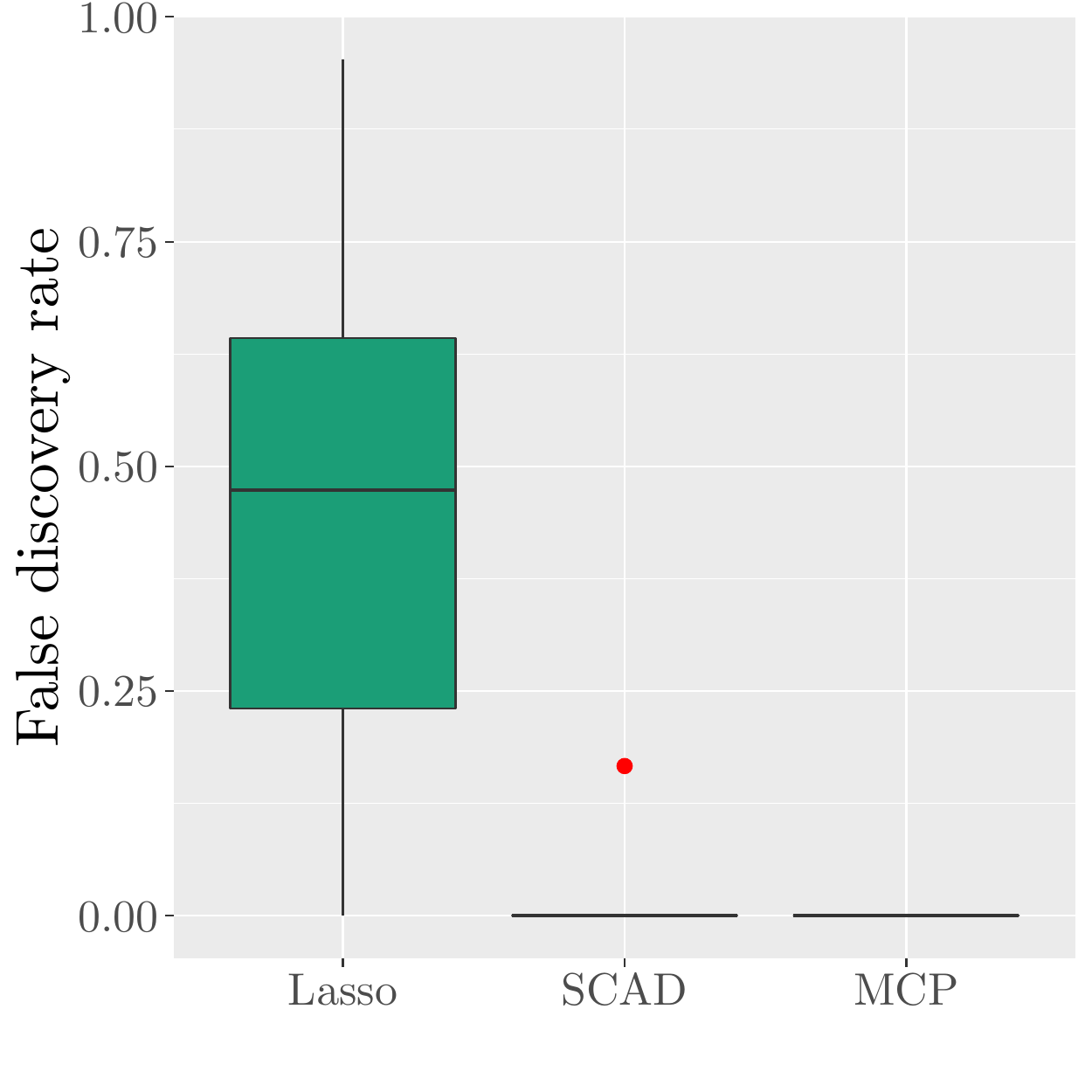}}  \qquad\quad 
  \subfigure[{\scriptsize $\ell_2$-error under heteroscedastic model}]{\includegraphics[scale=0.32]{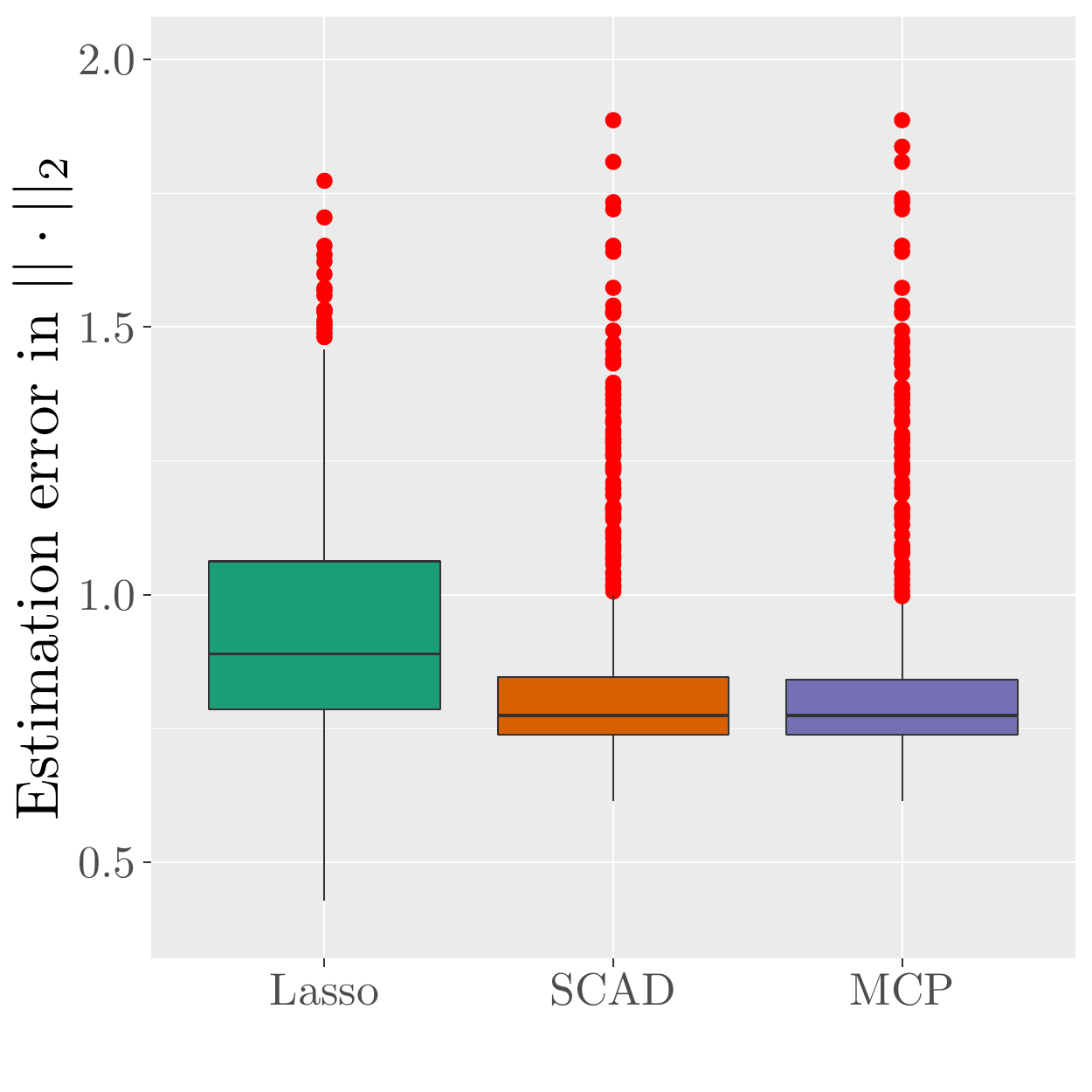}} \qquad\quad 
  \subfigure[{\scriptsize Runtime under heteroscedastic model}]{\includegraphics[scale=0.32]{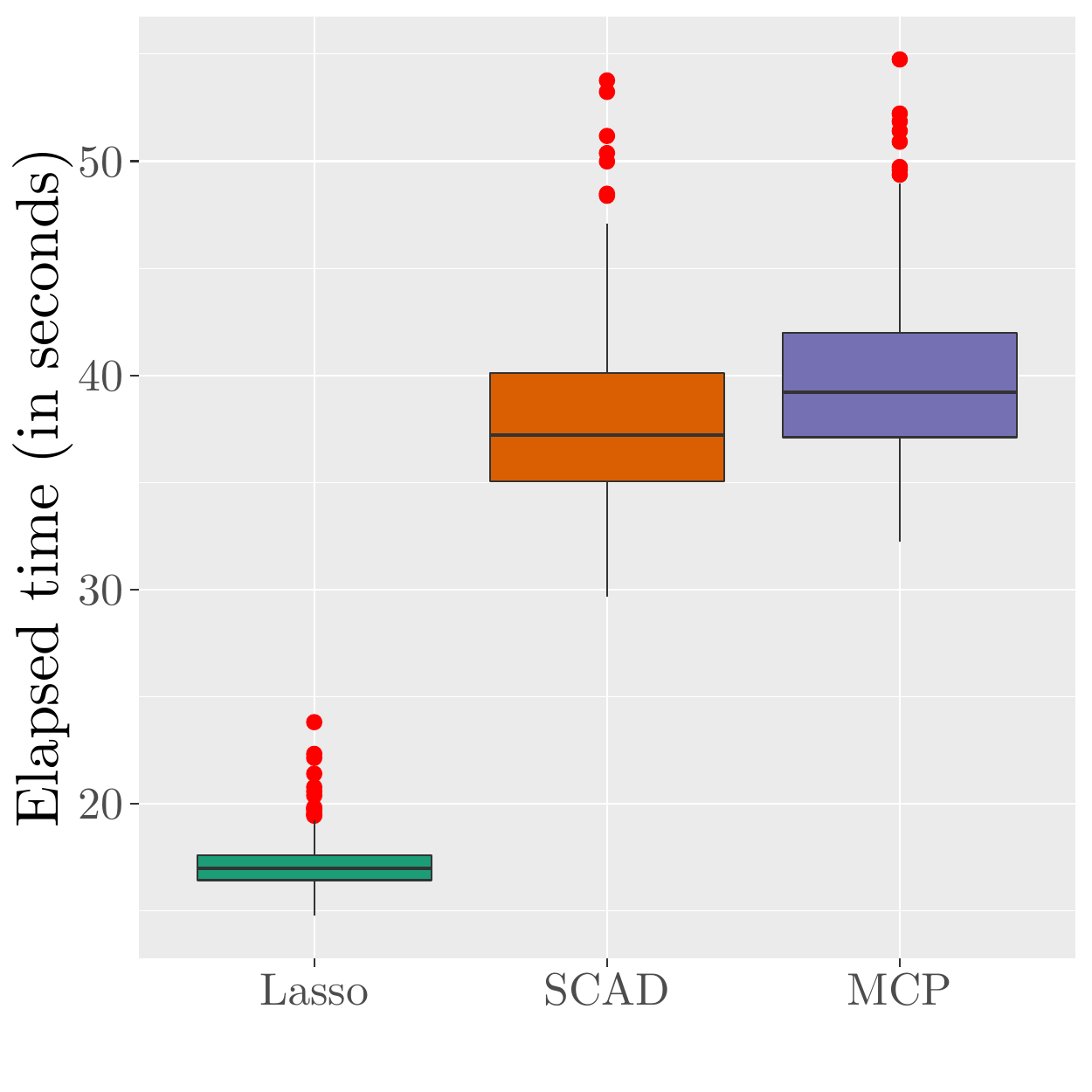}}
\caption{Box plots of the false discovery rate, $\ell_2$-error, and runtime for the $\ell_1$, SCAD, and MCP regularized smoothed CQR. The true positive rates (TPR) are not visually informative, and thus are reported as follows. For the homoscedastic model,  the average TPR are 1.00 for Lasso, 0.9996 for SCAD, and 0.9992 for MCP; for the heteroscedastic model, the average TPR are 0.9872 for Lasso, 0.919 for SCAD, and 0.917 for MCP. The censoring rates vary between $25\%$ and $30\%$.}
  \label{fig:simu.nonconvex}
\end{figure}

Results for the proposed method using different penalty functions, averaged over 500 replications when $(n, p) = (400, 1000)$, are reported in Figure~\ref{fig:simu.nonconvex}.
As expected,  $\ell_1$-penalized method tends to select larger models with many spurious variables, and thus has higher false discovery rates than SCAD and MCP. Under the heterogeneous model, both SCAD and MCP sometimes miss the first true signal and have lower TPR than Lasso. This is due to the fact that the first signal corresponds to the evolving quantile effect $Q_{t_2}(\tau)$ that vanishes as $\tau$ approaches $0.5$, and therefore is more likely to be missed by folded-concave regularization.

\begin{figure}[!htp]
  \centering
  \subfigure[{\scriptsize Runtime under homoscedastic model}]{\includegraphics[scale=0.32]{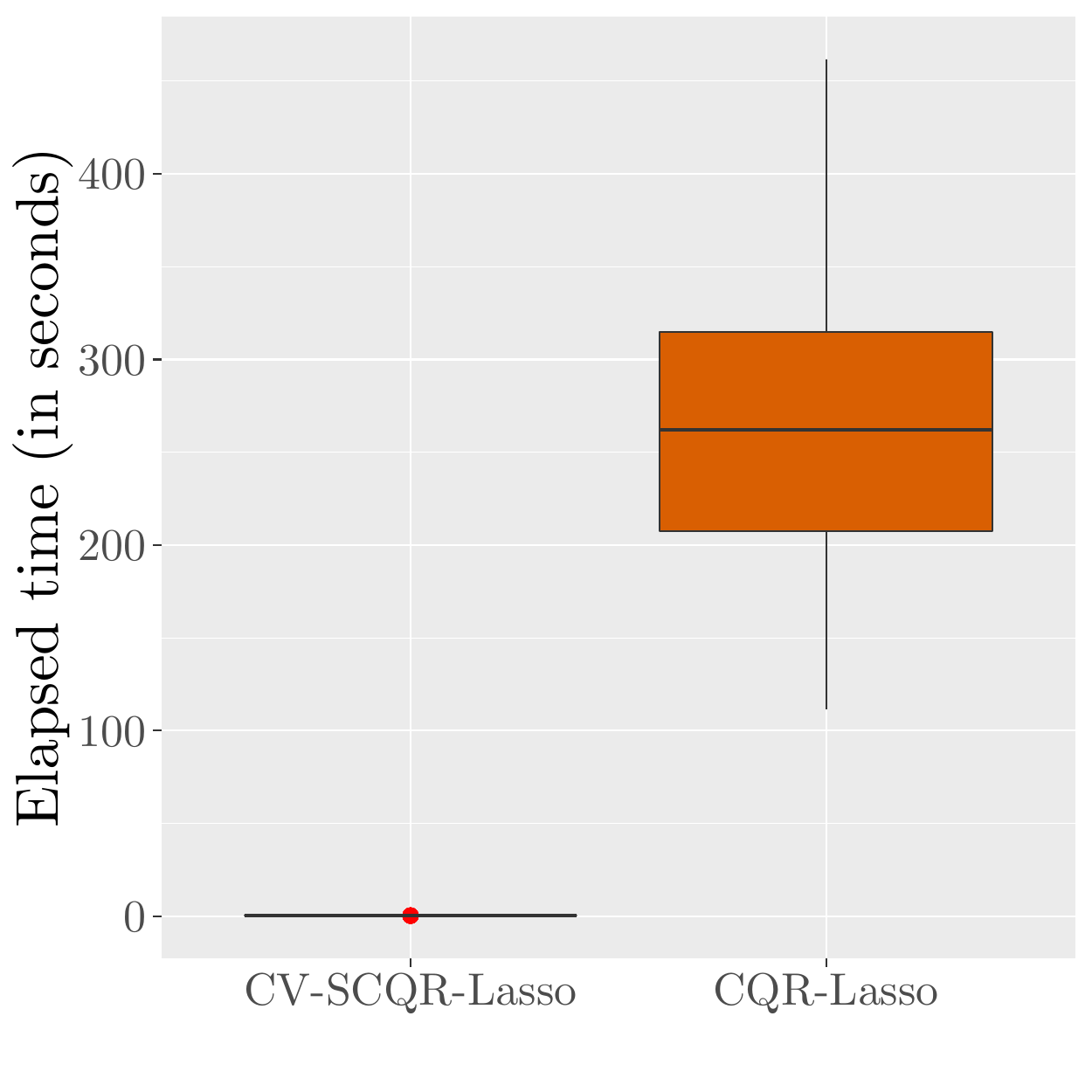}} \qquad\quad 
  \subfigure[{\scriptsize Runtime under heteroscedastic model}]{\includegraphics[scale=0.32]{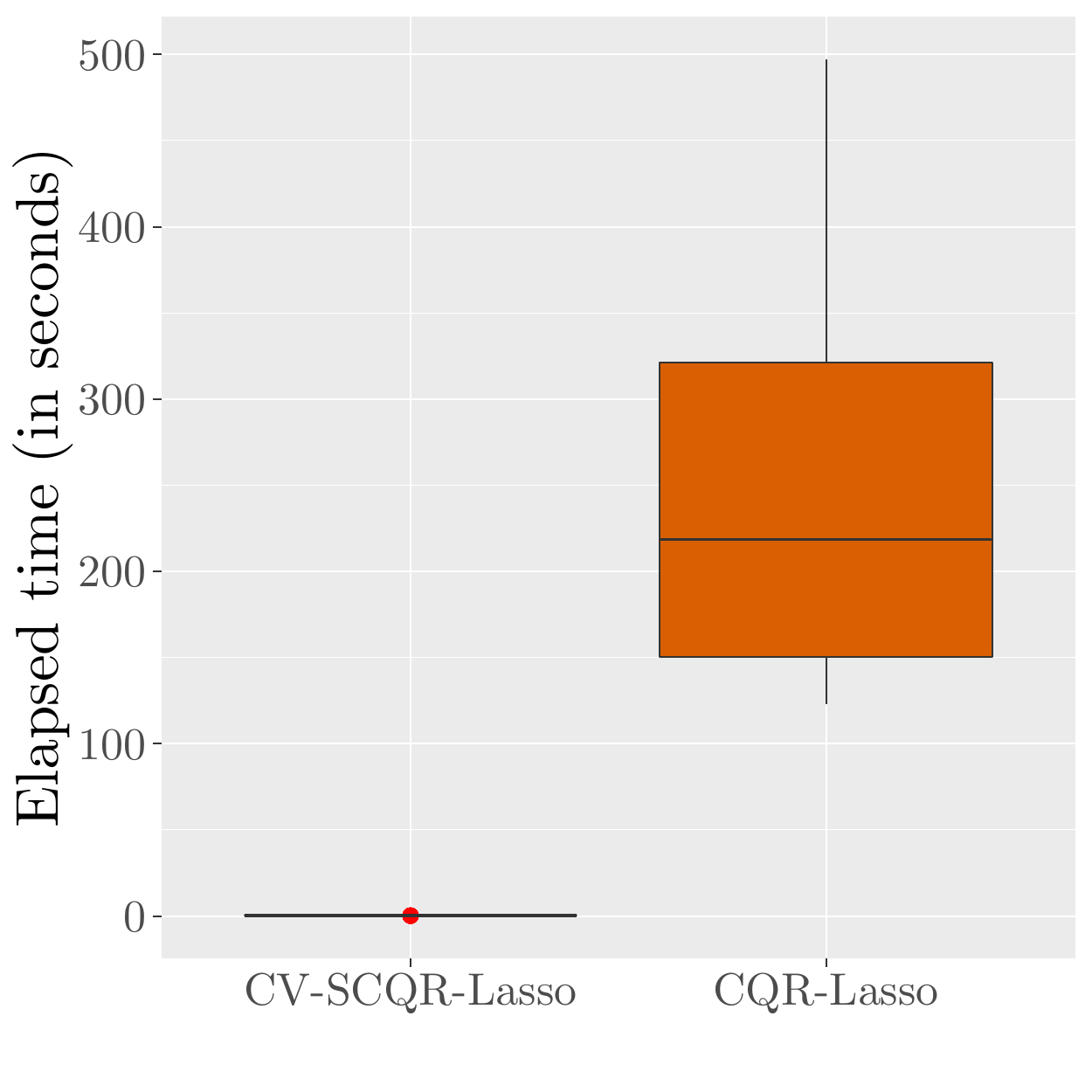}}
\caption{Box plots of runtime for $\ell_1$-penalized CQR (CQR-Lasso) and cross-validated $\ell_1$-penalized smoothed CQR (CV-SCQR-Lasso). The censoring rates vary from $25\%$ to $30\%$. 
CQR-Lasso is implemented using the \texttt{LASSO.fit} function in the package \texttt{rqPen}. The results on FDR and TPR are shown in Figure~F.4 in the supplementary material.}
  \label{fig:simu.lasso}
\end{figure}

To better demonstrate the computational efficiency of the proposed SEE method on large-scale data, we consider the $\ell_1$-penalized CQR (CQR-Lasso) method \citep{ZPH2018} as a benchmark.
As discussed in \cite{ZPH2018}, CQR-Lasso can be reformulated as a sequence of $\ell_1$-penalized median regressions with two pseudo observations, to which existing packages for penalized QR can be applied.
Moreover, \cite{ZPH2018} used cross-validation to choose $\lambda_0$ (the initial penalty level) and the increment $c>0$ by a two-dimensional grid search. In principle, we can apply this tuning scheme to both CQR-Lasso and its smoothed counterpart to achieve better variable selection performance.
From a computational point of view, we apply a simpler tuning method by only choosing $\lambda_0$ via cross-validation and focus on speed comparisons. To be specific, we first compute the cross-validated $\ell_1$-penalized smoothed CQR (SCQR-Lasso) and record its runtime, and then compute the CQR-Lasso  estimator using the same selected $\lambda$-sequence and record the runtime. For SCQR-Lasso, we apply the LAMM algorithm, described in Appendix~A.2 of the supplementary material,  to compute each $\hat \bbeta(\tau_k)$ defined in \eqref{def:l1cqr}; for CQR-Lasso, we use the \texttt{LASSO.fit} function in \texttt{rqPen} to fit the penalized median regression at each quantile level. The box plots of running time (in second) over 500 replications  are displayed in Figure~\ref{fig:simu.lasso}. On average, our implementation of the cross-validated SCQR-Lasso is more than 10 times faster than the CQR-Lasso implementation without cross-validation (18 seconds versus 250 seconds). 
The box plots of false discovery rates are shown in Figure~F.4 in the supplementary material. The code for the proposed method and our implementation of \cite{ZPH2018}'s method  is available at \href{https://github.com/XiaoouPan/scqr}{https://github.com/XiaoouPan/scqr}.

\section{Data Applications}
\label{sec:real.data}

As stated in the Introduction,  the data applications are conducted on a worker node with 2.5 GHz 32-core processor  and 512 GB of RAM in a high-performance computing cluster.

\subsection{Primary biliary cirrhosis data}

We apply the proposed method to the Mayo primary biliary cirrhosis dataset \citep{Detal1989}, a double-blinded randomized trial conducted by Mayo Clinic between 1974 and 1984.
Primary biliary cirrhosis is a rare but fatal chronic liver disease.
Our response of interest is the survival time on logarithmic scale, and an observation is censored if the patient stays alive by the end of the research. 
Five variables are included into our modeling: age in days, the presence of edema, serum bilirubin in mg/dl, albumin in gm/dl, and prothrombin time in seconds, with logarithmic transformations applied to the last three variables.   These features are  statistically significant in a multivariate Cox proportional hazards model \citep{Detal1989}. After removing data with missing covariates, the dataset contains 416 patients and a censoring rate of $61.5\%$.

\begin{figure}[!ht]
  \centering
  \subfigure[Variable: age]{\includegraphics[scale=0.32]{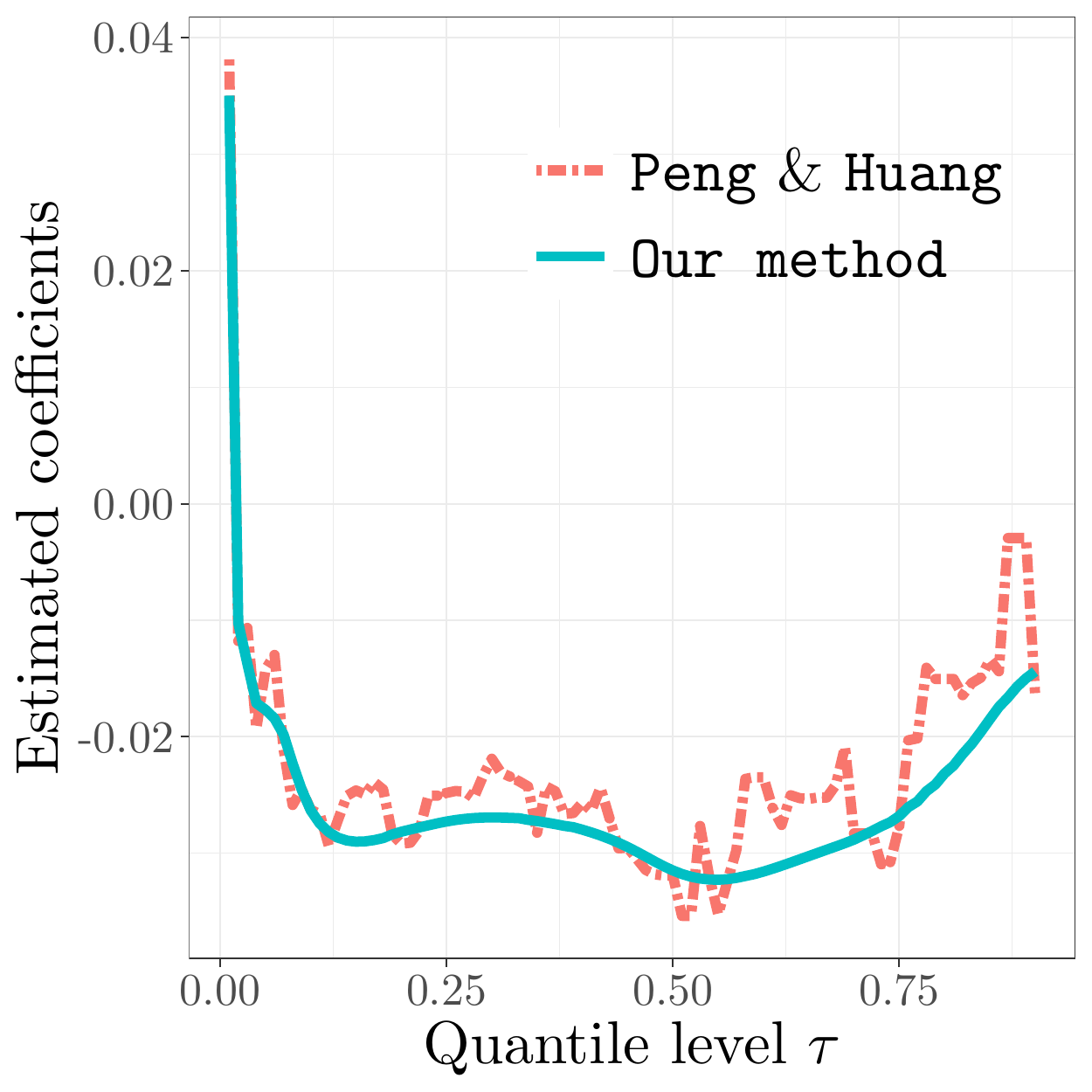}}  \qquad\quad 
  \subfigure[Variable: edema]{\includegraphics[scale=0.32]{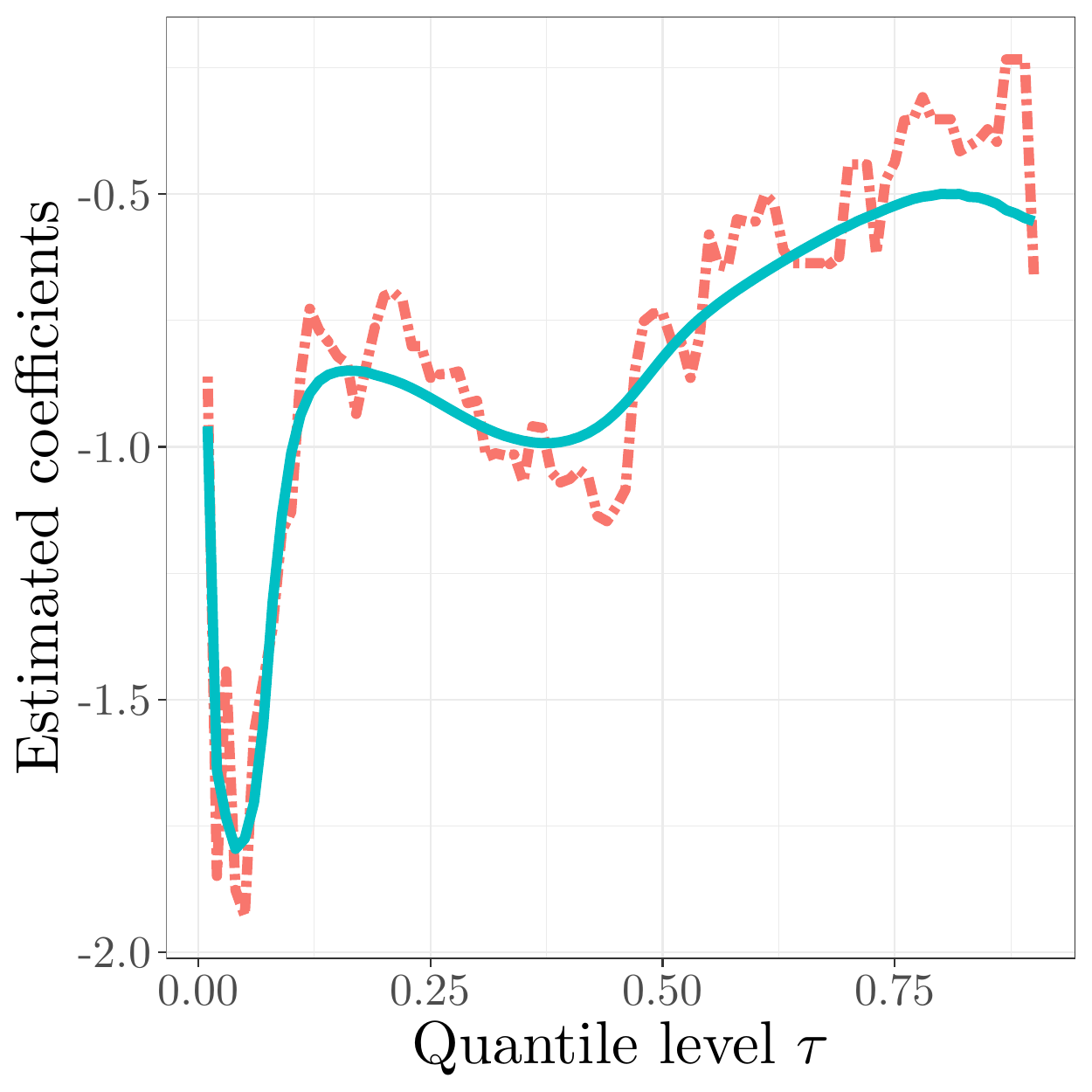}} \qquad\quad 
  \subfigure[Variable: $\log(\mathrm{bili})$]{\includegraphics[scale=0.32]{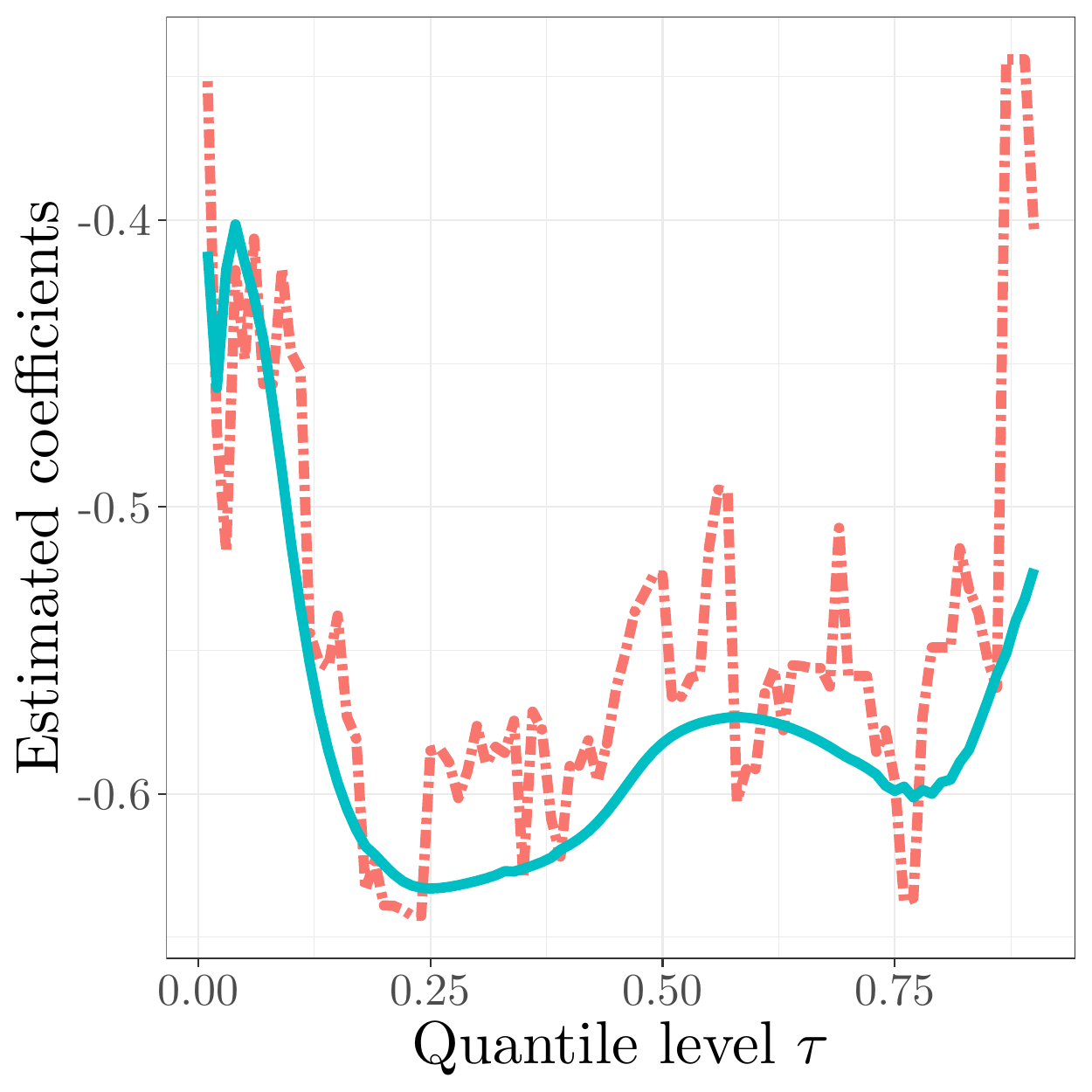}} \\
  \subfigure[Variable: $\log(\mathrm{albumin})$]{\includegraphics[scale=0.32]{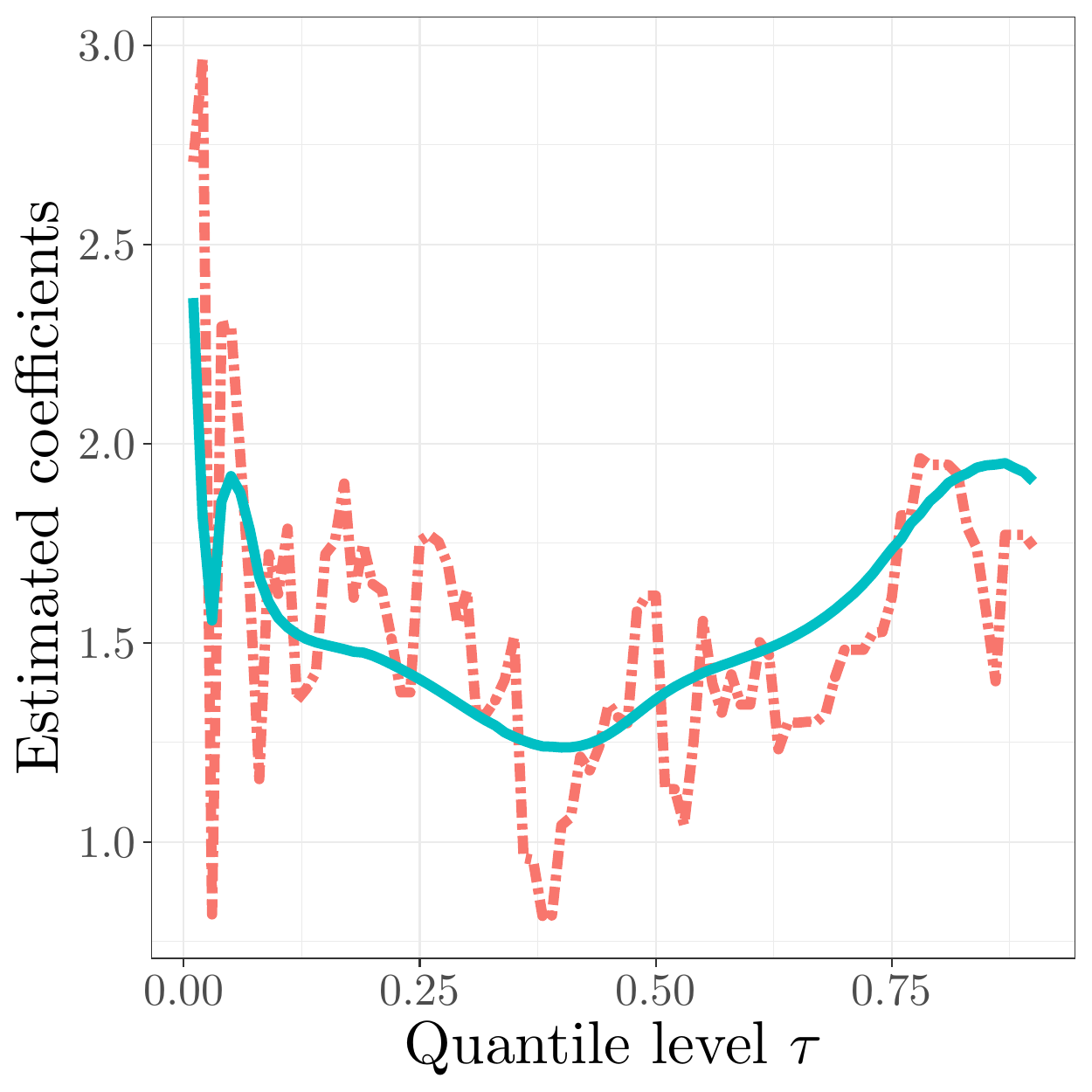}} \qquad\quad 
  \subfigure[Variable: $\log(\mathrm{protime})$]{\includegraphics[scale=0.32]{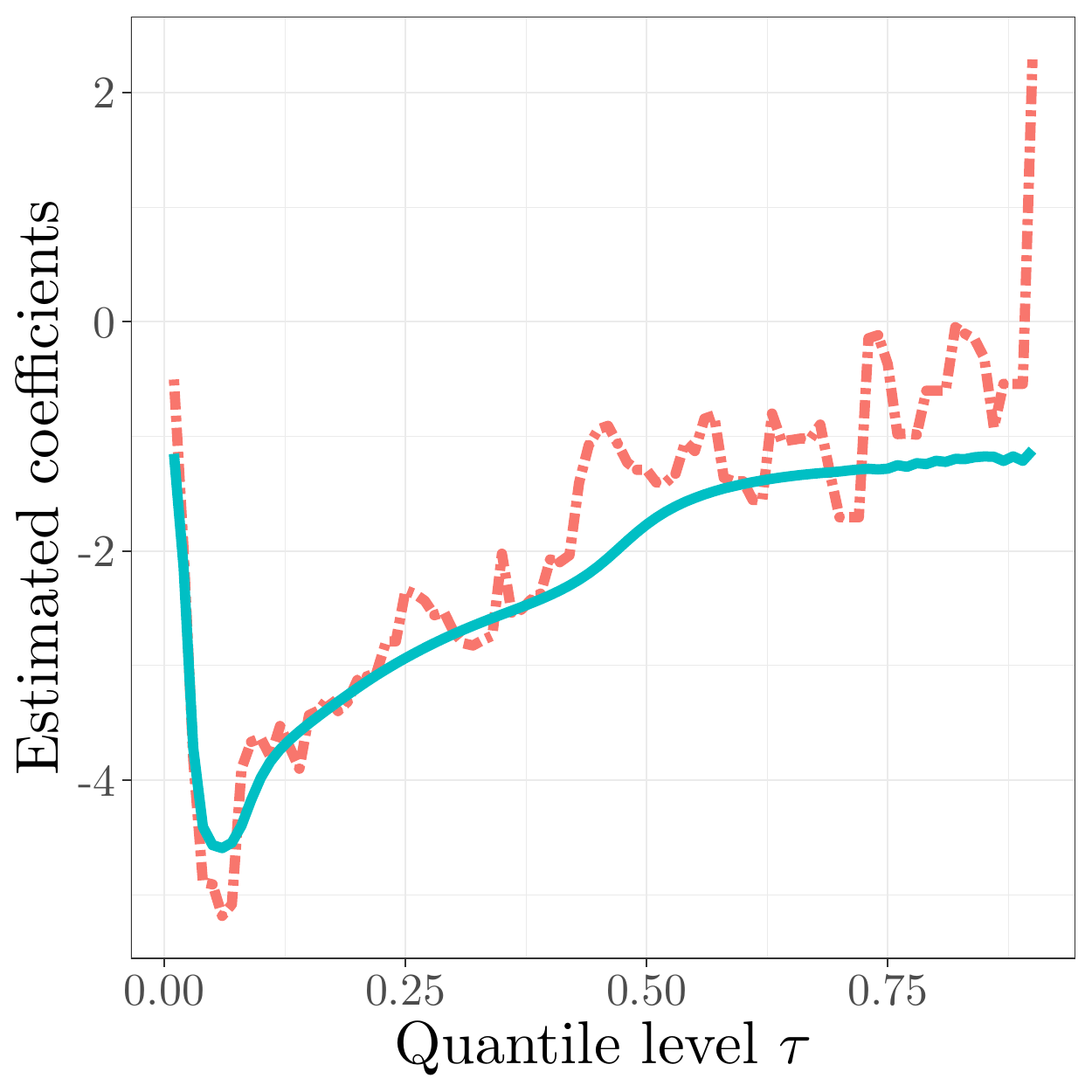}} 
\caption{Estimated regression coefficients over $\tau \in [0.01, 0.9]$ for five variables in the Mayo primary biliary cirrhosis data. Specifically, ``age'' stands for age in days, ``edema'' indicates the presence of edema, ``bili'' represents serum bilirubin in mg/dl, ``albumin'' means albumin in gm/dl, and ``protime'' refers to prothrombin time in seconds. Peng \& Huang's estimators are obtained via \texttt{crq} function in the \texttt{quantreg} package with \texttt{method = "PengHuang"}.}
  \label{fig:real.low}
\end{figure}

We apply both the classical and the proposed smoothed CQR methods to this dataset. 
The former is implemented by \texttt{crq(..., method = "PengHuang")} in the \texttt{quantreg} package over the quantile grid $\{0.01, 0.02, \dots, 0.90\}$.
The bandwidth parameter of our method is set to be $h = \max\{0.05 , \,  \{(p + \log n) / n\}^{2/5}\}$.
The estimated regression coefficients are plotted in Figure~\ref{fig:real.low} as functions of quantile levels. It is worth noting that our method leads to a fairly smooth estimated coefficient process, while there is much higher variability in the usual CQR estimator \citep{PH2008}.  Arguably this could be an advantage of the smoothed method because it produces more interpretable results.

Among the five covariates, age exhibits modest effects along the process, while albumin and prothrombin time possess varying effects with opposite signs, especially for short survivors.
Our findings echo the conclusions made in \cite{H2010}, and offer an alternative perspective to this dataset apart from \cite{TZW2005}, in which the regression coefficients are assumed to be different across the quantile levels.

\subsection{Microarray data for lung adenocarcinoma}

We now apply the proposed regularized smoothed CQR method to a gene expression-based data from a large retrospective study  for survival prediction in lung cancer \citep{Setal2008}.
The dataset provides gene expression profiling using microarray technologies, and has been briefly introduced in Section~\ref{sec:intro}. After removing observations with missing values, we have 22,283 genes from 442 lung adenocarcinomas samples, with a censoring proportion of $46.6\%$.

\begin{figure}[!ht]
  \centering
  \subfigure[Selected genes from the full data with $p = 22,283$]{\includegraphics[scale=0.32]{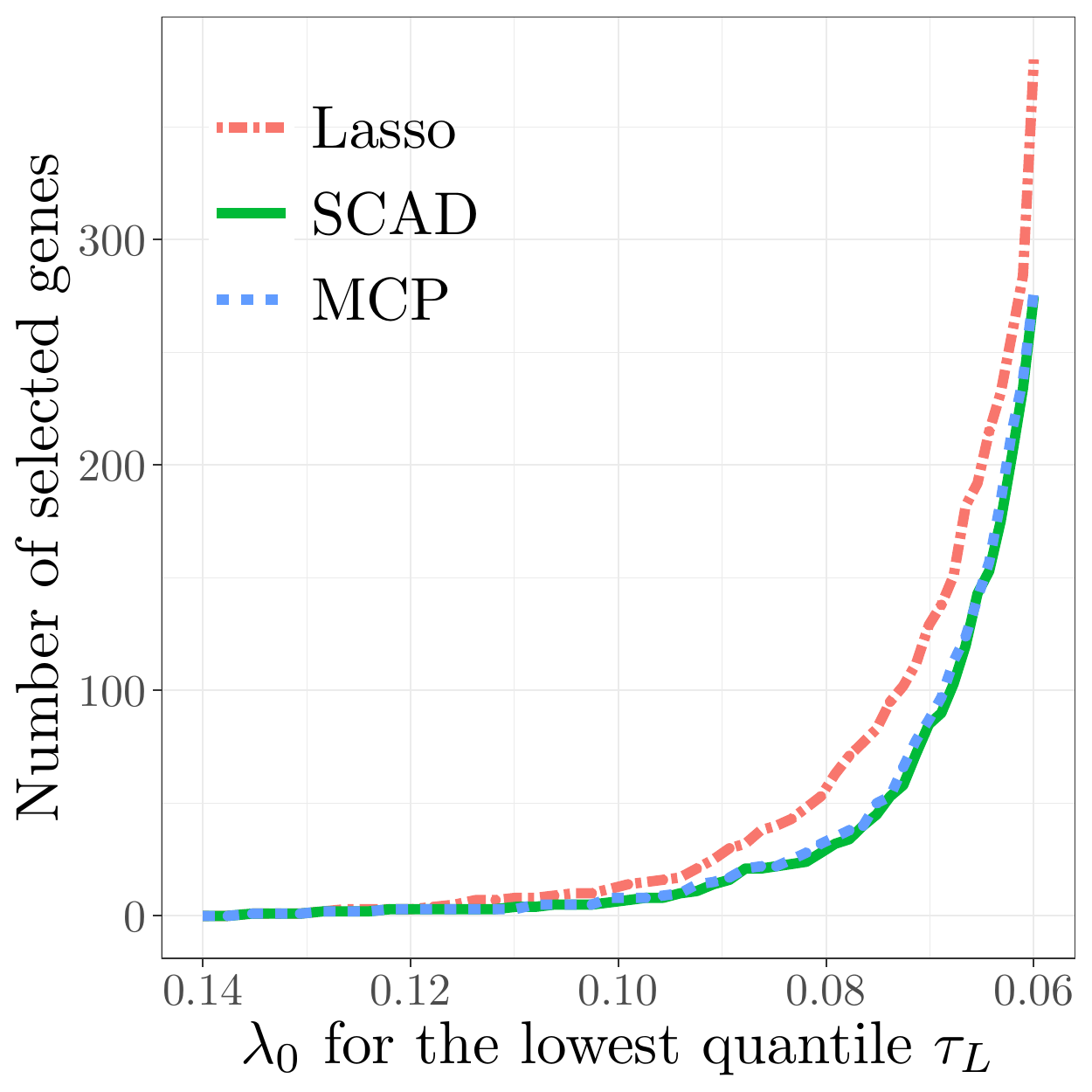}}  \qquad\quad 
  \subfigure[Selected genes from the preprocessed data with $p = 3,000$]{\includegraphics[scale=0.32]{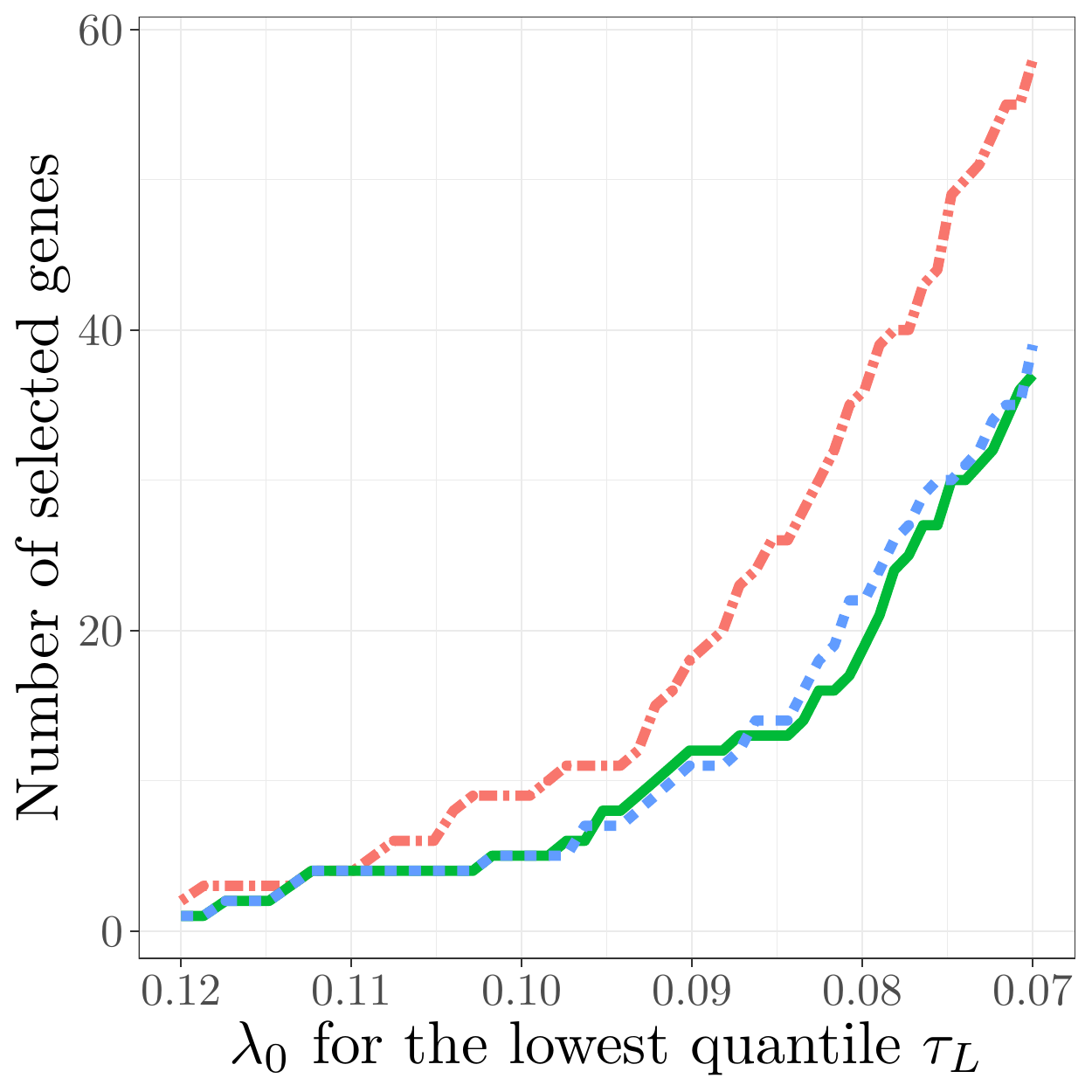}} \\
\caption{Number of selected genes from regularized smoothed CQR with Lasso, SCAD and MCP penalties, when $\lambda_0$ gradually decreases over a reasonable range. 
The left and right panels contain results for the entire data with $p = 22,283$ and the preprocessed data with $p = 3,000$, respectively.}
  \label{fig:real.gene}
\end{figure}

To demonstrate the scalability of our method, we first run regularized CQR on the whole dataset without any processing steps. Then, to roughly denoise the large dataset and to better interpret the results, we follow the preprocessing procedure carried out in \cite{ZPH2018} by selecting 3,000 genes with the largest variances, and further investigate the impact of these genes on lung cancer survival time.
For both analysis strategies, the quantile grid is set to be $\{0.1, 0.11, \dots, 0.7\}$, the bandwidth is set as $h=\{0.05 \vee 0.5\{ \log(p)/n \}^{1/4}\}$, and the tuning parameter is gradually dilating with $\lambda_k = \{ 1 + \log (\frac{1-\tau_L}{1-\tau_k} ) \} \lambda_0$ for $k=1,\ldots, m$, where $\lambda_0$ ranges over $50$ reasonable candidates.  Figure~\ref{fig:real.gene} contains the number of detected genes across various values of $\lambda_0$. Moreover, we report the first ten identified genes (denoted by their Affymetrix probe IDs)  in Table~\ref{tb.real.array}.

\begin{table}[!ht] 
\begin{center}
\setlength{\tabcolsep}{4pt}
\begin{tabular}{ c | c c c | c c c}
\hline
  & \multicolumn{3}{c|}{Whole data ($p = 22,283$)} & \multicolumn{3}{c}{Preprocessed data ($p = 3,000$)} \\
  \hline
 & Lasso & SCAD & MCP & Lasso & SCAD & MCP \\
\hline
   & 205394\_at & 205394\_at & 205394\_at & 213911\_s\_at & 213911\_s\_at & 213911\_s\_at \\ 
   & 220658\_s\_at & 220658\_s\_at & 220658\_s\_at & 217938\_s\_at & 217938\_s\_at & 217938\_s\_at \\ 
   & 221249\_s\_at & 221249\_s\_at & 221249\_s\_at & 201890\_at & 201890\_at & 201890\_at \\ 
   & 209825\_s\_at & 201250\_s\_at & 201250\_s\_at & 201250\_s\_at & 201250\_s\_at & 201250\_s\_at \\ 
 Identified & 217938\_s\_at & 40093\_at & 40093\_at & 200750\_s\_at & 200750\_s\_at & 200750\_s\_at \\ 
 genes  & 201250\_s\_at & 204728\_s\_at & 200750\_s\_at & 212951\_at & 201761\_at & 201761\_at \\ 
  & 40093\_at & 200750\_s\_at & 204728\_s\_at & 202503\_s\_at & 202503\_s\_at & 202503\_s\_at \\ 
  & 203967\_at & 218193\_s\_at & 209825\_s\_at & 209773\_s\_at & 212951\_at & 200786\_at \\ 
  & 210052\_s\_at & 203967\_at & 203967\_at & 201761\_at & 200786\_at & 209773\_s\_at \\ 
  & 218193\_s\_at & 219787\_s\_at & 219787\_s\_at & 204170\_s\_at & 209773\_s\_at & 212951\_at \\ 
\hline
  Time (in minutes) &  2.00 & 4.10 & 4.06 & 0.25 & 0.35 & 0.36 \\
\hline
\end{tabular}
\caption{The leading 10 identified genes (presented by their Affymetrix probe IDs) using regularized smoothed CQR with Lasso, SCAD and MCP penalties, as $\lambda_0$ gradually diminishes. The methods are applied to the whole data with $p = 22,283$, and the preprocessed data with $p = 3,000$. The last row indicates the average running time for each $\lambda_0$, and is recorded in minutes. The running time of $\ell_1$-penalized CQR for a single $\lambda_0$ is more than 7 days on the whole data, and 10.82 hours on the preprocessed data.}
\label{tb.real.array}
\end{center}
\end{table}

Our proposed method is computationally scalable and takes only 2 to 4 minutes for fitting the large microarray dataset with $p=22,283$.  As a reference, it takes 10.82 hours to run the $\ell_1$-penalized CQR even on the preprocessed data with $p=3,000$.  In addition, with the same data, the genes detected by Lasso and nonconvex penalties substantially overlap, and some are commonly identified regardless of the preprocessing step, e.g., ``201250\_s\_at'' and ``200750\_s\_at''. These genes may be potentially revealing and enlightening for survival prediction in lung cancer, and the intrinsic biological explanation can be a gripping topic for genetics research.

\begin{acks}[Acknowledgments]
The authors acknowledge two anonymous referees and an Associate Editor for their constructive comments that improved the quality and presentation of this paper.
\end{acks}

\begin{funding}
X. He was supported by NSF Grants DMS-1914496 and DMS-1951980.
K. M. Tan was supported by NSF Grants DMS-1949730 and DMS-2113356, and NIH Grant RF1-MH122833. W.-X. Zhou acknowledges the support of the NSF Grant DMS-2113409.
\end{funding}

\begin{supplement}
\stitle{Supplementary Material for ``Scalable Estimation and
Inference for Censored Quantile Regression Process''}
\sdescription{This supplementary material contains the proofs of all theoretical results in Sections~\ref{sec:theory} and \ref{sec:highdim}, along with the optimization algorithms and additional simulation studies.}
\end{supplement}




\end{document}